\documentclass[10pt]{amsart}

\usepackage{amsmath,amssymb,amsfonts,bbold,mathrsfs}
\usepackage{pinlabel}
\usepackage{pictex}

\DeclareMathAlphabet\oldmathcal{OMS}        {cmsy}{b}{n}
\SetMathAlphabet    \oldmathcal{normal}{OMS}{cmsy}{m}{n}
\DeclareMathAlphabet\oldmathbcal{OMS}       {cmsy}{b}{n}

\usepackage{eucal}
\newtheorem{thm}{Theorem}[section]
\newtheorem{prop}[thm]{Proposition}
\newtheorem{cor}[thm]{Corollary}
\newtheorem{lem}[thm]{Lemma}
\newtheorem{defn}[thm]{Definition}

\newtheorem{thmint}{Theorem}
\newtheorem{propint}[thmint]{Proposition}

\newenvironment{xpl}{\refstepcounter{thm} \medskip \noindent {\bf  Example \arabic{section}.\arabic{thm}}}
{\hfill$\diamondsuit$\mbox{}\bigskip}

\newenvironment{rmk}{\refstepcounter{thm} \medskip \noindent {\bf  Remark \arabic{section}.\arabic{thm}.}}{\hfill\mbox{}\bigskip}

\newcounter{num}

\newenvironment{thmlist}{\begin{list}{(\roman{num})}{\usecounter{num}\setlength{\leftmargin}{25pt}
\setlength{\itemindent}{0pt}\setlength{\labelwidth}{20pt}\setlength{\labelsep}{5pt}\setlength{\itemsep}{0in}}}{\end{list}}

\newcommand{\N}{\mathbb{N}}
\newcommand{\Z}{\mathbb{Z}}
\newcommand{\Q}{\mathbb{Q}}
\newcommand{\R}{\mathbb{R}}
\newcommand{\C}{\mathbb{C}}
\newcommand{\Ha}{\mathbb{H}}
\newcommand{\re}{\operatorname{Re}}

\newcommand{\rps}{\mathbb{R}P}
\newcommand{\cps}{\mathbb{C}P}
\newcommand{\qps}{\mathbb{H}P}
\newcommand{\ol}[1]{\bar{#1}}

\newcommand{\ASD}{\operatorname{ASD}}
\newcommand{\Aut}{\operatorname{Aut}}
\newcommand{\aut}{\operatorname{\mathfrak{aut}}}

\newcommand{\Divorb}{\operatorname{Div^{\text{orb}}}}

\newcommand{\End}{\operatorname{End}}

\newcommand{\Hol}{\operatorname{Hol}}
\newcommand{\Hom}{\operatorname{Hom}}

\newcommand{\ind}{\operatorname{Ind}}
\newcommand{\inj}{\operatorname{inj}}

\newcommand{\Isom}{\operatorname{Isom}}

\newcommand{\im}{\operatorname{Im}}

\newcommand{\Picorb}{\operatorname{Pic^{\text{orb}}}}
\newcommand{\rank}{\operatorname{rank}}
\newcommand{\Ric}{\operatorname{Ric}}
\newcommand{\Ricci}{\operatorname{Ricci}}
\newcommand{\SF}{\operatorname{SF}}
\newcommand{\Sing}{\operatorname{Sing}}

\newcommand{\Vol}{\operatorname{Vol}}
\newcommand{\weil}{\operatorname{Weil}}

\newcommand{\SL}{\operatorname{SL}}

\newcommand{\PSL}{\operatorname{PSL}}
\newcommand{\GL}{\operatorname{GL}}
\newcommand{\U}{\operatorname{U}}

\newcommand{\SO}{\operatorname{SO}}
\newcommand{\Sp}{\operatorname{Sp}}

\newcommand{\sm}[1]{\scriptscriptstyle{#1}}
\newcommand{\Sm}[1]{\scriptstyle{#1}}
\newcommand{\contr}{\,\lrcorner\,}

\title{Sasaki-Einstein 5-manifolds associated to toric 3-Sasaki manifolds}
\author{Craig van Coevering}
\address{Max-Planck-Institut fr Mathematik, Vivatsgasse 7, 53111 Bonn Germany}
\email{craigvan@mpim-bonn.mpg.de}
\date{June 15, 2012}

\begin{document}

\begin{abstract}
We give a correspondence between toric {3-Sasaki} 7-manifolds $\mathcal{S}$ and certain toric Sasaki-Einstein 5-manifolds $M$.
These 5-manifolds are all diffeomorphic to $\# k(S^2\times S^3),$ where $k=2b_2(\mathcal{S})+1$, and are given by a
pencil of Sasaki embeddings, where $M\subset\mathcal{S}$ is given concretely by the zero set of a component of the
{3-Sasaki} moment map.  It follows that there are infinitely many examples of these toric Sasaki-Einstein manifolds $M$ for each odd
$b_2(M)>1$.  This is proved by determining the invariant divisors of the twistor space $\mathcal{Z}$ of $\mathcal{S}$, and
showing that the irreducible such divisors admit orbifold K\"{a}hler-Einstein metrics.

As an application of the proof we determine the local space of anti-self-dual structures on a toric anti-self-dual Einstein orbifold.
\end{abstract}

\maketitle

\section*{Introduction}

Recall that a {3-Sasaki} manifold $\mathcal{S}$ is a Seifert $S^1$-fibration over its twistor space $\mathcal{Z}$, which is
complex contact K\"{a}hler-Einstein space, and $\mathcal{Z}$ is also the, usually singular, twistor space of an
quaternion-k\"{a}hler orbifold $\mathcal{M}$.  In this article we show that when $\mathcal{S}$ is 7-dimensional and
toric, i.e. has a two torus $T^2$ preserving the three Sasakian structures, there is a toric Sasaki-Einstein 5-manifold
$M$ naturally associated to $\mathcal{S}$.  There is a pencil of embeddings $M\subset\mathcal{S}$, which are equivariant
with respect to the $T^3$ action on $M$ and respect the respective Sasaki structures.  This is proved by
determining the $T^2_\C$-invariant divisors of $\mathcal{Z}$.  This gives a pencil with finitely many reducible elements.
Away from the reducible elements we get a toric surface $X\subset\mathcal{Z}$ whose orbifold singularities are inherited from
those of $\mathcal{Z}$.  The general picture is given in (\ref{diag:intro}), where the horizontal maps are embeddings and
vertical maps are orbifold fibrations.

There is an elementary construction of infinitely many toric {3-Sasaki} manifolds $\mathcal{S}$ for $b_2(\mathcal{S})$ any
positive integer due to C. P. Boyer, K. Galicki, B. Mann, and E. G. Rees~\cite{BGMR98}.  This is done by taking a {3-Sasaki}
version of a Hamiltonian reduction of $S^{4m+3}$ by a torus $T^{m-1}$.  In the case of 7-dimensional quotients simple numerical
criterion on the weight matrix $\Omega$ of the torus ensure that
\[\mathcal{S}_\Omega = S^{4m+3}/\negthickspace/{T^{m-1}} \]
is smooth.  It follows that the embedded $M\subset\mathcal{S}_\Omega$ is also smooth.  Thus we get infinitely many smooth
examples as in (\ref{diag:intro}).  Since $M$ is toric it is known from the classification of 5-manifolds that
it is diffeomorphic to $k\#(S^2 \times S^3)$ where in this case $k=2 b_2(\mathcal{S})+1$.

\begin{equation}\label{diag:intro}
\beginpicture
\setcoordinatesystem units <1pt, 1pt> point at 0 30
\put {$M$} at -15 60
\put {$\mathcal{S}$} at 15 60
\put {$X$} at -15 30
\put {$\mathcal{Z}$} at 15 30
\put {$\mathcal{M}$} at 15 0
\arrow <2pt> [.3, 1] from -6 60 to 6 60
\arrow <2pt> [.3, 1] from -6 30 to 6 30
\arrow <2pt> [.3, 1] from -15 51 to -15 39
\arrow <2pt> [.3, 1] from 15 51 to 15 39
\arrow <2pt> [.3, 1] from 15 21 to 15 9
\endpicture
\end{equation}

The results of this article are not only intimately related to the examples of~\cite{BGMR98}, but they also
provide examples of Einstein manifolds of positive scalar curvature exhibiting similar non-finiteness properties
in dimension 5 rather than 7.  In the above article it was shown that there are compact Einstein 7-manifolds of
positive scalar curvature with arbitrarily large total Betti number.
Also it was shown that there are infinitely many compact 7-manifolds which admit an Einstein metric of positive scalar curvature
but do not admit a metric with nonnegative sectional curvature.  The Sasaki-Einstein manifolds constructed here provide
examples of both phenomena in dimension 5.  In particular, we prove the following.

\begin{thmint}\label{thmint:main-S-E}
For each odd $k\geq 3$ there is a countably infinite number of toric Sasaki-Einstein structures on
$\# k(S^2\times S^3)$.
\end{thmint}
The next result shows that the moduli space of Einstein metrics on $\# k(S^2\times S^3)$ has infinitely many path
components.
\begin{propint}
For $M=\# k(S^2\times S^3)$ with $k>1$ odd, let $g_i$ be the sequence of Einstein metrics in the
theorem normalized so that $\Vol_{\sm g_i}(M)=1$.  Then we have $\Ric_{g_i}=\lambda_i g_i$ with the
Einstein constants $\lambda_i\rightarrow 0$ as $i\rightarrow\infty$.
\end{propint}
The result of M. Gromov~\cite{Gro} that a manifold which admits a metric of nonnegative sectional
curvature satisfies a bound on the total Betti number depending only on the dimension implies the following.
\begin{thmint}\label{thmint:Gromov}
There are infinitely many compact 5-dimensional Einstein manifolds of positive scalar curvature which
do not admit metrics on nonnegative sectional curvature.
\end{thmint}

The diagram (\ref{diag:intro}) gives a correspondence in the sense that from one of the given spaces the remaining
four are uniquely determined.  Furthermore, $M$ is smooth precisely when $\mathcal{S}$ is.  This is used in
proving the above theorems, as numerical criteria is for constructing a smooth {3-Sasaki} space $\mathcal{S}$ is
known from~\cite{BGMR98}.

In terms of toric geometry, the relation between $X$ and $M$ on the one hand and the righthand side of (\ref{diag:intro}) on the other
is elementary.  The $\ASD$ Einstein space $\mathcal{M}$ is a simply connected 4-orbifold with a $T^2$ action and is thus
characterized by the stabilizer groups along an exceptional set of 2-spheres.  And $P=\mathcal{M}/{T^2}$ is a polygon with edges
which can be labeled with $v_1, v_2, \ldots,v_\ell \in\Z^2$ which characterize the stabilizers.  Note that they are not assumed
to be primitive, as the metric may have a cone angle along the corresponding $S^2$.  It follows from the existence
of the positive scalar curvature $\ASD$ Einstein metric~\cite{CalSin06} that the vectors
$v_1, v_2, \ldots,v_\ell,-v_1, -v_2, \ldots,-v_\ell \in\Z^2$
are vertices of a convex polytope in $\R^2$.  Thus they define an augmented fan $\Delta^*$ which characterizes the toric Fano
orbifold surface $X_{\Delta^*}$ obtained above.

In Section~\ref{sec:Sasaki} we provide some necessary background on Sasaki and {3-Sasaki} manifolds and related
geometries.  In Section~\ref{sec:K-E-sym} we prove the existence of an orbifold K\"{a}hler-Einstein metric on the divisor $X$.
From this we get the Sasaki-Einstein structure on $M$.  More generally a proof is given that any symmetric toric Fano orbifold
admits a K\"{a}hler-Einstein metric.  Here \emph{symmetric} means that the normalizer $\mathcal{N}(T_{\C})\subset\Aut(X)$ of
the torus $T_\C$ acts on the characters of $T_\C$ fixing only the trivial character.  It is a result of
V. Batyrev and E. Selivanova~\cite{BatSel99} that a symmetric toric Fano manifold admits a K\"{a}hler-Einstein metric.
It was then proved by X. Wang and X. Zhu~\cite{WanZhu04} that every toric Fano manifold with vanishing Futaki invariant
has a K\"{a}hler-Einstein metric.  It was then shown by A. Futaki, H. Ono and G. Wang~\cite{FutOnWan09} that every
toric Sasaki manifold with $a\omega^T \in c_1(\mathscr{F}_\xi),\ a>0,$ where $c_1(\mathscr{F}_\xi)$ is the first Chern class
of the transversely holomorphic foliation, admits a transversal K\"{a}hler-Einstein metric.  This latter result includes
the orbifolds considered here.  But the proof included here gives a lower bound on the  Tian invariant, $\alpha_G(X)\geq 1$,
where $G\subset\mathcal{N}(T_{\C})$ is a maximal compact group.  For toric manifolds it is known that $\alpha_G(X)=1$
if and only if $X$ is symmetric~\cite{BatSel99,Son05}.

In Section~\ref{sec:Sasaki-sub} we construct the correspondence and embeddings in (\ref{diag:intro}).
In particular, in Section~\ref{subsec:twist-div} the existence of the pencil of embeddings $X\subset\mathcal{Z}$ is proved.
All of the $T^2_{\C}$-invariant divisors of $\mathcal{Z}$ are determined, and in effect the entire orbit structure of
$\mathcal{Z}$ is determined.  The reducible $T^2_{\C}$-invariant divisors $X$ represent the complex contact
bundle $\mathbf{L}=\mathbf{K}_{\mathcal{Z}}^{-\frac{1}{2}}$.

In Section~\ref{subsec:Sasak-Einst} we prove the Sasaki manifold $M$ admits a Sasaki-Einstein metric, and
prove the above theorems.

In Section~\ref{sec:asd-space} as an application of the results on the twistor space $\mathcal{Z}$ we prove that
\[\dim_{\C} H^1(\mathcal{Z},\Theta_{\mathcal{Z}})= \dim_{\C} H^1(\mathcal{Z},\Theta_{\mathcal{Z}})^{T^2} =b_2(\mathcal{Z})-2
=b_2(\mathcal{M})-1,\]
which give the dimension of the local deformation space of $\ASD$ conformal structures on $\mathcal{M}$.  This dimension
$b_2(\mathcal{M})-1 =\ell -3$, where $\ell$ is number of edges of the polygon $P=\mathcal{M}/{T^2}$ labeled by the 1-dimensional
stabilizers in $T^2$.  This is the same as the dimension of the space of deformations of $(\mathcal{M},g)$ preserving
the toric structure given by the Joyce ansatz~\cite{Joy95}.  In other words, locally every $\ASD$ deformation of the comformal
metric $[g]$ is a Joyce metric.  This is in contrast to the, in many respects similar, case of toric $\ASD$ structures on
$\# m\overline{\cps}^2$ as there are many examples of deformation preserving only an $S^1\subset T^2$~\cite{LeB91}.
It is known that the virtual dimension of the moduli space of $\ASD$ conformal structures on
$\# m\overline{\cps}^2$ is $7m-15$ plus the dimension of the conformal group.
Thus in general the expected dimension of the deformation space will be much greater than the $m-1$
dimensional space of Joyce metrics.  The deformations of $\mathcal{Z}$ are also of interest for other work of the author.
It is a consequence of results in~\cite{VCo12} that the existence of the K\"{a}hler-Einstein metric is open under deformations
of $\mathcal{Z}$.

\section{Sasaki manifolds}\label{sec:Sasaki}

We review the basics of Sasaki and {3-Sasaki} manifolds in this section.  See the monograph for more details~\cite{BoyGal08}.
The survey article~\cite{BoyGal99} is a good introduction to \mbox{3-Sasakian} geometry.  These references are a good source
of background on orbifolds and orbifold bundles which will be used in this article.  In a few places we will make use of
orbifold invariants $\pi_1^{orb}(X),\ H^*_{orb}(X)$ etc. which make use of local classifying spaces $B(X)$ for orbifolds.
An introduction to these topics can be found in the above references.

\subsection{Sasaki structures}\label{subsec:Sasak-st}

\begin{defn}
A Riemannian manifold $(M,g)$ is a \emph{Sasaki manifold}, or has a compatible Sasaki structure, if the metric cone
$(C(M),\ol{g})=(\R_{>0} \times M, dr^2 +r^2 g)$ is K\"{a}hler with respect to some complex structure $I$, where $r$ is the
usual coordinate on $\R_{>0}$.
\end{defn}
Thus $M$ is odd and denoted $n=2m+1$, while $C(M)$ is a complex manifold with $\dim_{\C} C(M) =m+1$.

Although, this is the simplest definition, Sasaki manifolds were originally defined as a special type of metric contact
structure.  We will identify $M$ with the $\{1\}\times M\subset C(M)$.  Let $r\partial_r$ be the Euler vector field on $C(M)$,
then it is easy to see that $\xi =Ir\partial_r$ is tangent to $M$.  Using the warped product formulae for the cone metric
$\ol{g}$~\cite{ONeil83} it is easy check that $r\partial_r$ is real holomorphic, $\xi$ is Killing with respect to both $g$ and
$\ol{g}$, and furthermore the orbits of $\xi$ are geodesics on $(M,g)$.
Define $\eta =\frac{1}{r^2}\xi\contr\ol{g}$, then we have
\begin{equation}
\eta =-\frac{I^* dr}{r} =d^c \log r,
\end{equation}
where $d^c =\sqrt{-1}(\ol{\partial} -\partial)$.  If $\omega$ is the K\"{a}hler form of $\ol{g}$, i.e.
$\omega(X,Y) =\ol{g}(IX,Y)$, then $\mathcal{L}_{r\partial_r} \omega =2\omega$ which implies that
\begin{equation}\label{eq:Kaehler-pot1}
\omega =\frac{1}{2}d(r\partial_r \contr\omega) =\frac{1}{2}d(r^2 \eta)=\frac{1}{4}dd^c(r^2).
\end{equation}
From (\ref{eq:Kaehler-pot1}) we have
\begin{equation}\label{eq:Kaehler-pot2}
\omega=rdr\wedge\eta +\frac{1}{2}r^2 d\eta.
\end{equation}

We will use the same notation to denote $\eta$ and $\xi$ restricted to $M$.  Then (\ref{eq:Kaehler-pot2}) implies that
$\eta$ is a contact form with Reeb vector field $\xi$, since $\eta(\xi)=1$ and $\mathcal{L}_{\xi} \eta =0$.
Let $D\subset TM$ be the contact distribution which is defined by
\begin{equation}
D_x =\ker\eta_x
\end{equation}
for $x\in M$.  Furthermore, if we restrict the almost complex structure to $D$, $J:=I|_D$, then $(D,J)$ is a strictly pseudoconvex CR structure on $M$.  We have a splitting of the tangent bundle $TM$
\begin{equation}
TM =D\oplus L_{\xi},
\end{equation}
where $L_{\xi}$ is the trivial subbundle generated by $\xi$.  It will be convenient to define a tensor $\Phi\in\End(TM)$ by
$\Phi|_D =J$ and $\Phi(\xi) =0$.  Then
\begin{equation}\label{eq:comp-tens}
\Phi^2 =-\mathbb{1} +\eta\otimes\xi.
\end{equation}
Since $\xi$ is Killing, we have
\begin{equation}
d\eta (X,Y) =2 g(\Phi(X),Y), \quad\text{where }X,Y\in\Gamma(TM),
\end{equation}
and $\Phi(X) =\nabla_X \xi$, where $\nabla$ is the Levi-Civita connection of $g$.  Making use of (\ref{eq:comp-tens}) we see
that
\[ g(\Phi X,\Phi Y) =g(X,Y) -\eta(X)\eta(Y), \]
and one can express the metric by
\begin{equation}\label{eq:metric}
g(X,Y)=\frac{1}{2}(d\eta)(X,\Phi Y) +\eta(X)\eta(Y).
\end{equation}

We will denote a Sasaki structure on $M$ by $(g,\eta,\xi,\Phi)$.  Although, the reader can check that merely specifying
$(g,\xi),\ (g,\eta),$ or $(\eta, \Phi)$ is enough to determine the Sasaki structure, it will be convenient to denote the
remaining structure.

The action of $\xi$ generates a foliation $\mathscr{F}_\xi$ on $M$ called the \emph{Reeb} foliation.
Note that it has geodesic leaves and is a Riemannian foliation, that is has a $\xi$ invariant Riemannian metric on the
normal bundle $\nu(\mathscr{F}_\xi)$.  But in general the leaves are not compact.  If the leaves are compact, or equivalently
$\xi$ generates an $S^1$-action, then $(g,\eta,\xi,\Phi)$ is said to be a \emph{quasi-regular} Sasaki structure, otherwise it is
\emph{irregular}.  If this $S^1$ action is free, then $(g,\eta,\xi,\Phi)$ is said to be \emph{regular}.  In this last case
$M$ is an $S^1$-bundle over a manifold $Z$, which we will see below is K\"{a}hler.  If the structure if merely quasi-regular, then
the leaf space has the structure of a K\"{a}hler orbifold $Z$.

The vector field $\xi -\sqrt{-1}I\xi =\xi +\sqrt{-1}r\partial_r$ is holomorphic on $C(M)$.  If we denote by $\tilde{\C}^*$ the
universal cover of $\C^*$, then $\xi +\sqrt{-1}r\partial_r$ induces a holomorphic action
of $\tilde{\C}^*$ on $C(M)$.  The orbits of $\tilde{\C}^*$ intersect $M\subset C(M)$ in the orbits of the Reeb foliation
generated by $\xi$.  We denote the Reeb foliation by $\mathscr{F}_\xi$.  This gives $\mathscr{F}_\xi$ a transversely holomorphic
structure.

We define a \emph{transversely K\"{a}hler} structure on $\mathscr{F}_\xi$ with K\"{a}hler form and metric
\begin{gather} \omega^T =\frac{1}{2}d\eta \\
g^T =\frac{1}{2}d\eta(\cdot,\Phi\cdot).
\end{gather}

Though in general it is not the case, the examples in this article will be quasi-regular.  Therefore, the transversely K\"{a}hler
leaf space of $\mathscr{F}_\xi$ will be a K\"{a}hler orbifold $Z$.  Up to a homothetic transformation all such examples are
as in the following example.

\begin{xpl}
Let $\mathbf{F}$ be a negative holomorphic orbifold line bundle on a complex orbifold $Z$ and $h$ an Hermitian connection
with negative curvature.  Define $r^2 =h(w,w)$ where $w$ is the fiber coordinate.  Then $\omega=\frac{1}{4}dd^c r^2$ is the
K\"{a}hler form of a cone metric on the total space minus the zero section $\mathbf{F}^\times$.  Then
$\eta=\frac{1}{2} d^c \log r^2$ is a contact form, and since $\mathbf{F}$ is negative
\[ \omega^T =\frac{1}{2}d\eta =\frac{1}{4} dd^c \log r^2 =-\frac{1}{2}\Theta_{\mathbf{F}} >0\]
gives the transversal K\"{a}hler metric.  In this case $\omega^T$ is an orbifold K\"{a}hler metric on $Z$.
\end{xpl}

The following follows from O'Neill tensor computations for a Riemannian submersion.  See~\cite{ONeil66} and~\cite[Ch. 9]{Bess87}.
\begin{prop}\label{prop:Sasaki-Ric}
Let $(M,g,\eta,\xi,\Phi)$ be a Sasaki manifold of dimension $n=2m+1$, then
\begin{thmlist}
\item  $\Ric_g (X,\xi) =2m\eta(X),\quad\text{for }X\in\Gamma(TM)$,\label{eq:submer-Ric-Reeb}
\item  $\Ric^T (X,Y) =\Ric_g (X,Y) +2g^T(X,Y),\quad\text{for }X,Y\in\Gamma(D),$
\item  $s^T =s_g +2m.$\label{eq:submer-scal}
\end{thmlist}
\end{prop}
\begin{defn}
A \emph{Sasaki-Einstein} manifold $(M,g,\eta,\xi,\Phi)$ is a Sasaki manifold with
\[ \Ric_g =2m\, g.\]
\end{defn}
Note that by (\ref{eq:submer-Ric-Reeb}) the Einstein constant must be $2m$, and $g$ is Einstein precisely when the cone $(C(M),\ol{g})$
is Ricci-flat.  Furthermore, the transverse K\"{a}hler metric is also Einstein
\begin{equation}
\Ric^T =(2m+2)\, g^T.
\end{equation}
Conversely, if one has a Sasaki structure $(g,\eta,\xi,\Phi)$ with $\Ric^T =\tau\, g^T$ with $\tau>0$, then after a
$D$-homothetic transformation one has a Sasaki-Einstein structure $(g',\eta',\xi',\Phi)$, where
$\eta'=a\eta,\ \xi'=a^{-1}\xi$, and $g'=ag +a(a-1)\eta\otimes\eta$, with $a=\frac{\tau}{2m+2}$.

\subsection{3-Sasaki and related structures}\label{subsec:3-Sasak}

Recall that a hyper-K\"{a}hler structure on a $4m$-dimensional manifold consists of a metric $g$ which is K\"{a}hler with respect to three complex structures $J_1, J_2, J_3$ satisfying the quaternionic relations
\[J_1^2 =J_2^2 =J_3^2 =-\mathbb{1},\ J_1 J_2 =-J_2 J_1 =J_3.\]

\begin{defn}\label{defn:3-Sasaki}
A Riemannian manifold $(\mathcal{S},g)$ is 3-Sasaki if the metric cone $(C(\mathcal{S}),\ol{g})$ is hyper-K\"{a}hler.
That is, $\ol{g}$ admits compatible
almost complex structures $J_i,\ i=1,2,3 $ such that $(\ol{g},J_1,J_2,J_3)$ is a hyper-K\"{a}hler structure on
$C(\mathcal{S})$.  Equivalently, $\Hol(C(\mathcal{S}))\subseteq\Sp(m)$.
\end{defn}
A consequence of the definition is that $(\mathcal{S},g)$ is equipped with three Sasaki structures
$(g,\eta_i,\xi_i,\phi_i),\ i=1,2,3$.
The Reeb vector fields $\xi_i =J_i(r\partial_r),\ i=1,2,3$ are orthogonal and satisfy $[\xi_i,\xi_j]=-2\varepsilon^{ijk}\xi_k$, where $\varepsilon^{ijk}$ is anti-symmetric in the indicies $i,j,k \in\{1,2,3\}$ and $\epsilon^{123}=1$.
The tensors $\phi_i ,\ i=1,2,3$ satisfy the identities
\begin{align}
\phi_i(\xi_j) & =\varepsilon^{ijk}\xi_k  \label{eq:3-Sasak-ident1} \\
\phi_i \circ\phi_j & =-\delta_{ij}\mathbb{1} +\epsilon^{ijk}\phi_k +\eta_j \otimes\xi_i \label{eq:3-Sasak-ident2}
\end{align}
It is easy to see that there is an $S^2$ of Sasaki structures with Reeb vector field
$\xi_\tau =\tau_1 \xi_1 +\tau_2 \xi_2 +\tau_3 \xi_3$ with $\tau\in S^2$.

The Reeb vector fields $\{\xi_1,\xi_2,\xi_3 \}$ generate a Lie algebra $\mathfrak{sp}(1)$, so there is an effective isometric action
of either $\SO(3)$ or $\Sp(1)$ on $(\mathcal{S},g)$.  Both cases occur in the examples in this article.
This action generates a foliation $\mathscr{F}_{\xi_1,\xi_2,\xi_3}$ with generic leaves either $\SO(3)$ or $\Sp(1)$.

If we set $D_i =\ker\eta_i \subset T\mathcal{S},\ i=1,2,3$ to be the contact subbundles, then the complex structures $J_i,\ i=1,2,3$ are
recovered by
\begin{equation}
J_i (r\partial_r) =\xi_i,\quad J_i|_{D_i} =\phi_i.
\end{equation}

Because a hyper-K\"{a}her manifold is always Ricci-flat we have the following.
\begin{prop}
A {3-Sasaki} manifold $(\mathcal{S},g)$ of dimension $4m+3$ is Einstein with Einstein constant $\lambda=4m+2$.
\end{prop}

We choose a Reeb vector field $\xi_1$ fixing a Sasaki structure, then the leaf space $\mathscr{F}_{\xi_1}$ is
a K\"{a}hler orbifold $\mathcal{Z}$ with respect to the transversal complex structure $J=\Phi_1$.
But it has addition has a complex contact structure and a fibering by rational curves which we now describe.
The 1-form $\eta^c =\eta_2 -\sqrt{-1}\eta_3$ is a $(1,0)$-form with respect to $J$. But it is not invariant under
the $\U(1)$ group generated by $\exp(t\xi_1)$.  We have $\exp(t\xi_1)^* \eta^c =e^{2\sqrt{-1}t}\eta^c$.
Let $\mathbf{L}= \mathcal{S}\times_{\U(1)}\C$, with $\U(1)$ action on $\C$ by be $e^{2\sqrt{-1}t}$.
This is a holomorphic orbifold line bundle; in fact $C(\mathcal{S})$ is either $\mathbf{L}^{-1}$ or $\mathbf{L}^{\frac{1}{2}}$
minus the zero section.  It is easy to see that each of these cases occur precisely where the Reeb vector fields generate
an effective action of $\SO(3)$ and $\Sp(1)$ respectively.  Then $\eta^c$ descends to an $\mathbf{L}$ valued holomorphic
1-form $\theta\in\Gamma\bigl(\Omega^{1,0}(\mathbf{L})\bigr)$.  It follows easily from identities (\ref{eq:3-Sasak-ident2})
that $d\eta^c$ restricted to $D_1 \cap\ker\eta^c$ is a non-degenerate type $(2,0)$ form.  Thus $\theta$ is a complex contact
form on $\mathcal{Z}$, and $\theta\wedge(d\theta)^m \in\Gamma\bigl(\mathbf{K}_{\mathcal{Z}}\otimes\mathbf{L}^{m+1}\bigr)$ is
a non-vanishing section.  Thus $\mathbf{L}\cong\mathbf{K}_{\mathcal{Z}}^{-\frac{1}{m+1}}$.

Each leaf of $\mathscr{F}_{\xi_1,\xi_2,\xi_3}$ descends to a rational curve in $\mathcal{Z}$.  Each curve is a $\cps^1$
but may have orbifold singularities for non-generic leaves.  We see that restricted to a leaf
$\mathbf{L}|_{\cps^1} =\mathcal{O}(2)$.

The element $\exp(\frac{\pi}{2}\xi_2)$ acts on $\mathcal{S}$ taking $\xi_1$ to $-\xi$, thus it descends to an anti-holomorphic
involution $\varsigma:\mathcal{Z}\rightarrow\mathcal{Z}$.  This \emph{real structure} is crucial to the twistor approach.
Note that $\varsigma^*\theta =\ol{\theta}$.

This all depends on the choice $\xi_1 \in S^2$ of the Reeb vector field.  But taking a different Reeb vector field gives an
isomorphic twistor space under the transitive action of $\Sp(1)$.

Taking the quotient of $\mathcal{S}$ by $\Sp(1)$ gives the leaf space of $\mathscr{F}_{\xi_1,\xi_2,\xi_3}$ an orbifold $\mathcal{M}$.
We now consider the orbifold $\mathcal{M}$ more closely.  Let $(\mathcal{M},g)$ be any $4m$ dimensional
Riemannian orbifold.  An \emph{almost quaternionic} structure on $\mathcal{M}$ is
a rank 3 V-subbundle $\mathcal{Q}\subset End(T\mathcal{M})$ which is locally spanned by almost complex
structures $\{J_i\}_{i=1,2,3}$ satisfying the quaternionic identities $J_i^2=-\mathbb{1}$ and
$J_1 J_2 =-J_2 J_1 =J_3$.  We say that $\mathcal{Q}$ is compatible with $g$ if $J_i^*g=g$ for $i=1,2,3$.
Equivalently, each $J_i,i=1,2,3$ is skew symmetric.
\begin{defn}\label{defn:quatern-kahler}
A Riemannian orbifold $(\mathcal{M},g)$ of dimension $4m, m>1$ is \emph{quaternion-K\"{a}hler} if there is an almost
quaternionic structure $\mathcal{Q}$ compatible with $g$ which is preserved by the Levi-Civita
connection.
\end{defn}
This definition is equivalent to the holonomy of $(\mathcal{M},g)$ being contained in
\linebreak $Sp(1)Sp(m)$.
For orbifolds this is the holonomy on $\mathcal{M}\setminus S_{\mathcal{M}}$ where $S_{\mathcal{M}}$ is
the singular locus of $\mathcal{M}$.
Notice that this definition always holds on an oriented Riemannian 4-manifold ($m=1$).
This case requires a different definition.  Consider the \emph{curvature operator}
\[\mathcal{R}:\Lambda^2\rightarrow\Lambda^2 \]
of an oriented Riemannian 4-manifold.
With respect to the decomposition $\Lambda^2=\Lambda^2_+ \oplus\Lambda^2_-$, we have
\begin{equation}\label{mat:curv-dec}
\mathcal{R}=
\left\lgroup
\begin{matrix}
W_g^+ +\frac{s_g}{12} & \overset{\circ}{\Ric}_g \\
\overset{\circ}{\Ric}_g & W_g^- +\frac{s_g}{12}
\end{matrix}\right\rgroup,
\end{equation}
where $W^+_g$ and $W^-_g$ are the self-dual and anti-self-dual pieces of the Weyl curvature and
$\overset{\circ}{\Ric}_g=\Ric_g -\frac{s_g}{4}g$ is the trace-free Ricci curvature.
An oriented 4 dimensional Riemannian orbifold $(\mathcal{M},g)$ is quaternion-K\"{a}hler if it is Einstein
and anti-self-dual, meaning that $\overset{\circ}{\Ric}_g=0$ and $W^+_g =0$.

One can prove that the $\{\Phi_i\}_{i=1,2,3}$ restricted to $D_1 \cap D_2 \cap D_3$, the horizontal space
to $\mathscr{F}_{\xi_1,\xi_2,\xi_3}$, define a quaternion-K\"{a}hler structure on the leaf space of $\mathscr{F}_{\xi_1,\xi_2,\xi_3}$.
\begin{thm}[\cite{BoyGal99}]\label{thm:QK-3-Sasak}
Let $(\mathcal{S},g)$ be a compact {3-Sasakian} manifold of dimension $n=4m+3$.  Then there is a natural
quaternion-K\"{a}hler structure on the leaf space of $\mathcal{F}_{\xi_1,\xi_2,\xi3}$, $(\mathcal{M},\check{g})$, such that
the orbifold map $\varpi:\mathcal{S}\rightarrow\mathcal{M}$ is a Riemannian submersion.
Furthermore, $(\mathcal{M},\check{g})$ is Einstein with scalar curvature \mbox{$s_{\check{g}}=16m(m+2))$}.
\end{thm}

\begin{figure}[tbp]\label{fig:3-Sasaki-geo}
\begin{center}
\setlength{\unitlength}{10pt}
\begin{picture}(12,12)
\put (5,1){\makebox(1,1){$\mathcal{M}$}}
\put (1.25,5){\makebox(1,1){$\mathcal{S}$}}
\put (8.75,5){\makebox(1,1){$\mathcal{Z}$}}
\put (5,9){\makebox(1,1){$C(\mathcal{S})$}}

\put (3,5.5){\vector(1,0){5}}
\put (4.5,8.5){\vector(-1,-1){2}}
\put (6.5,8.5){\vector(1,-1){2}}
\put (2.5,4.5){\vector(1,-1){2}}
\put (8.5,4.5){\vector(-1,-1){2}}

\put (2.5,7.8) {$\Sm{\R_+}$}
\put (8,7.8) {$\Sm{\C^*}$}
\put (5,6) {$\Sm{S^1}$}
\put (1.5,3) {$\Sm{\Sp(1)}$}\put (1.5,2.3) {$\Sm{\SO(3)}$}
\put (7.8,3) {$\Sm{\cps^1}$}

\end{picture}
\caption{Related geometries}
\end{center}
\end{figure}
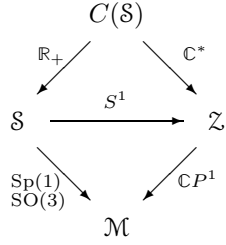

The geometries associated to a {3-Sasaki} manifold can be seen in Figure~\ref{fig:3-Sasaki-geo}.  Up to a finite cover, from
each space in Figure~\ref{fig:3-Sasaki-geo} the other three spaces can be recovered.  Unlike $\mathcal{S}$ the spaces $\mathcal{Z}$
and $\mathcal{M}$ are smooth in no more than finitely many cases for each $n\geq 1$.  Furthermore, it is proved~\cite{LeBSal94}
that the only smooth $\mathcal{M}^{4n}$ for $n=1,2$ are symmetric spaces.  That this is true for all $n$ is the famous
LeBrun-Salamon conjecture.

We will need to distinguish when the fibering $\mathcal{S}\rightarrow\mathcal{M}$ has generic fiber $\Sp(1)$.
The obstruction to this is the Marchiafava-Romani class.
An almost quaternionic structure $\mathcal{Q}$ is a reduction of the frame bundle
to an $Sp(1)Sp(m)$ bundle.  Let $\mathcal{G}$ be the sheaf of germs of smooth maps to
$Sp(1)Sp(m)$.  An almost quaternionic structure is an element $s\in H^1_{orb}(\mathcal{M},\mathcal{G})$.
Consider the exact sequence
\begin{equation}
0\rightarrow\Z_2\rightarrow Sp(1)\times Sp(m)\rightarrow Sp(1)Sp(m)\rightarrow 1.
\end{equation}
\begin{defn}\label{3-Sasak:defn-MRclass}
The \emph{Marchiafava-Romani class} is $\varepsilon =\delta(s)$, where
\[ \delta: H^1_{orb}(\mathcal{M},\mathcal{G})\rightarrow H^2_{orb}(\mathcal{M},\Z_2) \]
is the connecting homomorphism.
\end{defn}
One has that $\varepsilon$ is the orbifold Stiefel-Whitney class $w_2(\mathcal{Q})$.
Also, $\varepsilon$ is the obstruction to the existence of a square root $\mathbf{L}^{\frac{1}{2}}$
of $\mathbf{L}$. In the four-dimensional case $n=1$, $\varepsilon = w_2(\Lambda_+^2)=w_2(T\mathcal{M})$.
When $\varepsilon=0$ for the {3-Sasakian} space $\mathcal{S}$ associated to $(\mathcal{M},\check{g})$ we
will always mean the one with $Sp(1)$ generic fibres.

\subsection{Toric 3-Sasaki manifolds}\label{subsec:tor-3-Sasaki}

A {3-Sasaki} manifold $\mathcal{S}$ with $\dim\mathcal{S}=4m+3$ is \emph{toric} if it admits an effective action of
$T^{m+1}\subset\Aut(\mathcal{S},g)$, where $\Aut(\mathcal{S},g)$ is the group of 3-Sasaki automorphisms,
that is isometries preserving $(g,\eta_i,\xi_i,\phi_i),\ i=1,2,3$.  Equivalently, $C(\mathcal{S})$ is a toric hyper-K\"{a}hler
manifold~\cite{BieDan00}.
We will consider toric 3-Sasaki 7-manifolds which were constructed by a 3-Sasaki reduction procedure in~\cite{BGMR98}.
This constructs infinitely many smooth 3-Sasaki 7-manifolds for each $b_2 \geq 1$.
Subsequently it was proved by R. Bielawski~\cite{Bie99} that up to a finite cover all toric examples are obtained this way.

Let $\Aut(\mathcal{S},g)$ be the group 3-Sasaki automorphisms, that is isometries preserving $(g,\eta_i,\xi_i,\phi_i),\ i=1,2,3$.
Given a compact $G\subset\Aut(\mathcal{S},g)$ one can define the \emph{\mbox{3-Sasakian} moment map}
\begin{equation}
\mu_{\mathcal{S}}:\mathcal{S}\rightarrow\mathfrak{g}^*\otimes\R^3,
\end{equation}
where if $\tilde{X}$ is the vector field on $\mathcal{S}$ induced by $X\in\mathfrak{g}$ we have
\begin{equation}
\langle\mu^a_{\mathcal{S}}, X\rangle = \frac{1}{2}\eta^a(\tilde{X}),\quad  a=1,2,3\ \text{ for }X\in\mathfrak{g}.
\end{equation}
There is a quotient similar to the Marsden-Weinstein quotient of symplectic manifolds~\cite{BoyGalMan94}.
If a connected compact $G\subset\Aut(\mathcal{S},g)$ acts freely (locally freely) on $\mu_\mathcal{S}^{-1}(0)$,then
\[\mathcal{S}/\negthickspace/ G = \mu_\mathcal{S}^{-1}(0)/G \]
has the structure of a 3-Sasakian manifold (orbifold).

Consider the unit sphere $S^{4n-1}\subset\Ha^n$ with the round metric $g$ and the standard 3-Sasakian structure
induced by the right action of $Sp(1)$.
Then $\Aut(S^{4n-1},g)=Sp(n)$ acting by the standard linear representation on the left.
We have the maximal torus $T^n\subset Sp(n)$ and every representation of a subtorus
$T^k$ is conjugate to an inclusion $\iota_\Omega :T^k\rightarrow T^n$ which is
represented by a \emph{weight matrix} $\Omega=(a^i_j)\in\mathscr{M}_{k,n}(\Z)$, an integral $k\times n$ matrix.

Let $\{e_i\},i=1,\ldots,k$ be a basis for the dual of the Lie algebra of $T^k$, $\mathfrak{t}_k^*\cong\R^k$.
Then the moment map $\mu_\Omega:S^{4n-1}\rightarrow\mathfrak{t}_k^*\otimes\R^3$ can be written
as $\mu_\Omega =\sum_j\mu_\Omega^j e_j$ where in terms of complex coordinates $z_l+w_l j$ on $\Ha^n$ we have
\begin{equation}\label{eq:3-Sasak-mom}
\mu_\Omega^j(\mathbf{z},\mathbf{w})=i\sum_l a_l^j(|z_l|^2 - |w_l|^2)+2k\sum_l a_l^j\ol{w}_lz_l.
\end{equation}

We assume $\rank(\Omega)=k$ otherwise we just have an action of a subtorus of $T^k$.
Denote by
\begin{equation}
\Delta_{\alpha_1,\ldots,\alpha_k} = \det\left\lgroup
\begin{matrix}
a_{\alpha_1}^1 & \cdots & a_{\alpha_k}^1 \\
\vdots &         & \vdots \\
a_{\alpha_1}^k & \cdots & a_{\alpha_k}^k\\
\end{matrix}\right\rgroup
\end{equation}
the $\binom{n}{k}$ $k\times k$ minor determinants of $\Omega$.
\begin{defn}
Let $\Omega\in\mathscr{M}_{k,n}(\Z)$ be a weight matrix.
\begin{thmlist}
\item $\Omega$ is \emph{non-degenerate} if $\Delta_{\alpha_1,\ldots,\alpha_k}\neq 0$,
for all $1\leq\alpha_1 < \cdots < \alpha_k \leq n$.

\item  Let $\Omega$ be non-degenerate, and let $d$ be the $\gcd$ of all the $\Delta_{\alpha_1,\ldots,\alpha_k}$,
the kth determinantal divisor.  Then $\Omega$ is \emph{admissible} if
\[ \gcd(\Delta_{\alpha_2,\ldots,\alpha_{k+1}},\ldots,\Delta_{\alpha_1,\ldots,\hat{\alpha}_t,\ldots,\alpha_{k+1}},
\ldots,\Delta_{\alpha_1,\ldots,\alpha_k})=d\]
for all length $k+1$ sequences $1\leq\alpha_1 <\cdots <\alpha_t <\cdots <\alpha_{k+1}\leq n+1$.
\end{thmlist}
\end{defn}
The $\gcd$ $d_j$ of the jth row of $\Omega$ divides $d$.
We may assume that the $\gcd$ of each row of $\Omega$ is 1 by merely reparametrizing the coordinates
$\tau_j$ on $T^k$.  We say that $\Omega$ is in reduced form if $d=1$.

Choosing a different basis of $\mathfrak{t}_k$ results in an action on $\Omega$ by an element in $Gl(k,\Z)$.
We also have the normalizer of $T^n$ in
$Sp(n)$, the Weyl group $\mathscr{W}(Sp(n))=\Sigma_n \times\Z_2^n$ where $\Sigma_n$ is the permutation group.
$\mathscr{W}(Sp(n))$ acts on $S^{4n-1}$ preserving the \mbox{3-Sasakian} structure, and it acts on
weight matrices by permutations and sign changes of columns.
Thus the group $Gl(k,\Z)\times\mathscr{W}(Sp(n))$ acts on $\mathscr{M}_{k,n}(\Z)$, with the quotient only
depending on the equivalence class.
\begin{thm}[\cite{BGMR98}]
Let $\Omega\in\mathscr{M}_{k,n}(\Z)$ be reduced.
\begin{thmlist}
\item If $\Omega$ is non-degenerate, then $\mathcal{S}_\Omega$ is an orbifold.

\item  Supposing $\Omega$ is non-degenerate, $\mathcal{S}_\Omega$ is smooth if and only if $\Omega$ is admissible.
\end{thmlist}
\end{thm}
The quotient $\mathcal{S}_\Omega$ is toric, because its automorphism group contains $T^{n-k}\cong T^n/\iota_\Omega(T^k)$.

We are primarily interested in 7-dimensional toric quotients.  In this
case there are infinite families of distinct quotients.  We may take
matrices of the form
\begin{equation}\label{eq:3-Sasak-mat}
\Omega =\left\lgroup
\begin{matrix}
1 & 0 & \cdots & 0 & a_1 & b_1 \\
0 & 1 & \cdots & 0 & a_2 & b_2 \\
\vdots & \vdots & \ddots & \vdots & \vdots & \vdots \\
0 & 0 & \cdots & 1 & a_k & b_k
\end{matrix}
\right\rgroup.
\end{equation}

\begin{prop}[\cite{BGMR98}]\label{prop:3-Sasak-mat}
Let $\Omega\in\mathscr{M}_{k,k+2}(\Z)$ be as above.  Then $\Omega$ is admissible if and only if
$a_i,b_j,i,j=1,\ldots, k$ are all nonzero, $\gcd(a_i,b_i)=1$ for $i=1,\ldots,k$, and we do not have
$a_i= a_j$ and $b_i=b_j$, or $a_i=-a_j$ and $b_i=-b_j$ for some $i\neq j$.
\end{prop}

Proposition (\ref{prop:3-Sasak-mat}) shows that for $n=k+2$ there are infinitely many reduced admissible
weight matrices.  One can, for example, choose $a_i,b_j,i,j=1,\ldots k$ be all pairwise relatively prime.
We will make use of the cohomology computation of R. Hepworth~\cite{Hep07} to show that we have infinitely many
smooth 3-Sasakian 7-manifolds of each second Betti number $b_2\geq 1$.
Let $\Delta_{p,q}$ denote the $k\times k$ minor determinant of $\Omega$ obtained by deleting the $p^{th}$ and $q^{th}$
columns.

\begin{thm}[\cite{BGMR98,Hep07}]\label{thm:3-Sasak-coh}
Let $\Omega\in\mathscr{M}_{k,k+2}(\Z)$ be a reduced admissible weight matrix.  Then
$\pi_1(\mathcal{S}_\Omega)=e$.  And the cohomology of $\mathcal{S}_\Omega$ is

\renewcommand{\arraystretch}{1.5}
\begin{centering}
\begin{tabular}{l|cccccccc}
$p$ & $0$ & $1$ & $2$ & $3$ & $4$ & $5$ & $6$ & $7$\\ \hline
$H^p$ & $\Z$ & $0$ & $\Z^k$ & $0$ & $G_\Omega$ & $\Z^k$ & $0$ & $\Z$ \\
\end{tabular},\\
\end{centering}
where $G_\Omega$ is a torsion group of order
\[ \sum |\Delta_{s_1, t_1}|\cdots |\Delta_{s_{k+1},t_{k+1}}| \]
with the summand with index $s_1,t_1,\ldots, s_{k+1},t_{k+1}$ included if and only if
the graph on the vertices $\{1,\ldots, k+2\}$ with edges $\{s_i,t_i\}$ is a tree.
\end{thm}

If we consider weight matrices as in Proposition (\ref{prop:3-Sasak-mat}) then the order
of $G_\Omega$ is greater than $|a_1\cdots a_k|+|b_1\cdots b_k|$.
We have the following.
\begin{cor}[\cite{BGMR98,Hep07}]\label{cor:3-Sasak-coh}
There are smooth toric 3-Sasakian 7-manifolds with second Betti number $b_2=k$ for all $k\geq 0$.
Furthermore, there are infinitely many possible homotopy types of examples $\mathcal{S}_{\Omega}$ for each $k>0$.
\end{cor}

Note that the reduction procedure can be done on any of the four spaces in~\ref{fig:3-Sasaki-geo}.  In particular,
we have the $\ASD$ Einstein orbifold $\mathcal{M}_\Omega =\mathcal{S}_{\Omega}/\Sp(1)$, which is a quaternic-K\"{a}hler
quotient~\cite{GalLaw88} of $\qps^{n-1}$ by the torus $T^k$.  An $\ASD$ Einstein orbifold $\mathcal{M}$ is \emph{toric}
if it has an effective isometric action of $T^2$.

Recall the orbifold $\mathcal{M}_\Omega$ has an action of $T^2\cong T^{k+2}/\iota_\Omega(T^k)$,
and can be characterized as in~\cite{OrlRay70,HaeSal91} by its orbit space and stabilizer groups.
The stabilizers were determined in~\cite{BoyGalMan98,CalSin06}.
The orbit space is $Q_\Omega:= \mathcal{M}_\Omega /{T^2}$.
Then $Q_\Omega$ is a polygon with $k+2$ edges $C_1,C_2,\ldots, C_{k+2}$,
labeled in cyclic order with the interior of $C_i$ being orbits with stabilizer $G_i$, where
$G_i \subset T^2,\  i=1,\ldots, k+2$ are $S^1$ subgroups.
Choose an explicit surjective homomorphism
$\Phi:\Z^{k+2}\rightarrow\Z^2$ annihilating the rows of $\Omega$.  So
\begin{equation}
\Phi =\left\lgroup
\begin{matrix}
b_1 & b_2 & \cdots & b_{k+2} \\
c_1 & c_2 & \cdots & c_{k+2}
\end{matrix}
\right\rgroup
\end{equation}
It will be helpful to normalize $\Phi$.  After acting on the columns of $\Phi$ by $\mathscr{W}(Sp(k+2))$ and
on the right by $Gl(2,\Z)$ we may assume that $b_i>0$ for $i=1,\ldots,k+2$ and
$c_1/b_1 < \cdots <c_i/b_i <\cdots < c_{k+2}/b_{k+2}$.
Then the stabilizer groups $G_i\subset T^2$ are characterized by $(m_i,n_i)\in\Z^2$ where
\begin{equation}
(m_i,n_i)=\sum_{l=1}^{i}(b_l,c_l)-\sum_{l=i+1}^{k+2}(b_l,c_l),\quad i=1,\cdots k+2.
\end{equation}
It is convenient to take $(m_0,n_0)=-(m_{k+2},n_{k+2})$.

\section{K\"{a}hler-Einstein metric on symmetric Fano orbifolds}\label{sec:K-E-sym}

We prove the existence of the K\"{a}hler-Einstein metric on symmetric toric Fano orbifolds in this section.  The existence of
a K\"{a}hler-Einstein metric on a toric Fano manifold with vanishing Futaki invariant was proved by X. Wang and X. Zhu~\cite{WanZhu04}.
More generally, they proved the existence of a K\"{a}hler-Ricci soliton which is K\"{a}hler-Einstein if the Futaki invariant
vanishes.  Then A. Futaki, H. Ono, and G. Wang~\cite{FutOnWan09} proved an extension of that result, namely that any toric
Sasaki manifold, which satisfies the necessary positivity condition, admits a K\"{a}hler-Ricci soliton which is
Sasaki-Einstein if the transverse Futaki invariant vanishes.  This latter result includes the existence result proved here.
But the proof given here, as the proof in~\cite{BatSel99} for symmetric toric Fano manifolds, shows that the invariant
of G. Tian ~\cite{Ti89}, extended by J.-P. Demailly and J. Koll{\'a}r~\cite{DemKol01} to orbifolds,
satisfies $\alpha_G(X)\geq 1$, which in this case is an invariant of $X$ as a Fano orbifold.

\subsection{Symmetric Fano orbifolds}

Let $N\cong\Z^r$ be the free $\Z$-module of rank r and $M=\Hom_\Z(N,\Z)$ its dual.
We denote $N_\Q =N\otimes\Q$ and $M_\Q=M\otimes\Q$ with the natural pairing
\[\langle\ \, ,\ \rangle : M_\Q \times N_\Q \rightarrow\Q.\]
Similarly we denote $N_\R =N\otimes\R$ and $M_\R =M\otimes\R$.

Let $T_\C:=N\otimes_{\Z}\C^*\cong(\C^*)^n$ be the algebraic torus.
Each $m\in M$ defines a character $\chi^m :T_\C \rightarrow\C^*$ and each $n\in N$ defines a
one-parameter subgroup $\lambda_n :\C^*\rightarrow T_\C$.  In fact, this gives an isomorphism between
$M$ (resp. $N$) and the multiplicative group $\Hom_{\text{alg.}}(T_\C,\C^*)$
(resp. $\Hom_{\text{alg.}}(\C^*,T_\C)$).

An n-dimensional toric variety $X$ has $T_\C \subseteq\Aut(X)$ with an open dense orbit $U\subset X$.  Then $X$ is
defined by a fan $\Delta$ in $N_\Q$.  We denote this $X_\Delta$.  See~\cite{Ful93} or~\cite{Od88} for background on toric
varieties.  We denote by $\Delta(i)$ the set of $i$-dimensional cones in $\Delta$.

Recall that every element $\rho\in\Delta(1)$ is generated by a unique primitive element of $N$.  We will consider
non-primitive generators to encode an orbifold structure.
\begin{defn}
We will denote by $\Delta^*$ an \emph{augmented fan} by which we mean a fan $\Delta$ with elements
$n(\rho) \in N\cap\rho$ for every $\rho\in\Delta(1)$.
\end{defn}
\begin{prop}\label{prop:toric-orb}
For a complete simplicial augmented fan $\Delta^*$ we have a natural orbifold structure compatible with the action of $T_\C$
on $X_\Delta$.  We denote $X_\Delta$ with this orbifold structure by $X_{\Delta^*}$.
\end{prop}
\begin{proof}
Let $\sigma\in\Delta^*(n)$ have generators $p_1,p_2,\ldots,p_n$ as in the definition.  Let
$N^\prime\subseteq N$ be the sublattice $N^\prime =\Z\{p_1,p_2,\ldots,p_n\}$, and $\sigma^\prime$ the
equivalent cone in $N^\prime$.  Denote by $M^\prime$ the dual lattice of $N^\prime$ and $T^\prime_\C$ the torus.
Then $U_{\sigma^\prime}\cong\C^n$.  It is easy to see that
\[ N/N^\prime = \Hom_\Z(M^\prime/M,\C^*). \]
And $N/N^\prime$ is the kernel of the homomorphism
\[T^\prime_\C=\Hom_\Z(M^\prime,\C^*)\rightarrow T_\C=\Hom_\Z(M,\C^*).\]
Let $\Gamma=N/N^\prime$.  An element $t\in\Gamma$ is a homomorphism $t:M^\prime\rightarrow\C^*$ equal to 1
on $M$.  The regular functions on $U_{\sigma^\prime}$ consist of $\C$-linear combinations of $x^m$ for
$m\in\sigma^{\prime\vee}\cap M^\prime$.  And $t\cdot x^m =t(m)x^m$.  Thus the invariant functions are
the $\C$-linear combinations of $x^m$ for $m\in\sigma^{\vee}\cap M$, the regular functions of $U_\sigma$.
Thus $U_{\sigma^\prime}/\Gamma = U_\sigma$.  And the charts are easily seen to be compatible on intersections.
\end{proof}
Conversely one can prove that this definition gives all structures of interest.
\begin{prop}
Let $\Delta$ be a complete simplicial fan.  Suppose for simplicity that the local uniformizing groups are
abelian.  Then every orbifold structure on $X_\Delta$ compatible with the action of $T_\C$ arises from an
augmented fan $\Delta^*$.
\end{prop}
Note that these orbifold structures are not well formed, i.e. have complex codimension one singular sets.
That is, if some $n(\rho)=a_\rho p_\rho,\ a_\rho \in\N_{>1},$ is not primitive, then the divisor $D_\rho$ has a cone angle
of $2\pi/a_\rho$.  This has no significance for the complex structure, but compatible metrics will have this cone singularity.

We modify the usual definition of a support function to characterize orbifold line bundles on $X_{\Delta^*}$.
We will assumed from now on that the fan $\Delta$ is simplicial
and complete.
\begin{defn}
A real function $h:N_\R \rightarrow\R$ is a \emph{$\Delta^*$ -linear support function} if for
each $\sigma\in\Delta^*$ with given $\Q$-generators $p_1,\ldots,p_r$ in $N$, there is an
$l_\sigma\in M_\Q$ with $h(s)=\langle l_\sigma ,s\rangle$ and $l_\sigma$ is $\Z$-valued on the sublattice
$\Z\{p_1,\ldots,p_r\}$.  And we require that $\langle l_\sigma, s\rangle=\langle l_\tau ,s\rangle$ whenever
$s\in\sigma\cap\tau$.  The additive group of $\Delta^*$-linear support functions will be denoted by
$\SF(\Delta^*)$.
\end{defn}
Note that $h\in\SF(\Delta^*)$ is completely determined by the integers $h(n(\rho))$ for all
$\rho\in\Delta(1)$.  And conversely, an assignment of an integer to $h(n(\rho))$ for all
$\rho\in\Delta(1)$ defines $h$.  Thus
\[ \SF(\Delta^*)\cong\Z^{\Delta(1)}. \]

\begin{defn}
Let $\Delta^*$ be a complete augmented fan.  For $h\in\SF(\Delta^*)$,
\[ \Sigma_h :=\{m\in M_\R : \langle m,n\rangle \geq h(n),\text{ for all }n\in N_\R\}, \]
is a, possibly empty, convex polytope in $M_\R$.
\end{defn}

Recall that a certain subset of $\Q$-Weil divisors correspond to orbifold line bundles.
\begin{defn}
A \emph{Baily divisor} is a $\Q$-Weil divisor $D\in\weil(X)\otimes\Q$ whose inverse image
$D_{\tilde{U}}\in\weil(\tilde{U})$ in every local uniformizing chart $\pi:\tilde{U}\rightarrow U$
is Cartier.  The additive group of Baily divisors is denoted $\Divorb(X)$.
\end{defn}
A Baily divisor $D$ defines a holomorphic orbifold line bundle $[D]\in\Picorb(X)$ in a way completely analogous
to Cartier divisors.  We denote the Baily divisors invariant under $T_{\C}$ by $\Divorb_{T_\C}(X)$.
We denote the group of isomorphism equivariant orbifold line bundles by $\Picorb_{T_\C}(X)$.  Then likewise
we have $[D]\in\Picorb_{T_\C}(X)$ whenever $D\in\Divorb_{T_\C}(X)$.

A straight forward generalization of~\cite[Prop. 2.1]{Od88} to this situation gives the following.
\begin{prop}\label{prop:toric-vbund}
Let $X=X_{\Delta^*}$ be compact with the standard orbifold structure, i.e. $\Delta^*$ is simplicial and complete.
\begin{thmlist}
\item  There is an isomorphism $\SF(\Delta^*)\cong\Divorb_{T_\C}(X)$ obtained by sending \label{prop:toric-vbund-div}
$h\in\SF(\Delta^*)$ to
\[ D_h :=-\sum_{\rho\in\Delta(1)}h(n(\rho)) D_\rho,\]
where $D_\rho$ is the divisor of $X$ associated to $\rho\in\Delta(1)$.
\item  There is a natural homomorphism $\SF(\Delta^*)\rightarrow\Picorb_{T_\C}(X)$ which associates
an equivariant orbifold line bundle $\mathbf{L}_h$ to each $h\in\SF(\Delta^*)$.
\item  Suppose $h\in\SF(\Delta^*)$ and $m\in M$ satisfies
\[ \langle m,n\rangle\geq h(n)\text{  for all  }n\in N_\R,\]
then $m$ defines a section $\psi:X\rightarrow\mathbf{L}_h$ which has the equivariance  property
$\psi(tx)=\chi^m(t)(t\psi(x))$.
\item  The set of sections $H^0(X,\mathcal{O}(\mathbf{L}_h))$ is the finite dimensional
$\C$-vector space with basis $\{x^m :m\in\Sigma_h \cap M\}$.
\item  Every Baily divisor is linearly equivalent to a $T_\C$-invariant Baily divisor.
Thus for $D\in\Picorb(X)$, $[D]\cong[D_h]$ for some $h\in\SF(\Delta^*)$.
\item  If $\mathbf{L}$ is any holomorphic orbifold line bundle, then $\mathbf{L}\cong\mathbf{L}_h$
for some $h\in\SF(\Delta^*)$.  The homomorphism in part i. induces an isomorphism
$\SF(\Delta^*)\cong\Picorb_{T_\C}(X)$ and we have the exact sequence
\[ 0\rightarrow M\rightarrow\SF(\Delta^*)\rightarrow\Picorb(X)\rightarrow 1. \]
\end{thmlist}
\end{prop}
\begin{rmk}
The notation is a bit deceptive that in (\ref{prop:toric-vbund-div}) it appears that $D_h$ is a $\Z$-Weil divisor.
But they are written with their coefficients in the uniformizing chart, and the components in ramification divisors of
the chart are generally fractional.
\end{rmk}

For $X=X_{\Delta^*}$ there is a unique $k\in\SF(\Delta^*)$ such that $k(n(\rho))=1$ for all $\rho\in\Delta(1)$.
The corresponding Baily divisor
\begin{equation}
D_k :=-\sum_{\rho\in\Delta(1)}D_\rho
\end{equation}
is the \emph{orbifold canonical divisor}.  The corresponding orbifold line bundle is $\mathbf{K}_X$, the orbifold bundle
of holomorphic n-forms.  This will in general be different from the canonical sheaf in the algebraic geometric
sense.

\begin{defn}
Consider support functions as above but which are only required to be $\Q$-valued on $N_\Q$, denoted $\SF(\Delta,\Q)$.
$h$ is \emph{strictly upper convex} if $h(n+n^\prime )\geq h(n)+h(n^\prime)$ for all $n,n^\prime \in N_\Q$
and for any two $\sigma,\sigma^\prime \in\Delta(n)$, $l_\sigma$ and $l_{\sigma^\prime}$ are different linear functions.
\end{defn}
Given a strictly upper convex support function $h$, the polytope $\Sigma_h$ is the convex hull in $M_\R$ of the
vertices $\{l_\sigma : \sigma\in\Delta(n)\}$.  Each $\rho\in\Delta(1)$ defines a facet by
\[ \langle m, n(\rho)\rangle\geq h(n(\rho)). \]
If $n(\rho)=a_\rho n^\prime$ with $n^\prime \in N$ primitive and $a_\rho\in\N$ we may label the face with
$a_\rho$ to get the labeled polytope $\Sigma^*_h$ which encodes the orbifold structure.  Conversely, from
a rational convex polytope $\Sigma^*$ we associate a fan $\Delta^*$ and a support function $h$.
\begin{prop}[\cite{Od88,Ful93}]
There is a one-to-one correspondence between the set of pairs $(\Delta^*,h)$ with $h\in\SF(\Delta,\Q)$ strictly
upper convex, and rational convex marked polytopes $\Sigma^*_h$.
\end{prop}

\begin{defn}\label{defn:toric-Fano}
Let $X=X_{\Delta^*}$ be a compact toric orbifold.  We say that $X$ is \emph{Fano} if $-k\in\SF(\Delta^*)$,
which defines the anti-canonical orbifold line bundle $\mathbf{K}^{-1}_X$, is strictly upper convex.
\end{defn}
These toric variety aren't necessarily Fano in the usual sense, since $\mathbf{K}^{-1}_X$ is the
\emph{orbifold} anti-canonical class.  This condition is equivalent to $\{n\in N_\R : k(n)\leq 1 \}\subset N_\R$
being a convex polytope with vertices $n(\rho), \rho\in\Delta(1)$.  We will use $\Delta^*$ to denote both the
augmented fan and this polytope in this case.

\subsection{Symmetric toric varieties}

Let $X_\Delta$ be an $n$-dimensional toric variety.
Let $\mathcal{N}(T_{\C})\subset\Aut(X)$ be the normalizer of $T_{\C}$.  Then
$\mathcal{W}(X):=\mathcal{N}(T_{\C})/T_{\C}$ is isomorphic to the finite group of all symmetries of $\Delta$,
i.e. the subgroup of $GL(n,\Z)$ of all $\gamma\in GL(n,\Z)$ with $\gamma(\Delta)=\Delta$.
Then we have the exact sequence.
\begin{equation}\label{toric:wely-exact1}
1\rightarrow T_{\C}\rightarrow\mathcal{N}(T_{\C}) \rightarrow\mathcal{W}(X)\rightarrow 1.
\end{equation}
Choosing a point $x\in X$ in the open orbit, defines an inclusion $T_{\C}\subset X$.  This also
provides a splitting of (\ref{toric:wely-exact1}).
Let $\mathcal{W}_0(X)\subseteq\mathcal{W}(X)$ be the subgroup which are also automorphisms of $\Delta^*$;
$\gamma\in\mathcal{W}_0(X)$ is an element of $\mathcal{N}(T_{\C})\subset\Aut(X)$ which preserves the orbifold
structure.  Let $G\subset\mathcal{N}(T_{\C})$ be the compact subgroup generated by $T^n$, the maximal compact subgroup
of $T_\C$, and $\mathcal{W}_0(X)$.
Then we have the, split, exact sequence
\begin{equation}\label{weyl:exact2}
1\rightarrow T^n\rightarrow G\rightarrow\mathcal{W}_0(X)\rightarrow 1.
\end{equation}
\begin{defn}
A \emph{symmetric Fano toric orbifold} $X$ is a Fano toric orbifold with
$\mathcal{W}_0$ acting on $N$ with the origin as the only fixed point.
Such a variety and its orbifold structure is characterized by the convex polytope $\Delta^*$ invariant under $\mathcal{W}_0$.
We call a toric orbifold \emph{special symmetric} if $\mathcal{W}_0(X)$ contains the
involution $\sigma: N\rightarrow N$, where $\sigma(n)=-n$.
\end{defn}
Conversely, given an integral convex polytope $\Delta^*$, inducing a simplicial fan $\Delta$,
invariant under a subgroup $\mathcal{W}_0\subset GL(n,\Z)$ fixing
only the origin, we have a symmetric Fano toric orbifold $X_{\Delta^*}$.
\begin{defn}\label{defn:toric-ind}
The \emph{index} of a Fano orbifold $X$ is the largest positive integer $m$ such that there is a holomorphic
$V$-bundle $\mathbf{L}$ with $\mathbf{L}^m \cong\mathbf{K}_X^{-1}$. The index of $X$ is denoted
$\ind(X)$.
\end{defn}
Note that $c_1(X)\in H^2_{orb}(X,\Z)$, and $\ind(X)$ is the greatest positive integer $m$ such that
$\frac{1}{m}c_1(X)\in H^2_{orb}(X,\Z)$.
\begin{prop}\label{prop:toric-ind-sp}
Let $X_{\Delta^*}$ be a special symmetric toric Fano orbifold.  Then $\ind(X)=1$ or $2$.
\end{prop}
\begin{proof}
We have $\mathbf{K}^{-1}\cong\mathbf{L}_{-k}$ with $-k\in\SF(\Delta^*)$ where $-k(n_\rho)=-1$
for all $\rho\in\Delta(1)$.  Suppose we have $\mathbf{L}^m\cong\mathbf{K}^{-1}$.
By proposition (\ref{prop:toric-vbund}) there is an $h\in\SF(\Delta^*)$ and $f\in M$ so
that $mh=-k+f$.  For some $\rho\in\Delta(1)$,
\begin{gather*}
mh(n_\rho) = -1 + f(n_\rho) \\
mh(-n_\rho) = -1 - f(n_\rho).
\end{gather*}
Thus $m(h(n_\rho)+h(-n_\rho))=-2$, and $m=1$ or $2$.
\end{proof}

\begin{figure}[tb]
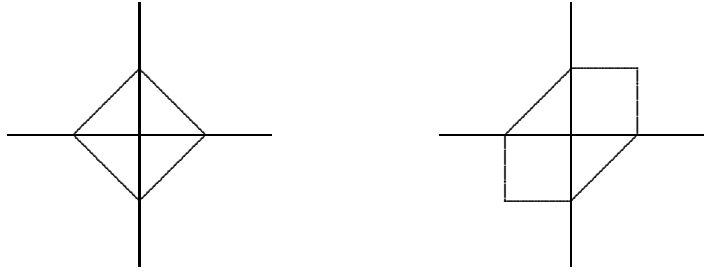
\label{fig:smooth-sur}
\centering

\mbox{
\beginpicture
\setcoordinatesystem units <.5pt,.5pt> point at 0 0
\setplotarea x from -100 to 100, y from -100 to 100
\axis bottom shiftedto x=0 /
\axis left shiftedto y=0 /
\setlinear
\plot 0 50   50 0  0 -50  -50 0  0 50 /
\endpicture

\hspace{60pt}

\beginpicture
\setcoordinatesystem units <.5pt,.5pt> point at 0 0
\setplotarea x from -100 to 100, y from -100 to 100
\axis bottom shiftedto x=0 /
\axis left shiftedto y=0 /
\setlinear
\plot 0 50  50 50  50 0  0 -50  -50 -50  -50 0  0 50 /
\endpicture}

\caption{Smooth examples}
\end{figure}

In the in the subsequent sections we will be interested in special symmetric toric Fano surfaces.
Figure~\ref{fig:smooth-sur} gives the polytopes $\Delta^*$ for the two smooth such examples.

\begin{figure}[tb]
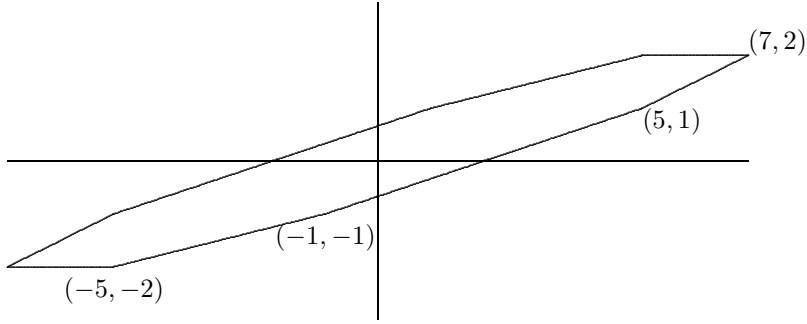
\label{fig:typ-xpl}
\centering
\mbox{
\beginpicture
\setcoordinatesystem units <.4pt,.4pt> point at 0 0
\setplotarea x from -350 to 350, y from -150 to 150
\axis bottom shiftedto x=0 /
\axis left shiftedto y=0 /
\setlinear
\plot 50 50  250 100  350 100  250 50  -50 -50  -250 -100  -350 -100  -250 -50  50 50 /
\put {$(7,2)$} [lb] at 350 100
\put {$(5,1)$} [lt] at 250 50
\put {$(-1,-1)$} [ct] at -50 -60
\put {$(-5,-2)$} [ct] at -250 -110
\endpicture}

\caption{Example with 8 point singular set and $\mathcal{W}_0=\Z_2$}
\end{figure}

\subsection{K\"{a}hler-Einstein metric}

Any compact toric orbifold associated
to a polytope admits a K\"{a}hler metric.  In particular, we need a K\"{a}hler metric with K\"{a}hler form $\omega$ satisfying
$[\omega]\in 2\pi c_1^{orb}(X)=-c_1(\mathbf{K}_X)$.  The Hamiltonian reduction procedure of~\cite{Gui94,Gui95} and~\cite{CalDaGau03}
provides an explicit metric on the toric orbifold associated to the marked polytope $\Sigma^*_h$.
Let $X_{\Sigma^*_{-k}}$ be Fano, it will follow that this metric will satisfy $[\omega]\in 2\pi c_1^{orb}(X)$.

Let $\Sigma^*$ be a convex polytope in $M_\R \cong{\R^n}^*$ defined by the inequalities
\begin{equation}
\langle x,u_i\rangle\geq\lambda_i, \quad i=1,\ldots, d,
\end{equation}
where $u_i \in N\subset N_\R\cong\R^n$ and $\lambda_i \in\R$.
If $\Sigma^*_h$ is associated to $(\Delta^*,h)$, then
the $u_i$ and $\lambda_i$ are the set of pairs $n(\rho)$ and $h(n(\rho))$ for $\rho\in\Delta(1)$.
We allow the $\lambda_i$ to be real but require any set $u_{i_1},\ldots, u_{i_n}$ corresponding to a
vertex to form a $\Q$-basis of $N_\Q$.

Let $(e_1,\ldots, e_d)$ be the standard basis of $\R^d$ and $\beta:\R^d\rightarrow\R^n$ be the map
which takes $e_i$ to $u_i$.  Let $\mathfrak{n}$ be the kernel of $\beta$, so we have the exact sequence
\begin{equation}\label{eq:toric-exact}
0\rightarrow\mathfrak{n}\overset{\iota}{\rightarrow}\R^d \overset{\beta}{\rightarrow}\R^n \rightarrow 0,
\end{equation}
and the dual exact sequence
\begin{equation}\label{eq:toric-exact-dual}
0\rightarrow{\R^n}^*\overset{\beta^*}{\rightarrow}{\R^d}^*\overset{\iota^*}{\rightarrow}\mathfrak{n}^*\rightarrow 0.
\end{equation}
Since (\ref{eq:toric-exact}) induces an exact sequence of lattices, we have an exact sequence
\begin{equation}
1\rightarrow N\rightarrow T^d \rightarrow T^n \rightarrow 1,
\end{equation}
where the connected component of the identity of $N$ is an $(d-n)$-dimensional torus.
The standard representation of $T^d$ on $\C^d$ preserves the K\"{a}hler form
\begin{equation}
\frac{i}{2}\sum_{k=1}^d dz_k \wedge d\ol{z}_k,
\end{equation}
and is Hamiltonian with moment map
\begin{equation}
\mu(z)=\frac{1}{2}\sum_{k=1}^d |z_k|^2 e_k +c,
\end{equation}
unique up to a constant $c$.  We will set $c=\sum_{k=1}^d \lambda_k e_k$.  Restricting
to $\mathfrak{n}^*$ we get the moment map for the action of $N$ on $\C^d$
\begin{equation}
\mu_N (z)=\frac{1}{2}\sum_{k=1}^d |z_k|^2 \alpha_k +\lambda,
\end{equation}
with $\alpha_k =\iota^*e_k$ and $\lambda=\sum\lambda_k\alpha_k$.

We have the Marsden-Weinstein quotient
\[ X_{\Sigma^*} =\mu_N^{-1}(0)/N \]
with a canonical metric with K\"{a}hler form $\omega_0$.  We have an action of $T^n =T^d/N$ on $X_{\Sigma^*}$ which is Hamiltonian for $\omega$.
The map $\nu$ is $T^d$ invariant, and it descends to a map, which we also call $\nu$,
\begin{equation}\label{eq:moment}
\nu:X_{\Sigma^*}\rightarrow {\R^n}^*,
\end{equation}
which is the moment map for this action.  The above comments show that
$\im(\nu)=\Sigma^*$.  The action $T^n$ extends to the complex torus $T^n_\C$ and one can show
that as an analytic variety and orbifold $X_{\Sigma^*}$ is the toric variety constructed from
$\Sigma^*$ in the previous section.

It follows from results of~\cite{Gui94,Gui95} that
\[ [\omega_0] = -2\pi\sum_{i=1}^d \lambda_i c_i ,\]
where $c_i \in H^(X,\R)$ is dual to the divisor $D_i \subset X$ associated with the face $\langle x,u_i\rangle =\lambda_i$
of $\Sigma^*$.  In particular, if $X_{\Sigma^*}$ is Fano, then $\lambda_i =-1,\ i=1,\ldots, d$ and
\[  [\omega_0] = 2\pi\sum_{i=1}^d c_i =2\pi c_1^{orb}(X).\]

From now on we assume that $X_{\Sigma^*}$ is symmetric and Fano, and we have a metric $g_0$ invariant under the compact
group $G\subset\Aut(X)$ with K\"{a}hler form $\omega_0$ representing $2\pi c_1^{orb}(X)$.
Finding a K\"{a}hler-Einstein metric on $X_{\Sigma^*}$ is equivalent to solving the complex Monge-Amp\`{e}re for
$\phi\in C^\infty(X)$:
\begin{equation}\label{eq:Monge-Amp}
(\omega_0 +\sqrt{-1}\partial\ol{\partial}\phi)^n =\omega_0^n e^{f-t\phi},\quad t\in [0,1],
\end{equation}
where $f\in C^\infty(X)$ is defined by
\[ \sqrt{-1}\partial\ol{\partial}f =\Ricci(\omega_0)-\omega_0\text{  and  } \int_X e^f\, d\mu_{g_0} =\Vol_{g_0}(X).\]
If $\phi$ is a solution to (\ref{eq:Monge-Amp}) for $t=1$, then
\[ \omega =\omega_0  +\sqrt{-1}\partial\ol{\partial}\phi\]
is K\"{a}hler-Einstein.  It is well known that a solution to (\ref{eq:Monge-Amp}) exists for $t\in [0,\epsilon)$ for
$\epsilon$ small, and the existence of a solution at $t=1$ is equivalent to an \emph{a priori} $C^0$
estimate on $\phi$.

We recall the definition of the invariant $\alpha_G(X)$ introduced by G. Tian~\cite{Ti89}.
Define
\[ P_G(X,g_0) :=\Bigl\{\phi\in C^2(X)^G\ |\ \omega_0 +\sqrt{-1}\partial\ol{\partial}\phi\geq 0\text{ and }\underset{X}{\sup}\ \phi=0 \Bigr\}.\]
The \emph{Tian invariant} $\alpha_G(X)$ is the supremum of $\alpha >0$ such that
\[ \int_X e^{-\alpha\phi} d\mu_{g_0} \leq C(\alpha),\quad\forall\, \phi\in P_G(X,g_0),\]
where $C(\alpha)$ depends only on $\alpha, X$ and $g_0$.

G. Tian proved the following sufficient condition for an \emph{a priori} $C^0$ estimate on (\ref{eq:Monge-Amp}).  It was
shown to also suffice for orbifolds in~\cite{DemKol01}.
\begin{thm}
Let $X$ be a Fano orbifold and $G\subset\Aut(X)$ a compact subgroup such that
\[ \alpha_G(X) >\frac{n}{n+1},\]
then $X$ admits a K\"{a}hler-Einstein metric.
\end{thm}

Choosing a point $x_0 \in U\subset X$ gives identifications $\mathcal{W}_0(X)\subset\Aut(X)$, $U\cong T_\C$, and
$U/T \cong N_\R$, which identifies $Tx_0$ with $0\in N_R$.  Thus $\mathcal{W}_0(X)$ acts linearly on $N_\R$.  And if we
choose an integral basis $e_1,\ldots, e_n$ of $N$, then we have identifications $N_\R \cong\R^n$, $M_R \cong\R^n$, and
$T_\C \cong (\C^*)^n$.  And we introduce logarithmic coordinates $x_i =\log|t_i|^2$ on $N_\R$, where $t_1,\ldots,t_n$ are the usual
holomorphic coordinates on $(\C^*)^n$.  Thurs $t_i =e^{\frac{1}{2}x_i +\sqrt{-1}\theta_i}$, where $0\leq\theta_i \leq 2\pi$.
We will denote the dual coordinates on $M_\R$ by $y_1,\ldots,y_n$.  We define $l_k(y) =\langle u_k, y\rangle -\lambda_k,\
k=1,\ldots,d$.  So $\Sigma$ is defined by $\cap_{k=1}^d \{l_k \geq 0\}$.

Since the action of $T$ on $U$ is Hamiltonian for $\omega_0$, the orbits of $T$ are isotropic and $\omega_0 |_U$ is exact.
Furthermore, since $H^{0,k}(U) =H^{0,k}(U)^T =0$, we easily get the following.
\begin{lem}
The K\"{a}hler form $\omega_0$ restricted to $U$ has $T$-invariant potential function.  That is, there is a $F\in C^\infty(N_\R)$
with
\[ \omega_0 |_U =\sqrt{-1}\partial\ol{\partial} F.\]
\end{lem}
It was observed in~\cite{Gui94} that up to a constant the moment map (\ref{eq:moment}) is
\[ \nu :N_\R \rightarrow M_\R,\]
\begin{equation}\label{eq:Leg-moment}
\nu(x_1,\ldots,x_n) =\Bigl(\frac{\partial F}{\partial x_1}(x),\cdots,\frac{\partial F}{\partial x_n}(x)  \Bigr).
\end{equation}
By replacing $F$ with $F+\sum_k c_k x_k$ if necessary, we have that (\ref{eq:Leg-moment}) coincides with (\ref{eq:moment})
restricted to $U$.  Therefore, it is a diffeomorphism of $N_\R$ onto the interior of $\Sigma$.

It was shown in~\cite{Gui94} that the symplectic potential $G$ of the metric $\omega_0$ is
\begin{equation}\label{eq:symp-pot}
G=\frac{1}{2}\sum_{k=1}^d l_k(y)\log l_k(y),
\end{equation}
where it was also shown that $F$ and $G$ are related by the Legendre transform.  As a consequence we get
\begin{equation}
F=\nu^*\Bigr(\frac{1}{2}\sum_{k=1}^d \lambda_k \log l_k +l_\infty \Bigr),
\end{equation}
where $l_\infty (y) =\langle u_\infty ,y\rangle,\ u_\infty =\sum_{k=1}^d u_k$.  It is easy to see that the symmetric condition on
$X_{\Sigma^*}$ implies $u_\infty =0$.

Let $\sigma_j ,\ j=1,\ldots,e$ be the vertices of $\Sigma$.  Thus for each element of $\Delta^*(n)$ defined by
$u_{j_1},\ldots, u_{j_n}$ one has that $\sigma_j$ is the unique linear function with $\sigma_j (u_{j_i})=-1,\ i=1,\ldots,n$.
Recall that $\lambda_i =-1,\ i=1,\ldots,d$.
We define the piecewise linear function on $N_\R$
\begin{equation}
\ol{w}(x):= \underset{j=1,\ldots,e}{\sup}\langle\sigma_j, x\rangle.
\end{equation}

\begin{lem}\label{lem:bound}
There exists a constant $C>0$, depending only on $\Sigma^*$, so that
\[ |F-\ol{w}|\leq C.\]
\end{lem}
\begin{proof}
We prove this on momentum coordinates on $M_\R$.  The moment map is inverted by $x_i =\frac{\partial G}{\partial y_i}$ for
$i=1,\ldots, n$.  And one computes
\[ \begin{split}
\frac{\partial G}{\partial y}(y) & = \frac{1}{2}\sum_{k=1}^d u_k \log l_k(y) +u_\infty \\
                                & = \frac{1}{2}\sum_{k=1}^d u_k \log l_k(y)
\end{split}\]
Thus on the interior of $\Sigma$,
\begin{equation}
\ol{w}(y) = \underset{j=1,\ldots,e}{\sup}\frac{1}{2}\sum_{k=1}^d \langle\sigma_j,u_k\rangle\log l_k(y).
\end{equation}

Fix a $j\in\{1,\ldots, e\}$, then
\[\begin{split}
F(y) -\frac{1}{2}\sum_{k=1}^d \langle\sigma_j,u_k\rangle\log l_k(y) & =\frac{1}{2}\sum_{k=1}^d (-1-\langle\sigma_j,u_k\rangle)\log l_k(y)\\
                            & \geq C_j,
\end{split}\]
for some constant $C_j$, because each term $(-1-\langle\sigma_j,u_k\rangle)\log l_k(y)$ is bounded below.  Recall that
$\langle\sigma_j,u_k\rangle \geq -1,\ \forall k=1,\ldots,d$.

Taking the infimum $C$ of the $C_j,\ j=1,\ldots,e$, we get $F-\ol{w} \geq C$.

To prove the inequality $C'\geq  F-\ol{w}$ we define subsets of $\Sigma$.
Define $V_i =\{y\ |\ l_i(y) \leq\epsilon\}\cap\Sigma$, where $\epsilon>0$ is chosen small enough that the polytope
$\cap_{k=1}^d \{y\ |\ l_i(y) \geq\epsilon\}\subset\Sigma$ has the same faces as $\Sigma$.  Recall that a face
of $\Sigma$ is given by a multi-index $i_1,\ldots,i_\ell$ with $l_{i_1} =\cdots =l_{i_\ell} =0$.
For each face define $V_{i_1\cdots i_\ell}=\cap_{k=1}^\ell V_{i_k}$.

For each face of $\Sigma$ we define a subset $W_{i_1,\ldots,i_\ell}$ as follows.
\begin{eqnarray*}
  && W_0 =\Sigma -\bigcup_{j=1}^d V_j\\
  && W_i =V_i - V_i \cap\Bigl(\bigcup_{j\neq i} V_j \Bigr) \\
  && W_{i_1 i_2} =V_{i_1 i_2} - V_{i_1 i_2}\cap\Bigl(\bigcup_{j\neq i_1,i_2} V_j \Bigr) \\
  && \cdots \\
  && W_{i_1 i_2 \cdots i_\ell} =V_{i_1 i_2 \cdots i_\ell} - V_{i_1 i_2 \cdots i_k}\cap\Bigl(\bigcup_{j\neq i_1,i_2,\ldots,i_\ell} V_j \Bigr) \\
  && \cdots \\
  && W_{i_1 \cdots i_n} =V_{i_1 \cdots i_n}
\end{eqnarray*}

\[\begin{split}
F(y)-\ol{w}(y) & =\frac{1}{2}\sum_{k=1}^d -\log l_k(y) -underset{j=1,\ldots,e}{\sup}\Bigl[\frac{1}{2}\sum_{k=1}^d \langle\sigma_j,u_k\rangle\log l_k(y)\Bigr] \\
    & \leq C_0,
\end{split}\]
on $W_0$ for some $C_0$, because it is continuous and $\ol{W}_0$ is compact.

For $W_i$ choose a $\sigma_j$ with $\sigma_j(u_i) =-1$.  Then
\[\begin{split}
F(y)-\ol{w}(y) & \leq \frac{1}{2}\sum_{k=1}^d -\log l_k(y) -\frac{1}{2}\sum_{k=1}^d \langle\sigma_j,u_k\rangle\log l_k(y) \\
               & =\frac{1}{2}\sum_{\substack{k=1 \\k\neq i}} ^d -\log l_k(y)-\frac{1}{2}\sum_{\substack{k=1\\k\neq i}}^d \langle\sigma_j,u_k\rangle\log l_k(y) \\
               & \leq C_i,
\end{split}\]
for some constant $C_i$, because the remaining terms are continuous on $\ol{W}_i$.  In general, for the set
$W_{i_1 \cdots i_\ell}$ choose $\sigma_j$ so that $\sigma_j(u_{i_k})=-1,\ k=1,\ldots,\ell$.  Then as before
\[\begin{split}
F(y)-\ol{w}(y) \leq\frac{1}{2}\sum_{\substack{k=1 \\k\neq i_1,\ldots,i_\ell}} ^d -\log l_k(y)
-\frac{1}{2}\sum_{\substack{k=1\\k\neq i_1,\ldots,i_\ell}}^d \langle\sigma_j,u_k\rangle\log l_k(y) \\
               & \leq C_{i_1 \cdots i_\ell},
\end{split}\]
where the constant $C_{i_1 \cdots i_\ell}$ exists because all the expression is continuous on $\ol{W}_{i_1 \cdots i_\ell}$.

Letting $C'$ be the supremum of the $C_{i_1 \cdots i_\ell}$, we have $C\leq F-\ol{w}\leq C'$.
\end{proof}

Given a $G$-invariant $\phi\in C^\infty(X)$, we will denote its descent to a $\mathcal{W}_0$-invariant smooth function
on $N_\R$ by $\tilde{\phi}$.  We define
\[ P_G(N_\R,F) =\Bigl\{\tilde{\phi}\in C^2(N_\R)^{\mathcal{W}_0}\ |\ \frac{\partial^2(F+\tilde{\phi})}{\partial x_i\partial x_j} \geq 0,\
\underset{N_\R}{\sup}{\, \tilde{\phi}} =0,\ \text{and } |\tilde{\phi}|\ \text{is bounded on}\ N_\R \Bigr\}.\]
The following proposition was proved in~\cite{BatSel99}.
\begin{prop}\label{prop:alpha-tilde}
Let $X$ be a toric Fano orbifold with $G\subset\Aut(X)$ as above.  Let $dx$ be the volume form on $N_\R \cong\R^n$ corresponding
to the Haar measure normalized by the lattice $N\subset N_\R$.  Let $\tilde{\alpha}_G(X)$ be the supremum of all
$\alpha>0$ such that
\[ \int_{N_\R} e^{-\alpha\tilde{\phi}-F}\, dx \leq\tilde{C}(\alpha),\quad\forall\tilde{\phi}\in P_G(N_\R,F).\]
Then
\[ \tilde{\alpha}_G(X)\leq\alpha_G(X).\]
\end{prop}
The proof in~\cite{BatSel99} works here, so we omit it.  It follows easily from the following observation.
 As in the smooth case, we have that
\[ e^{-F} \frac{dt_1 \wedge d\ol{t}_1\wedge\ldots\wedge dt_n \wedge d\ol{t}_n}{|t_1|^2 \cdots|t_n|^n}
= e^{-F} dx_1\wedge\ldots\wedge dx_n\wedge d\theta_1\wedge\ldots\wedge\theta_n \]
can be extended to a non-vanishing volume form on $X$.  Therefore it is related to the volume form of $g_0$ by
\[ e^h d\mu_{g_0} = e^{-F} dx_1\wedge\ldots\wedge dx_n\wedge d\theta_1\wedge\ldots\wedge\theta_n \]
for $h\in C^\infty(X)$, where $h$ differs from $f$ defined after (\ref{eq:Monge-Amp}) by a constant.

\begin{lem}\label{lem:int-bound}
Let $\lambda >0$.  Then $\int_{N_\R} e^{-\lambda F}\, dx \leq C(\lambda)$.
\end{lem}
\begin{proof}
By Lemma~\ref{lem:bound} we have
\begin{equation}\label{eq:int-bd}
 \int e^{-\lambda F}\, dx \leq\int e^{-\lambda C -\lambda\ol{w}}\, dx =e^{-\lambda C}\int e^{-\lambda\ol{w}}\, dx.
\end{equation}
Let $\tau\in\Delta^*(n)$ be spanned by $u_{i_1},\ldots, u_{i_n}$.  Then restricted to the cone
$-\tau=\R_{\geq 0}\{-u_{i_1},\ldots, -u_{i_n}\}$
we have $\ol{w}=\sigma$, where $\sigma$ is the linear function with $\sigma(-u_{i_k}) =1$.
Therefore
\[ \begin{split}
e^{-\lambda C}\int_{-\tau} e^{-\lambda\ol{w}}\, dx & =
\frac{e^{-\lambda C}}{|\Gamma_\tau|}\int_{\R_{\geq 0}^n} e^{-\lambda(x_1 +\cdots+ x_n)}\, dx_1 \cdots dx_n \\
    & =\frac{e^{-\lambda C}}{|\Gamma_\tau|} \prod_{i=1}^n\Bigl( \int_{\R_{\geq 0}} e^{-\lambda x_i}\, dx_i \Bigr)\\
    & = \frac{e^{-\lambda C}}{|\Gamma_\tau|} \frac{1}{\lambda^n},
\end{split}\]
where $|\Gamma_\tau|$ is the order of the orbifold group $\Gamma_\tau$ associated to $\tau$, and combining with (\ref{eq:int-bd})
completes the proof since $N_\R =\cup_{\tau\in\Delta} -\tau$.
\end{proof}
\begin{lem}\label{lem:lower-bound}
There exists a constant $C$ so that for any $\tilde{\phi}\in P_G(N_\R,F)$ we have
\[ F(x) +\tilde{\phi} \geq C,\quad\forall x\in N_\R.\]
\end{lem}
\begin{proof}
Given an arbitrary $\tilde{\phi}\in P_G(N_\R,F)$ we consider the moment map
\[ \nu_{F+\tilde{\phi}} :N_\R \rightarrow M_\R,\]
\[ \nu_{F+\tilde{\phi}}(x) :=\Bigl(\frac{\partial(F+\tilde{\phi})}{\partial x_1}(x),\cdots,\frac{\partial(F+\tilde{\phi})}{\partial x_n}(x)  \Bigr).\]

We will first show that $\nu_{F+\tilde{\phi}}(N_\R)\subset\Sigma$.  Let $y_0 =\nu_{F+\tilde{\phi}}(x_0)$. By the convexity
of $F+\tilde{\phi}$,
\[ F(x) +\tilde{\phi}(x)\geq\langle y_0,x-x_0 \rangle +F(x_0)+\tilde{\phi}(x_0).\]
Thus $F(x) +\tilde{\phi}(x) -\langle y_0,x\rangle$ has a global minimum at $x_0$.  By Lemma~\ref{lem:bound}
and the fact that $\tilde{\phi}$ is globally bounded, $\ol{w}-\langle y_0,x\rangle \geq c$ for some constant $c$.
Since this is a piecewise linear function, we have
\begin{equation*}\label{eq:pl-bound}
\ol{w}-\langle y_0,x\rangle\geq 0,
\end{equation*}
and this implies that $y_0 \in\Sigma$.

Since $\underset{N_\R}{\sup}\tilde{\phi}=0$, we choose a sequence $\{p_k\}$ in $N_\R$ so that
$-1/k \leq\tilde{\phi}(p_k)\leq 0$.  Set $q_k =\nu_{F+\tilde{\phi}}(p_k)$.  Since $\Sigma$ is compact by passing to a subsequence
if necessary, we may assume that
\[ \underset{k}{\lim}\, q_k =q\in\Sigma.\]
The convexity of $F+\tilde{\phi}$ implies that
\[ F(x)+\tilde{\phi}(x) -\langle q_k,x\rangle \geq F(p_k)+\tilde{\phi}(p_k) -\langle q_k,p_k\rangle. \]
By Lemma~\ref{lem:bound} there is a constant $C$ so that
\[ F(p_k) +C \geq\ol{w}(p_k) \geq \langle q_k,p_k\rangle,\]
where the second inequality holds because $q_k \in\Sigma$.  Therefore
\[ F(x)+\tilde{\phi}(x) -\langle q_k,x\rangle \geq -C-\frac{1}{k}, \]
and taking $k\rightarrow\infty$
\begin{equation}\label{eq:lower-bd}
F(x)+\tilde{\phi}(x) -\langle q,x\rangle \geq -C.
\end{equation}
Since $\mathcal{W}_0$ is a finite group and $F,\ \tilde{\phi}$ are $\mathcal{W}_0$-invariant, one can average (\ref{eq:lower-bd})
to get
\begin{equation}
F(x)+\tilde{\phi}(x) -\langle\ol{q},x\rangle \geq -C.
\end{equation}
Here $\ol{q} =\frac{1}{|\mathcal{W}_0|}\sum_{g\in\mathcal{W}_0} g^* q$ is $\mathcal{W}_0$-invariant, and therefore $\ol{q}=0$.
\end{proof}

We can now prove the main theorem of the section.
\begin{thm}\label{thm:K-E}
Let $X_{\Sigma^*}$ be a symmetric toric Fano orbifold with $G\subset\Aut(X)$ as above, then $\alpha_G(X)\geq 1$.
Therefore, $X$ admits a $G$-invariant K\"{a}hler-Einstein metric.
\end{thm}
\begin{proof}
Let $0<\alpha<1$ and $\tilde{\phi}\in P_G(N_\R,F)$, then
\[\begin{split}
\int_{N_\R} e^{-\alpha\tilde{\phi}-F}\, dx & =\int_{N_\R} e^{-\alpha(\tilde{\phi}+F)}\, e^{(\alpha-1)F}\, dx \\
            & \leq e^{-\alpha C}\int_{N_\R} e^{(\alpha-1)F}\, dx \quad\text{(Lemma~\ref{lem:lower-bound})} \\
            & \leq e^{-\alpha C} C(1-\alpha), \quad\text{(Lemma~\ref{lem:int-bound})}.
\end{split}\]
Thus $\tilde{\alpha}_G \geq 1$, and the theorem follows from Proposition~\ref{prop:alpha-tilde}.
\end{proof}

\section{Corresponding Sasaki-Einstein spaces and embeddings}\label{sec:Sasaki-sub}

In this section we prove the correspondence in (\ref{diag:intro}).  First we obtain the toric surface $X$ and
Sasaki-Einstein space $M$ from $\mathcal{S}$ only using toric geometry.  It is an elementary result of the toric
geometry of a toric $\ASD$ Einstein space $\mathcal{M}$ that there is a toric Fano orbifold surface $X_{\delta^*}$
associated to it.  The Sasaki-Einstein space $M$ is not necessarily smooth.  In the following section we prove the embeddings
in (\ref{diag:intro}) from which it follows that $M$ is smooth precisely when the 3-Sasaki space $\mathcal{S}$ associated
to $\mathcal{M}$ is.

\subsection{Toric surfaces and $\ASD$ Einstein orbifolds}\label{subsec:toric-sur-asd}

We will consider toric anti-self-dual Einstein orbifolds $\mathcal{M}$ in greater detail.
By the previous Section~\ref{subsec:tor-3-Sasaki} quaternion-K\"{a}hler reduction gives us
infinitely many examples.  By reducing $\qps^{k+1}$ by a subtorus $T^k \subset Sp(k+2)$ defined by
an admissible matrix $\Omega$ we get a toric $\ASD$ Einstein orbifold $\mathcal{M}_\Omega$
with $b_2(\mathcal{M})=k$.  The orbifold $\mathcal{M}$ is characterized by a polygon
$Q_\Omega=\mathcal{M}/{T^2}$ with $k+2$ edges labeled in cyclic order with
$(m_0,n_0),(m_1,n_1),\ldots, (m_{k+2},n_{k+2})$ in $\Z^2$ with $(m_0,n_0)=-(m_{k+2},n_{k+2})$.
These vectors satisfy the following:
\renewcommand{\theenumi}{\alph{enumi}}
\begin{enumerate}
\renewcommand{\labelenumi}{(\theenumi)}
\item  The sequence $m_i, \ i=0,\ldots k+2$ is strictly increasing.

\item  The sequence $(n_i - n_{i-1})/(m_i -m_{i-1}), \ i=1,\ldots k+2$ is strictly increasing.
\end{enumerate}
We will make use of the following classification result of D. Calderbank and M. Singer~\cite{CalSin06}.
\begin{thm}\label{thm:asd-class}
Let $\mathcal{M}$ be a compact toric 4-orbifold with $\pi^{orb}_1(\mathcal{M})=e$ and $k=b_2(\mathcal{M})$.
Then the following are equivalent.
\begin{thmlist}
\item  One can arrange that the isotropy data of $\mathcal{M}$ satisfy (a) and (b) above by
cyclic permutations, changing signs, and acting by $Gl(2,\Z)$.
\item  $\mathcal{M}$ admits a toric $\ASD$ Einstein metric unique up to homothety and
equivariant diffeomorphism.  Furthermore, $(\mathcal{M},g)$ is isometric to the quaternionic K\"{a}hler
reduction of $\qps^{k+1}$ by a torus $T^k\subset Sp(k+2)$.
\end{thmlist}
\end{thm}
It is well known that the only possible smooth compact $\ASD$ Einstein spaces with positive scalar curvature
are $S^4$ and $\overline{\cps}^2$ with the round and Fubini-Study metrics~\cite{Hit81,FriKur82}.
Note that the stabilizer vectors $v_0 =(m_0,n_0),v_1 =(m_1,n_1),\ldots ,v_{k+2}=(m_{k+2},n_{k+2})$
form half a convex polygon with edges of increasing slope.

\begin{thm}
There is a one to one correspondence between compact toric anti-self-dual Einstein orbifolds $\mathcal{M}$ with
$\pi^{orb}_1(\mathcal{M})=e$ and
special symmetric toric Fano orbifold surfaces $X$ with $\pi^{orb}_1(X)=e$.
By theorem (\ref{thm:K-E}) $X$ has a K\"{a}hler-Einstein metric of positive scalar curvature.
Under the correspondence if $b_2(\mathcal{M})=k$, then $b_2(X)=2k+2$.
\end{thm}
\begin{proof}
Suppose $\mathcal{M}$ has isotropy data $v_0,v_1,\ldots,v_{k+2}$.  Then it is immediate that
$v_0,v_1,\ldots,v_{k+2},-v_1,-v_2,\ldots,-v_{k+1}$ are the vertices of a convex polygon
in $N_\R =\R^2$, which defines an augmented fan $\Delta^*$ defining $X$.  The symmetry of $X$ is clear.

Suppose $X$ is a special symmetric toric Fano surface.  Then $X$ is characterized by a convex polygon
$\Delta^*$ with vertices $v_0,v_1,\ldots, v_{2k+4}$ with $v_{2k+4}=v_0$.  Choose a primitive
$p=(u,w)\in\Z\times\Z, w>0$ which is not proportional to any $v_i -v_{i-1}, i=1,\ldots, k+2$.
Choose $s,t\in\Z$ with $su+tw=1$.  Then let $v^\prime_i ,i=0,\ldots,2k+4$ be the images of the $v_i$ under
$\begin{bmatrix}
w & -u \\
s & t
\end{bmatrix}$
There is a $v^{\prime}_j =(m^\prime_j,n^\prime_j)$ with $m^\prime_j$ smallest.
And $v^\prime_j, v^\prime_{j+1},\ldots,v^\prime_{j+k+2}$, where the subscripts are mod $2k+4$, satisfy
a.\ and b.
Such a toric orbifold is simply connected if and only if the isotropy data span $\Z\times\Z$.
One can show that the correspondence does not depend on the particular isotropy data.
\end{proof}

In the next section we will prove a more useful geometric correspondence between toric
$\ASD$ Einstein orbifolds and symmetric toric K\"{a}hler-Einstein surfaces.

\begin{xpl}
Consider the admissible weight matrix
\[\Omega=\left\lgroup
\begin{matrix}
1 & 0 & 1 & 1 \\
0 & 1 & 1 & 2
\end{matrix}\right\rgroup.\]
Then the 3-Sasakian space $\mathcal{S}_\Omega$ is smooth and
$b_2(\mathcal{S}_\Omega)=b_2(\mathcal{M}_\Omega)=2$.
And the anti-self-dual orbifold $\mathcal{M}_\Omega$ has isotropy data
\[v_0 =(-7,-2),(-5,-2),(-1,-1),(5,1),(7,2)=v_4.\]
The singular set of $\mathcal{M}$ consists of two points with stabilizer group $\Z_3$ and two with
$\Z_4$.  The associated toric K\"{a}hler-Einstein surface is that in figure (\ref{fig:typ-xpl}).
\end{xpl}

\begin{prop}\label{prop:ind-fano}
Let $X$ be the symmetric toric Fano surface associated to the $\ASD$ Einstein orbifold $\mathcal{M}$.
Then $\ind(X)=2$ if and only if $w_2(\mathcal{M})=0$.
In other words, $\mathbf{K}^{-1}_X$ has a square root if and only if the contact line bundle on $\mathcal{Z}$, $\mathbf{L}$, does.
\end{prop}
Recall that $w_2(\mathcal{M})$ is equal to the Marchiafava-Romani class $\varepsilon$.
Thus the vanishing of $w_2(\mathcal{M})$ is equivalent to the existence of a square root
$\mathbf{L}^{\frac{1}{2}}$ of the contact line bundle $\mathbf{L}$ on $\mathcal{Z}$.

\begin{proof}
Suppose $\ind(X)=2$ which is equivalent to $w_2(X)=0$, where $w_2$ denotes the orbifold Seifel-Whitney class.
Recall that the orbit space of $\mathcal{M}$ is a $k+2$-gon $W$ with labeled edges $C_1,\ldots,C_{k+2}$.
Since $\pi^{orb}_1(\mathcal{M})=e$, there exists an edge $C_i$ for which the orbifold uniformizing group $\Gamma$ has
odd order.  Let $U$ be a tubular neighborhood of an orbit in $C_i$.  So $U\cong S^1\times I\times D/\Gamma$, where
$I$ is an open interval and $D$ is a 2-disk.  And let $V$ be a neighborhood homotopically equivalent to
$\mathcal{M}\setminus U$ with $U\cup V=\mathcal{M}$.  Consider the exact homology sequence in
$\Z_2$-coefficients,
\begin{multline}\label{eq:asd-m-v}
\cdots\rightarrow H_2(BU)\oplus H_2(BV)\rightarrow H_2(B\mathcal{M})\rightarrow H_1(B(U\cap V))\\
\rightarrow H_1(BU)\oplus H_1(BV)\rightarrow 0.
\end{multline}
We have $BU\cong S^1\times I\times EO(4)/\Gamma$.  Since $EO(4)$ is contractible,
$H_*(EO(4)/\Gamma,A)=H_*(\Gamma,A)$ for any abelian group $A$.  In particular,
$H^n(\Gamma,\Z_2)=0$ for all $n>0$, since $|\Gamma|$ is odd.
Thus $H_2(BU,\Z_2)=0$ and $H_1(BU,\Z_2)=\Z_2$.  Similarly, it not hard to show that
$H_1(B(U\cap V),\Z_2)=\Z_2$.  From the exact sequence (\ref{eq:asd-m-v}) the inclusion $j:V\rightarrow\mathcal{M}$
induces a surjection $j_* :H_2(BV,\Z_2)\rightarrow H_2(B\mathcal{M},\Z_2)$.
Considering the orbit spaces one sees that there is a smooth embedding $\iota:V\rightarrow X$.
The tangent V-bundle $T\mathcal{M}$ lifts to a genuine vector bundle on $B\mathcal{M}$ which will also be denoted $T$.
Then
\[w_2(\mathcal{M})=w_2(T\mathcal{M})\in H^2(B\mathcal{M},\Z_2)=\Hom(H_2(B\mathcal{M},\Z_2),\Z_2).\]
Let $\alpha\in H_2(B\mathcal{M},\Z_2)$.  Then there exists a $\beta\in H_2(BV,\Z_2)$ with
$j_*\beta=\alpha$.  Then
\[ w_2(T\mathcal{M})(\alpha)=w_2(TV)(\beta)=w_2(TX)(\iota_*\beta)=0.\]
Thus $w_2(\mathcal{M})=0$.

The converse statement will follow from the main result of the next section.
\end{proof}

\subsection{Twistor space and divisors}\label{subsec:twist-div}

We will consider the twistor space $\mathcal{Z}$ of an $\ASD$ positive scalar curvature Einstein
orbifold $\mathcal{M}$.
For now suppose $(\mathcal{M},[g])$ is an anti-self-dual, i.e. $W_g^+ \equiv 0$,
conformal orbifold.  There exists a complex three dimensional orbifold
$\mathcal{Z}$ with the following properties:
\renewcommand{\theenumi}{\alph{enumi}}
\begin{enumerate}
\renewcommand{\labelenumi}{(\theenumi)}
\item  There is a $C^\infty$ orbifold bundle $\varpi:\mathcal{Z}\rightarrow\mathcal{M}$.
\item  The general fiber of $P_x =\varpi^{-1}(x),\ x\in\mathcal{M}$ is a projective line
$\cps^1$ with normal bundle $N\cong\mathcal{O}(1)\oplus\mathcal{O}(1)$, which holds over singular fibers
with $N$ an orbifold bundle.
\item  There exists an anti-holomorphic involution $\varsigma$ of $\mathcal{Z}$ leaving the fibers
$P_x$ invariant.
\end{enumerate}

Let $T$ be an oriented real 4-dimensional vector space with inner product $g$.  Let
$C(T)$ be set of orthogonal complex structures inducing the orientation, i.e.
if $r,s\in T$ is a complex basis then $r,Jr,s,Js$ defines the orientation.  One has
$C(T)= S^2\subset\Lambda^2_+(T)$, where $S^2$ is the sphere of radius $\sqrt{2}$.
Now take $T$ to be $\Ha$.
Recall that $Sp(1)$ is the group of unit quaternions.  Let
\begin{equation}
Sp(1)_+ \times Sp(1)_-
\end{equation}
act on $\Ha$ by
\begin{equation}
w\rightarrow gw{g^\prime}^{-1},\text{  for  }w\in\Ha\text{ and } (g,g^\prime)\in Sp(1)_+ \times Sp(1)_-.
\end{equation}
Then we have
\begin{equation}
Sp(1)_+ \times_{\Z_2} Sp(1)_- \cong SO(4),
\end{equation}
where $\Z_2$ is generated by $(-1,-1)$.
Let
\begin{equation}
\begin{split}
C & =\{ai+bj+ck: a^2+b^2+c^2=1, a,b,c\in\R\}\\
  & =\{g\in Sp(1)_+: g^2=-1\}\cong S^2.
\end{split}
\end{equation}
Then $g\in C$ defines an orthogonal complex structure by
\[ w\rightarrow gw,\text{  for  } w\in\Ha,\]
giving an identification $C=C(\Ha)$.
Let $V_+ =\Ha$ considered as a representation of $Sp(1)_+$ and a right $\C$-vector space.
Define $\pi:V_+\setminus\{0\}\rightarrow C$ by $\pi(h)=-hih^{-1}$.  Then the fiber of
$\pi$ over $hih^{-1}$ is $h\C$.  Then $\pi$ is equivariant if $Sp(1)_+$ acts on $C$ by
$q\rightarrow gqg^{-1}, g\in Sp(1)_+$.  We have a the identification
\begin{equation}
C=V_+\setminus\{0\}/{\C^*} =\mathbb{P}(V_+).
\end{equation}

Fix a Riemannian metric $g$ in $[g]$.  Let $\phi:\tilde{U}\rightarrow U\subset\mathcal{M}$
be a local uniformizing chart with group $\Gamma$.
Let $F_{\tilde{U}}$ be the bundle of oriented orthonormal frames on $\tilde{U}$.
Then
\begin{equation}\label{eq:twistor-chart}
F_{\tilde{U}}\times_{SO(4)}\mathbb{P}(V_+) =F_{\tilde{U}}\times_{SO(4)} C
\end{equation}
defines a local uniformizing chart for $\mathcal{Z}$ mapping to
\[F_{\tilde{U}}\times_{SO(4)}\mathbb{P}(V_+)/\Gamma =F_{\tilde{U}}/\Gamma \times_{SO(4)}\mathbb{P}(V_+).\]
Right multiplication by $j$ on $V_+ =\Ha$
defines the anti-holomorphic involution $\varsigma$ which is fixed point free on (\ref{eq:twistor-chart}).
We will denote a neighborhood as in (\ref{eq:twistor-chart}) by $\tilde{U}_{\mathcal{Z}}$.

An almost complex structure is defined as follows.  At a point $z\in\tilde{U}_{\mathcal{Z}}$ the Levi-Civita
connection defines a horizontal subspace $H_z$ of the real tangent space $T_z$ and we have a
splitting
\begin{equation}\label{eq:twistor-comp}
T_z =H_z\oplus T_z P_x = T_x\oplus T_z P_x,
\end{equation}
where $\varpi(z)=x$ and $T_x$ is the real tangent space of $\tilde{U}$.
Let $J_z$ be the complex structure on $T_x$ given by $z\in P_x=C(T_x)$, and let
$J_z^\prime$ be complex structure on $T_x\oplus T_z P_x$ arising from the natural complex structure
on $P_x$.  Then the almost complex structure on $T_z$ is the direct sum of $J_z$ and $J_z^\prime$.
This defines a natural almost complex structure on $Z_{\tilde{U}}$ which is invariant under $\Gamma$.
We get an almost complex structure on $\mathcal{Z}$ which is integrable precisely when $W_+\equiv 0$.

Assume that $\mathcal{M}$ is $\ASD$ Einstein with non-zero scalar curvature.  Then
$\mathcal{Z}$ has a complex contact structure $D\subset T^{1,0}\mathcal{Z}$ with
holomorphic contact form \mbox{$\theta\in\Gamma(\Lambda^{1,0}\mathcal{Z}\otimes\mathbf{L})$} where
$\mathbf{L}=T^{1,0}\mathcal{Z}/D$.

The group of isometries $\Isom(\mathcal{M})$ lifts to an action on $\mathcal{Z}$ by real holomorphic transformations.
Real means commuting with $\varsigma$.  This extends to a holomorphic action of the complexification
$\Isom(\mathcal{M})_{\C}$.
For $X\in\mathfrak{Isom}(\mathcal{M})\otimes\C$, the Lie algebra of $\Isom(\mathcal{M})_{\C}$,
we will also denote by $X$
the holomorphic vector field induced on $\mathcal{Z}$.
Then $\theta(X)\in H^0(\mathcal{Z},\mathcal{O}(\mathbf{L}))$.
By a well known twistor correspondence the map $X\rightarrow\theta(X)$ defines an isomorphism
\begin{equation}\label{eq:twistor-trans}
\mathfrak{Isom}(\mathcal{M})\otimes\C\cong H^0(\mathcal{Z},\mathcal{O}(\mathbf{L})),
\end{equation}
which maps real vector fields to real sections of $\mathbf{L}$.

Suppose for now on that $\mathcal{M}$ is a toric $\ASD$ Einstein orbifold with twistor space $\mathcal{Z}$.
We will assume that $\pi_1^{orb}(\mathcal{M})=e$ which can always be arranged by taking the orbifold cover.
Then as above $T^2$ acts on $\mathcal{Z}$ by holomorphic transformations.  And the action extends to
$T^2_\C =\C^*\times\C^*$, which in this case is an algebraic action.
Let $\mathfrak{t}$ be the Lie algebra of $T^2$ with $\mathfrak{t}_\C$ the Lie algebra of $T^2_\C$.  Then we have
from (\ref{eq:twistor-trans}) the pencil
\begin{equation}
P=\mathbb{P}(\mathfrak{t}_\C)\subseteq |L|,
\end{equation}
where for $t\in P$ we denote $X_t=(\theta(t))$ the divisor of the section $\theta(t)\in H^0(\mathcal{Z},\mathcal{O}(\mathbf{L}))$.
Note that $P$ has an equator of real divisors.  Also, since $T^2_\C$ is abelian, every $X_t, t\in P$ is
$T^2_\C$ invariant.

Consider again the $T^2$-action on $\mathcal{M}$.
Let $K_x$ denote the stabilizer of $x\in\mathcal{M}$.
Recall the set with non-trivial stabilizers of the $T^2$-action on $\mathcal{M}$ is
$B=\bigcup_{i=1}^{k+2}B_i$ where $B_i$ is topologically a 2-sphere.  Denote
$x_i=B_i\cap B_{i+1}$, $B_i^\prime =B_i\setminus\{x_i,x_{i-1}\}$ and $B^\prime =\bigcup_{i=1}^{k+2}B_i^\prime$.
And denote the stabilizer of $B_i^\prime =B_i\setminus\{x_i,x_{i-1}\}$ by
$K_i =S^1(m_i,n_i)$.  The stabilizer of $x_i$ is $K=T^2$.
We will first determine the singular set $\Sigma\subset\mathcal{Z}$ for the $T^2$-action on $\mathcal{Z}$.
\begin{lem}\label{lem:twistor-fix}
For $x\in B$ there exists on $P_x$ precisely two fixed points $z^+,z^-$ for the action of $K_x$ which
are $\varsigma$ conjugate.  For $x\in B'$, the stabilizer group in $T^2$ of any other $z\in P_x$
is trivial.
\end{lem}
\begin{proof}
Let $\phi:\tilde{U}\rightarrow U$ be a uniformizing chart centered at $x$ with group $\gamma$.
We may assume that  $\tilde{K}_x$ acts on $\tilde{U}$ with $\gamma\subset\tilde{K}_x$ and $\tilde{K}_x/\gamma =K_x$.
Then the uniformized tangent space splits
\begin{equation}\label{eq:twist-iso}
T_{\tilde{x}}=T_1\oplus T_2.
\end{equation}
When $x\in B^\prime$ we take $T_1$ to be the space on which $\tilde{K}_i$ acts trivially and $T_2$ on which
$\tilde{K}_i$ act faithfully.  When $x=x_i$, $\tilde{K}_x =\tilde{K}_i\times\tilde{K}_{i+1}$ assume
$\tilde{K}_i$ acts faithfully on $T_1$ and trivially on $T_2$, and $\tilde{K}_{i+1}$ trivially on $T_1$ and
faithfully on $T_2$.

We determine the action of $\tilde{K}_x$ on $z\in \tilde{P}_x$.
Identify (\ref{eq:twist-iso}) with $\Ha=\C\oplus j\C$, considered as a right $\C$-vector space.
The action of $\tilde{K}_x =S^1(t)$ in the first case is
\[ (x,y)\rightarrow (x, ty),\]
and the action of $\tilde{K}_x =S^1(s)\times S^1(t)$ in the second is
\[ (x,y)\rightarrow (sx, ty).\]
If $(u,v)\in S^1\times S^1\subset Sp(1)_+\times Sp(1)_-$, then the action of $(u,v)$ on $T_x$ is
\[ (x,y)\rightarrow (uv^{-1}x, (uv)^{-1}y).\]
In the first case the action of $\tilde{K}_i$ is realized by the subgroup $\{(u,u)\}$ with
$t=u^{-1}$.  Considering the representation of $Sp(1)_+$ on $V_+=\Ha$, $u$ acts by
\[(w,z)\rightarrow (uw,u^{-1}z).\]
One sees that the only fixed points on $\tilde{P}_x=\mathbb{P}(V_+)$ are
$[1:0]$ and $[0:1]$.  It is easy to see that $\tilde{K}_i$ acts freely on every other point of $\tilde{P}_x$.
This also proves the statement for $x=x_i$.
\end{proof}

Denote the two $K_x$ fixed points on $P_x$ for $x=x_i$ by $z_i^{\pm}$.
We will denote $P_i :=P_{x_i}, i=1,\ldots, k+2$.
The next result is an easy consequence of the last lemma.
\begin{lem}
There exist two irreducible rational curves $C_i^{\pm}, i=1,\ldots, k+2$ mapped diffeomorphically to
$B_i$ by $\varpi$.  Furthermore, $\varsigma(C_i^{\pm})=C_i^{\mp}$.
\end{lem}

The singular set for the $T^2$-action on $\mathcal{Z}$ is the union of rational curves
\begin{equation}
\Sigma=\Bigl(\cup_{i=1}^{k+2}P_i \Bigr)\bigcup\Bigl(\cup_{i+1}^{k+2}C_i^+\cup C_i^-\Bigr).
\end{equation}
The fixed points for $T^2$ are $z_i^{\pm}, i=1,\ldots,k+2$.  And the stabilizer group of
$C^{\prime\pm}_i =C^\pm_i\setminus\{z_i^\pm ,z_{i-1}^\pm\}$ is $K_i$.
If $S_{\mathcal{Z}}$ is the orbifold singular set, then $S_{\mathcal{Z}}\subset\Sigma$.  In this case
$S_{\mathcal{Z}}=\Sing(\mathcal{Z})$, the singular set of $\mathcal{Z}$ as an analytic variety.

We will denote the union of the curves $C_i^\pm$ by
\[ C=\bigcup_{i=1}^{k+2}(C_i^+\cup C_i^-). \]
Then either $C$ is a connected cycle, or it consists of two $\varsigma$-conjugate cycles.
It will turn out that $C$ is always connected.  Thus it may be more convenient to denote its
components by $C_i ,i=1,\ldots, 2n$, where $n=k+2$, and the points $z_i^\pm$ by $z_i$ and $z_{i+n}$
such that
\[ z_i=C_i \cap C_{i+1}, i=1,\ldots, 2n, \]
where we take the index to be mod $2n$.

We now consider the action of $T^2_\C$ on $\mathcal{Z}$.
The stabilizer group of $z\in\mathcal{Z}$ in $T^2_\C$ will be denoted $G_z$.
Let $G_i\subset T^2_\C$ be the complexification
of $K_i$.
\begin{lem}
For $z\in C_i^\prime , i=1,\ldots, 2n$, the stabilizer group $G_z$ coincides with $G_i$.
\end{lem}
\begin{proof}
We have $G_i\subset G_z$ with $\dim G_i=1$.  Suppose $G_i\neq G_z$ then $G_z/G_i$ is a discrete subgroup
of $T^2_\C/{G_i}\cong\C^*$.  It is easy to see that $G_z/G_i$ is an infinite cyclic subgroup of $T^2_\C/G_i$.
Then the orbit of $z$, $C^\prime_i\cong T^2_\C/G_z$ must be a one dimensional complex torus, which is a
contradiction.
\end{proof}

Recall that a parametrization of a stabilizer group $K_i =S^1(m_i,n_i)$, of $B_i$, $i=1,\ldots, n$, is only fixed
up to sign.  This amounts to a choice of orientation of $B_i$.  In view of proposition (\ref{prop:twistor-iso})
for the stabilizer group $G_i$ of $C_i^\prime$, $i=1,\ldots, 2n$, there is a fixed parametrization
$\rho_i :\C^*\rightarrow T^2_\C$.  One picks one of two possibilities by the rule:
For $z$ in a sufficiently small neighborhood of a point of $C_i^\prime$ one has
\[ \underset{t\rightarrow 0}{\lim}\, \rho_i(t)z\in C_i.\]
\begin{lem}
We have $\rho_i =-\rho_{i+n}$ for $i=1,\ldots, n$, where we consider the $\rho_i$ to be elements of the $\Z^2$
lattice of one parameter subgroups of $T^2_\C$.
\end{lem}
\begin{proof}
Let $x\in B_i'$.  And consider the action of $G_i$ on the twistor line $P_i$ as described in
the proof of lemma (\ref{lem:twistor-fix}).
If $z\in P_i$, then $\lim_{t\rightarrow 0}\rho_i(z)=z_+ \in C_i$ implies
$\lim_{t\rightarrow 0}\rho_i^{-1}(z)=z_+ \in C_{i+n}$.
\end{proof}

We now consider the isotropy representations of $G_z$.  The proof of the following is straight forward.
\begin{prop}\label{prop:twistor-iso}
Let $z\in C$ with $\varpi(z)=x$.
And let $\phi:\tilde{U}\rightarrow U$ be a $\tilde{K}_x$-invariant local uniformizing chart with group $\gamma\subset\tilde{K}_x$.
Also $\tilde{G}_i$ denotes the complexification of $\tilde{K}_x$.
\begin{thmlist}
\item  Let $z\in C^\prime_i, i=1,\ldots,2n$.  Then there are $\C$-linear coordinates $(u,v,w)$ on
$T_{\tilde{z}}\tilde{U}_\mathcal{Z}$ and an identification $\tilde{G}_z\cong\C^*(t)$ so that
$\tilde{G}_z$ acts by
\[ (u,v,w)\rightarrow (u,tv,tw).\]
And the subspace $v=w=0$ maps to the tangent space of $C^\prime_i$ at $z$.
\item  Let $z=z_i$ for $i=1,\ldots, 2n$.  Then there are $\C$-linear coordinates $(u,v,w)$ on
$T_{\tilde{z}}\tilde{U}_\mathcal{Z}$ and an identification $\tilde{G}_z\cong\C^*(s)\times\C^*(t)$ so that
$\tilde{G}_z$ acts by
\[ (u,v,w)\rightarrow (stu,sv,tw).\]
And the uniformized tangent space of $P_i$ (resp. $C_i$, and $C_{i+1}$) at $z$ is the subspace
$v=w=0$ (resp. $u=v=0$ and $u=w=0$).
\end{thmlist}
\end{prop}

We will determine the $T^2_\C$-action in a neighborhood of $C$.
Let $z\in C$, and let $\tilde{U}_{\mathcal{Z}}$ be a $\tilde{K}_z$-invariant uniformizing
neighborhood as above with local group
$\gamma\subset\tilde{K}_z$.  Then there is
\renewcommand{\theenumi}{\roman{enumi}}
\begin{enumerate}
\renewcommand{\labelenumi}{\theenumi.}
\item  a $\tilde{K}_z$-invariant neighborhood $W$ of the origin in $T_{\tilde{z}}\tilde{U}_\mathcal{Z}$,

\item  a $\tilde{K}_z$-invariant neighborhood $V$ of $\tilde{z}$ in $\tilde{U}_{\mathcal{Z}}$, and

\item  a $\tilde{K}_z$-invariant biholomorphism $\varphi: W\rightarrow V$, i.e.
\begin{equation}\label{eq:twistor-equi}
\varphi(gx)=g\varphi(x), \text{  for }x\in W, g\in\tilde{K}_z.
\end{equation}
\end{enumerate}
This is a well known; see for example~\cite{BoM}.

This linear action extend locally to $\tilde{G}_z$, where $\tilde{G}_z$ is the complexification of
$\tilde{K}_z$.
Let $W_0\subset W$ be a connected relatively compact neighborhood of the origin.
And define the open set $A=\{(g,w)\in\tilde{G}_z\times W_0 : gw\in W \}$, and let
$A_0\subset A$ be the connected component containing $\tilde{K}_z \times W_0$.
Then for any $(g,w)\in A_0$, we have $g\varphi(w) \in V$ and (\ref{eq:twistor-equi}).

We now describe the local action of $T^2_\C$ around a point $z\in C$.
There are two cases, i. and ii., distinguished as in proposition (\ref{prop:twistor-iso}).
In case i. $z\in C_i^\prime$ for some $i=1,\ldots, 2n$.  And in case ii. $z=z_i$ for some
$i=1,\ldots, 2n$.  we will use proposition (\ref{prop:twistor-iso}) and the above remarks to produce a
neighborhood $U$ of $z$ as follows.

Case (i).  Suppose $z\in C_i^\prime$.  There exists an equivariant uniformizing neighborhood $\phi:\tilde{U}\rightarrow U$
centered at $z$ with group $\gamma\subset\tilde{K}_i$.  One can lift the corresponding one parameter group
$\tilde{\rho}_i:\C^*(t)\rightarrow\tilde{T}_\C$ with image $\tilde{G}_i$.  Let $\tilde{G}^\prime=\C^*(s)$ be a compliment
to $\tilde{G}_i$ in $\tilde{T}_\C$.  There exists coordinates $(u,v,w)$ in $\tilde{U}$ so that
\begin{equation}
\tilde{U}=\{(u,v,w):|u-1|<\epsilon, |v|<1, |w|<1\},\epsilon>0, \tilde{z}=(1,0,0).
\end{equation}
And $v=w=0$ is the subset mapped to $C$ and $\tilde{G}_i$ acts by
\begin{equation}
(u,v,w)\rightarrow (u,tv,tw), \text{  for }|t|\leq 1.
\end{equation}
The action of $\tilde{G}^\prime$ is given by $(u,v,w)\rightarrow (su,v,w)$ for
$|su-1|<\epsilon$.

Case (ii).  Suppose $z=z_i$ for some $i=1,\ldots, 2n$.  There exists an equivariant uniformizing neighborhood $\phi:\tilde{U}\rightarrow U$
centered at $z$ with group $\gamma\subset\tilde{K}_z =\tilde{T}^2$.  And one can lift the one parameter groups to
$\tilde{\rho}_i$ and $\tilde{\rho}_{i+1}$ to give an isomorphism
\[ \tilde{\rho}_i\times\tilde{\rho}_{i+1}: \C^*(s)\times\C^*(t)\rightarrow\tilde{T}^2_\C,\]
where $\tilde{T}^2_\C$ is the complexification of $\tilde{T}^2$.  There exists coordinates $(u,v,w)$ in
$\tilde{U}$ so that
\begin{equation}
\tilde{U}=\{(u,v,w):|u|<1, |v|<1, |w|<1\}, \tilde{z}=(0,0,0),
\end{equation}
where the equations $u=v=0$, $u=w=0$, and $v=w=0$ are the equations defining the subsets mapped to
$C_i$, $C_{i+1}$, and $P_i$ respectively.  And the action of $(s,t)\in\C^*(s)\times\C^*(t)$ is given by
\begin{equation}
(u,v,w)\rightarrow (stu,sv,tw), \text{  for }|s|\leq 1, |t|\leq 1.
\end{equation}

We will call such a neighborhood $U$ of a point of $C$ an \emph{admissible neighborhood}, and
$\phi:\tilde{U}\rightarrow U$ with group $\gamma$ an \emph{admissible uniformizing system}.
Let $U$ be an admissible neighborhood.  We set
\[ U^\prime :=U\setminus\Sigma.\]
Denote by $\tilde{U}'$ the preimage of $U^\prime$ in $\tilde{U}$.
We will define subsets $\tilde{U}_{ab}', \tilde{U}_{01}',$ and $\tilde{U}_{01}''$ of $\tilde{U}'$.

In case (i), for $(a,b)\neq 0$, define
\[ \tilde{U}_{ab}' :=\{(u,v,w)\in\tilde{U}' :av=bw \}.\]

In case (ii), for $(a,b)$ with $a\neq 0$, we define
\[ \tilde{U}_{ab}' :=\{(u,v,w)\in\tilde{U}' : au=bvw \},\]
and the two subsets
\[ \tilde{U}_{01}' :=\{(u,v,w)\in\tilde{U}' :v=0\}, \text{ and }\tilde{U}_{01}'' :=\{(u,v,w)\in\tilde{U}' :w=0\}.\]
\begin{lem}
The subsets defined above are connected closed submanifolds of $\tilde{U}'$ and each consists of a single
local $\tilde{T}^2_\C$-orbit with these being all the orbits.  And the closure of each orbit is an analytic submanifold of
$\tilde{U}$.
\end{lem}
This follows from the above description of the $\tilde{T}^2_\C$-action.
Note that $\gamma$ preserves the orbits so this gives a description of the local orbits of $T^2_\C$ in $U$.
We will denote by $U_{ab}', U_{01}'$, and $U_{01}''$ the corresponding local orbits in $U$.

We have the local leaf structure of the orbits in an admissible neighborhood.  In most cases this gives
the global leaf structure.
\begin{lem}\label{lem:twistor-leaf}
Let $U$ be an admissible neighborhood.  Let $E,F\subset U'$ be separate local leaves not both
being of type $U_{01}'$ or $U_{01}''$.  Then $E$ and $F$ are not contained in the same
$T^2_\C$-orbit.
\end{lem}
\begin{proof}
After acting by an element of $T^2_\C$ we may assume $U$ is an admissible neighborhood as in case
(i). with coordinates $(u,v,w)$ and $v=w=0$ defining $C_i\cap U$.  Let $z\in E$ and $z'\in F$ both have $u=1$.
There is a $g\in T^2_\C$ with $gz=z'$.  Let
$z_0 =\lim_{t\rightarrow 0}\rho_i(t)z =\lim_{t\rightarrow 0}\rho_i(t)z'$.  Then
\[ gz_0 = g\left(\lim_{t\rightarrow 0}\rho_i(t)z\right)= \lim_{t\rightarrow 0}\rho_i(t)gz
=\lim_{t\rightarrow 0}\rho_i(t)z' =z_0.\]
So $g\in G_i$, and $g=\rho_i(t_0)$.  If $|t_0|\leq 1$, then $g$ preserves the local leaves.
If $|t_0|>0$, the equation $z=g^{-1}z'$ gives a contradiction.
\end{proof}

\begin{lem}\label{lem:twistor-sta}
For any $z\in U'$, an admissible neighborhood, the stabilizer group $G_z$ is the identity.
\end{lem}
\begin{proof}
If $g\in G_z$, then $g$ fixes the entire $T^2_\C$-orbit of $z$.  Therefore $g$ fixes the entire
set $U_{ab}'$ containing $z$.  But the closure of $U_{ab}'$ intersects either $C_i$ or $C_{i+1}$.
So $g$ is contained in either $G_i$ or $G_{i+1}$.  But from the above description of the action on $U'$,
we see that $g=e$.
\end{proof}

\begin{lem}\label{lem:twistor-linest}
Let $z$ be any point of $P_i' =P_i\setminus\{z_i,z_{i+n}\}$.  And let $U$ be an admissible neighborhood of
$z_i$ or $z_{i+n}$.  Then there exists a neighborhood $V$ of $z$ and $g\in T^2_\C$ so that
$g(V)\subset U$.
\end{lem}
\begin{proof}
The stabilizer group of $P_i'$ is the image of the one parameter group
$\rho_i{\rho_{i+1}}^{-1}:C^*(s)\rightarrow T^2_\C$.  Then the orbit of $z$ by $G_i$ for example is
$P_i'$.  So a suitable element $g\in G_i$ will work.
\end{proof}

By lemmas (\ref{lem:twistor-sta}) and (\ref{lem:twistor-linest}) there is a small neighborhood $W$ of
$\Sigma\subset\mathcal{Z}$, so that if we set $W':=W\setminus\Sigma$, the stabilizer of every point
or $W'$ in $T^2_\C$ is the identity.

Our goal is to determine the structure of the divisors in the pencil $P$.
As before we will consider the one parameter groups $\rho_i \in N=\Z\times\Z$, where $N$ is the lattice
of one parameter $\C^*$-subgroups of $T^2_\C$.
Also, we will identify the Lie algebra $\mathfrak{t}$ of $T^2$ with $N\otimes\R$ and
the Lie algebra $\mathfrak{t}_\C$ of $T^2_\C$ with $N\otimes\C$.
Since $\mathbf{L}|_{P_x} =\mathcal{O}(2)$ a divisor $X_t\in P$ intersects a generic twistor line
$P_x$ at two points.
\begin{lem}
For any $X_t \in P$, we have $C \subset X_t$.
\end{lem}
\begin{proof}
Let $x\in B_i$.  Suppose that $z\in P_x'=P_x\setminus\{z^+,z^-\}$ and $z\in X_t$.
Then $G_i$ preserves the twistor line $P_x$, and the orbit of $z$ by $G_i$ is
$P_x'$.  Since $X_t$ is $T^2_\C$-invariant $P_x\subset X_t$.
Therefore, we either have $P_x\cap X_t=\{z^+,z^-\}$ or $P_x\subset X_t$.
\end{proof}

\begin{thm}\label{thm:twistor-div}
Let $\mathcal{M}$ be a compact $\ASD$ Einstein orbifold with $b_2(\mathcal{M})=k$ and
$\pi^{orb}_1(\mathcal{M})=e$.  Let $n=k+2$.
Then there are distinct real points $t_1,t_2,\ldots,t_{n}\in P$ so that for
$t\in P\setminus\{t_1,t_2,\ldots,t_{n}\}$, $X_t\subset\mathcal{Z}$ is a suborbifold.
And $X_t$ is a special symmetric toric Fano surface.  The anti-canonical cycle of $X_t$
is $C_1,C_2,\ldots, C_{2n}$, and the corresponding stabilizers are
$\rho_1, \rho_2, \ldots, \rho_{2n}$ which define the vertices in $N=\Z\times\Z$ of
$\Delta^*$ with $X_t =X_{\Delta^*}$.

For $t_i\in P$, $X_{t_i} =D+\ol{D}$, where $D,\ol{D}$ are irreducible degree one divisors with
$\varsigma(D)=\ol{D}$.  The $D,\ol{D}$ are suborbifolds of $\mathcal{Z}$ and are toric Fano
surfaces.  We have $D\cap\ol{D}=P_{i}$ and the elements
$\pm(\rho_1,\ldots,\rho_i,-\rho_i +\rho_{i+1},\rho_{n+i+1},\ldots, \rho_{2n})$
define the augmented fans for $D$ and $\ol{D}$.
\end{thm}
\begin{proof}
Let $z\in W'$, so the stabilizer of $z$ is the identity.
Let O be the $T^2_\C$-orbit of $z$.  Since $T^2_\C$ has only one end $\ol{O}\setminus O$ is connected.
Since $\ol{O}\cap\Sigma\neq\emptyset$, and the stabilizer of every point of $W'$ is the identity,
$\ol{O}\setminus O\subset\Sigma$.

Define elements $t_i\in P$ by $t_i =\rho_{i+1}-\rho_i, i=1,\ldots,n$.  Recall
that the stabilizer of $P_i$ is $\rho_i\rho^{-1}_{i+1}:\C^*\rightarrow T^2_\C$.
If $t\in P\setminus\{t_1,t_2,\ldots,t_{n}\}$, then a vector field induced by $t$ is tangent to,
and non-vanishing on,
$P_i, i=1,\ldots,n$.  Since the contact structure $D=\ker\theta$ is transverse to the twistor lines,
$P_i\cap X_t =\{z^+,z^-\}$.  Let $z\in X_t$ be in an admissible neighborhood of $C$.
Then the $T^2_\C$-orbit O of $z$ satisfies $\ol{O}\setminus O\subset C$.
The intersection of O with any admissible neighborhood is a leaf $U_{ab}'$ which has
analytic closure.  Let $Y=\ol{O}$, then $Y$ is an analytic subvariety.

Suppose $C$ consists of two disjoint cycles with $Y\cap C=\bigcup_{i=1}^n C_i$.
Then $Y$ is a degree one divisor, i.e. intersecting a generic twistor line at one point.
If $\ol{Y}=\varsigma(Y)$, then $Y\cap\ol{Y}=\bigcup_{i=1}^m P_{x_i}$, a disjoint union of twistor lines
with $x_i\notin B, i=1,\ldots, m$.  Since $Y\cap\ol{Y}$ is $T^2_\C$-invariant, we must have
$Y\cap\ol{Y}=\emptyset$.  Thus $Y$ intersects each twistor line at one point.
This is impossible. $Y$ defines a, positively oriented, almost complex structure $J$ on $\mathcal{M}$.
Then if $c_1 =c_1^{orb}(\mathcal{M},J)$, $c_1^2 =2\chi_{orb} +3\tau_{orb}$ where $\chi_{orb}$ and
$\tau_{orb}$ are defined by the same Gauss-Bonnet formulae as on smooth 4-manifolds~\cite{Bess87}.  We have
\[2\chi_{orb}+3\tau_{orb}=\frac{1}{4\pi^2}\int_\mathcal{M} \frac{s^2}{24}d\mu >0. \]
But a familiar Bochner argument shows the intersection form is negative definite.
Therefore $C\subset Y$, $Y$ is a degree two divisor, and $X_t =Y$.  From
the description of the admissible uniformizing systems and the local leaves, we see that
$X_t$ is a suborbifold.  Since $X_t$ is the closure of an orbit isomorphic to $T^2_\C$ it
is a toric variety and has the anti-canonical cycle $C$ and stabilizers $\rho_i$ defining
$\Delta^*$.

The adjunction for $X=X_t$ formula gives
$\mathbf{K}_X \cong\mathbf{K}_\mathcal{Z}\otimes[X]|_X =\mathbf{K}_\mathcal{Z}^{\frac{1}{2}}|_X$.
Thus $\mathbf{K}^{-1}_X >0$.  The orbifold version of the Kodaira embedding theorem~\cite{Bai57} implies
that $\mathbf{K}^{-m}_X$ is very ample for $m\gg 0$.  From basic properties of toric varieties it follows that
$-k\in\SF(\Delta^*)$ is strictly upper convex.  Thus $X$ is Fano and
$\Delta^*$ is a convex polytope.  It follows that $t_1,\ldots,t_n\subset P$ form a cycle of distinct
points.

Suppose $t=t_i, i=1,\ldots,n$.  Then $X_t\cap\Sigma =C\cup P_i$.  Let $z\in X_t$ be in an
admissible neighborhood of type i. with orbit O.  Let $D=\ol{O}$.  Then
$D\setminus O\subset C\cup P_i$.  And $P_i\subset D$, for otherwise we would have $D=X_t$ as in the last paragraph.
For an admissible neighborhood $U$ of $z_i$ or $z_{i+n}$
O must intersect $U$ in a leaf $U_{01}'$ or $U_{01}''$.  This can be seen from lemma
(\ref{lem:twistor-leaf}).  We must have either $D\cap\Sigma =C\cup P_i$ or
a cycle of the form $C_1,\ldots,C_i,P_i,C_{i+n+1},\ldots, C_{2n}$.
In the first case $D=X_t$ is irreducible.   Since $X_t$ is a real divisor, arguments as in
~\cite{PP} show that $X_t$ must be a suborbifold, i.e. smooth on a uniformizing neighborhood.
But $X_t$ has a crossing singularity along $P_i$, a contradiction.
Therefore, $D\cap\Sigma=C_1,\ldots,C_i,P_i,C_{i+n+1},\ldots, C_{2n}$, and $D$ is an analytic
subvariety, and a suborbifold.  Since $D=\ol{O}$ it is a toric variety.  Since
$X_t$ is real, $\ol{D}\subset X_t$.  And $D\cup\ol{D}=X_t$ as both are degree two.
\end{proof}

Note that if the isotropy data of $\mathcal{M}$ is normalized to satisfy conditions a.\ and b.\
before (\ref{thm:asd-class}), then we have the identification
\begin{multline}
\rho_1 =(m_1,n_1),\ldots,\rho_{k+2}=(m_{k+2},n_{k+2}),\rho_{k+3}=-(m_1,n_1),\ldots \\
\ldots, \rho_{2k+4}=-(m_{k+2},n_{k+2})=(m_0,n_0).
\end{multline}
Here, as above, we identify $\rho_i$ with a lattice point in $N=\Z\times\Z$.

\subsection{Sasaki-embeddings}

Associated to each compact toric $\ASD$ Einstein orbifold $\mathcal{M}$ with
$\pi_1^{orb}(\mathcal{M})=e$ is the twistor space $\mathcal{Z}$ and a family of
embeddings $X_t\subset\mathcal{Z}$ where $t\in P\setminus\{t_1,t_2,\ldots,t_{k+2}\}$
and $X=X_t$ is the symmetric toric Fano surface canonically associated to $\mathcal{M}$.
We denote the family of embeddings by
\begin{equation}
\iota_t: X\rightarrow\mathcal{Z}.
\end{equation}

Let $M$ be the total space of the $S^1$-Seifert bundle associated to $\mathbf{K}_X$
or $\mathbf{K}_X^{\frac{1}{2}}$, depending on whether $\ind(X)=1$ or $2$.
\begin{thm}\label{thm:subman}
Let $\mathcal{M}$ be a compact toric $\ASD$ Einstein orbifold with
$\pi_1^{orb}(\mathcal{M})=e$.  There exists a Sasakian structure $(\tilde{g},\tilde{\eta},\xi,\tilde{\Phi})$ on
$M$.  So that if $(X,\tilde{h})$ is the K\"{a}hler structure making $\pi:M\rightarrow X$ a
Riemannian submersion, then we have the following diagram, where the horizontal maps are
isometric embeddings and $(\tilde{g},\tilde{\eta},\xi,\tilde{\Phi})$ is the pull-back of the Sasaki structure
$(g,\eta_1,\xi_1,\Phi_1)$ under $\ol{\iota}_t$.
\begin{equation}\label{diag:subman-fund}
\beginpicture
\setcoordinatesystem units <1pt, 1pt> point at 0 30
\put {$M$} at -15 60
\put {$\mathcal{S}$} at 15 60
\put {$X$} at -15 30
\put {$\mathcal{Z}$} at 15 30
\put {$\mathcal{M}$} at 15 0
\put {$\ol{\iota}_t$} at 0 65
\put {$\iota_t$} at 0 35
\arrow <2pt> [.3, 1] from -6 60 to 6 60
\arrow <2pt> [.3, 1] from -6 30 to 6 30
\arrow <2pt> [.3, 1] from -15 51 to -15 39
\arrow <2pt> [.3, 1] from 15 51 to 15 39
\arrow <2pt> [.3, 1] from 15 21 to 15 9
\endpicture
\end{equation}
Furthermore, the image is $\ol{\iota}_t(M) =\{\eta^c(\ol{X}_t) =0\}\subset\mathcal{S}$, where $\ol{X}_t\in\aut(\mathcal{S},g)$
denotes the induced vector field on $\mathcal{S}$, for $t\in P\setminus\{t_1,t_2,\ldots,t_{k+2}\}$.

If the 3-Sasakian space $\mathcal{S}$ is smooth, then so is $M$.  If $M$ is smooth, then
\[M\underset{\text{diff}}{\cong}\# k(S^2\times S^3), \text{  where  } k=2b_2(\mathcal{S})+1. \]
\end{thm}
\begin{proof}
The adjunction formula gives
$\mathbf{K}_X \cong\mathbf{K}_\mathcal{Z}\otimes[X]|_X =\mathbf{K}_\mathcal{Z}\otimes\mathbf{K}_\mathcal{Z}^{-\frac{1}{2}}|_X
=\mathbf{K}_\mathcal{Z}^{\frac{1}{2}}|_X.$
Let $h$ be the K\"{a}hler-Einstein metric on $\mathcal{Z}$ related to the the 3-Sasakian metric
$g$ on $\mathcal{S}$ by Riemannian submersion.  So $\Ric_{\sm h} =8h$.
Recall that $\mathcal{S}$ is the total space of the $S^1$-Seifert bundle associated to
$\mathbf{L}^{-1}$, or $\mathbf{L}^{-\frac{1}{2}}$ iff $w_2(\mathcal{M})=0$.
Also $M$ is the total space of the $S^1$-Seifert bundle associated to either
$\mathbf{K}_X$ or $\mathbf{K}_X^{\frac{1}{2}}$.  By the above adjunction isomorphism
we lift $\iota_t$ to $\ol{\iota}_t$.  Then we have
\[ g=\eta\otimes\eta + \pi^*h,\]
and $\eta=\frac{d}{8}\theta$ with $\theta$ a connection on
$\mathbf{L}^{-1}$ or $\mathbf{L}^{-\frac{1}{2}}$ and where $d=\ind(\mathcal{Z})=2$ or $4$ respectively.
Then it is not difficult to see that by pulling the connection back by
$\iota_t^*\mathbf{L}^{-1}\cong\mathbf{K}_X$ (or $\iota_t^*\mathbf{L}^{-\frac{1}{2}}\cong\mathbf{K}_X^{\frac{1}{2}}$)
we can pull $\eta$ back to $\ol{\eta}$ on $M$.  And define $\tilde{h} =\iota_t^*(h)$.
Then
\[ \tilde{g}=\tilde{\eta}\otimes\tilde{\eta} +\pi^*\tilde{h} \]
is a Sasakian metric on $M$.

If $\mathcal{S}$ is smooth, then locally the orbifold groups of $\mathcal{Z}$ act on
$\mathbf{L}^{-1}$ (or $\mathbf{L}^{-\frac{1}{2}}$)
without non-trivial stabilizers.  Thus this holds for the bundle
$\mathbf{K}_X$ (or $\mathbf{K}_X^{\frac{1}{2}}$) on $X$.

By a theorem in~\cite{HaeSal91} $\pi_1^{orb}(X)=e$ follows from $\pi^{orb}_1(\mathcal{M})=e$.
Given a 4-dimensional orbifold $X$ with an effective 2-torus $T^2$ action, let Let $X_0$ be the open dense
subset of 2-dimensional orbits.  Then $W_0 =X_0/{T^2}$ is a 2-orbifold.  The only other possible
orbits are of dimensions $1$ and $0$, that is, with stabilizers of dimensions
$1$ and $2$, respectively.  Then $W=X/{T^2}$ is a compact connected oriented 2-orbifold with edges and corners with
each edge labeled with a $\Lambda_i$, where $\Lambda_i$ is a rank 1 sublattice of $\Lambda$, the integral lattice such that
$T^2 =\mathfrak{t}/\Lambda$, such that the two sublattices at a corner are linearly independent.
Then we have the exact sequence
\begin{equation}
\pi_2^{orb}(W_0)\rightarrow \Lambda/{\scriptstyle\sum_i \Lambda_i}\rightarrow\pi_1^{orb}(X)\rightarrow
\pi_1^{orb}(W_0)\rightarrow e.
\end{equation}
Since for both $\mathcal{M}$ and $X$, we have that $W$ is a polygon with no orbifold singularities,
$\pi_1^{orb}(X)=\pi_1^{orb}(\mathcal{M})= \Lambda/{\scriptstyle\sum_i \Lambda_i}$.

Suppose $M$ is smooth.  Since $\pi_1^{orb}(X)=e$, the $\pi_i(M)$ must be finite.  This is because
on a Sasaki manifold only always has $H^1(M/\mathscr{F}_\xi) =H^1(M,\R)$, where $H^1(M/\mathscr{F}_\xi)=H^1(X,\R)$
is basic cohomology.  Since $\pi_1^{orb}(X)=e$, the universal cover $\tilde{M}\rightarrow M$ must be a root
of the $S^1$-Seifert bundle over $X$.  By Proposition~\ref{prop:ind-fano} this is a trivial cover, so $\pi_1(M)=e$.

It is a result of H. Oh~\cite{Oh83} that a simply connected 5-manifold with an effective $T^3$ action has $H_2(M,\Z)=\Z^{\ell-3}$,
where $\ell$ is number of edges of $W$.  Since $M$ is spin, the S. Smale classification of 5-manifolds give
the diffeomorphism.

Recall the 1-form $\eta^c=\eta_2 -\sqrt{-1}\eta_3$ of section (\ref{subsec:3-Sasak}) which is $(1,0)$ with
respect to the CR structure $\Phi_1$.
For $t\in\mathfrak{t}$ let $X_t$ denote the killing vector field on
$\mathcal{Z}$ with lift $\ol{X}_t\in\aut(\mathcal{S},g)$. Then
$\theta(X_t)\in H^0(\mathcal{Z},\mathcal{O}(\mathbf{L}))$ which defines a holomorphic function on
$\mathbf{L}^{-1}$.  The $S^1$ subbundles of $\mathbf{L}^{-1}$ is identified with $\mathcal{S}$.
In this way we get $\theta(X_t)=\eta(X_t)$ as holomorphic functions on $C(\mathcal{S})$.
Complexifying gives the same equality for $t\in\mathfrak{t}_{\C}$.
Thus for $t\in P\setminus\{t_1,t_2,\ldots,t_{k+2}\}$, we have
$M_t :=\ol{\iota}_t(M)=\{\eta^c(X_t)=0\}\subset\mathcal{S}$.

Note that here we are setting 2/3 s of the moment map to zero.
\end{proof}

\section{Consequences}\label{sec:conseq}

\subsection{Sasaki-Einstein metrics}\label{subsec:Sasak-Einst}

In this section we present the new infinite families of Sasakian-Einstein 5-manifolds.

\begin{thm}\label{thm:main-smooth}
Let $(\mathcal{S},g)$ be a toric 3-Sasakian 7-manifold with $\pi_1(\mathcal{S})=e$.
Canonically associated to $(\mathcal{S},g)$ are a special symmetric toric Fano surface $X$ and a
toric Sasakian-Einstein 5-manifold $M$ which fit in the commutative diagram (\ref{diag:subman-fund}).
We have $\pi_1^{orb}(X)=e$ and $\pi_1(M)=e$.  And
\[M\underset{\text{diff}}{\cong}\# k(S^2\times S^3), \text{  where  } k=2b_2(\mathcal{S})+1 \]
Furthermore $(\mathcal{S},g)$ can be recovered from either $X$ or $M$ with their torus actions.
\end{thm}
\begin{proof}
The homotopy sequence
\[\cdots\rightarrow\pi_1(G)\rightarrow\pi_1(\mathcal{S})\rightarrow\pi^{orb}_1(\mathcal{M})\rightarrow e,\]
where $G=SO(3)$ or $Sp(1)$, shows that $\pi_1^{orb}(\mathcal{M})=e$.
The surface $X$ is uniquely determined by theorem~\ref{thm:twistor-div}.
It follows from the proof of Theorem~\ref{thm:subman} that $\pi_1^{orb}(X)=e$ and we have the above
diffeomorphism.  An application of Theorem~\ref{thm:K-E} and the remarks at the end of
Section~\ref{subsec:Sasak-st} give the Sasaki-Einstein structure on $M$.
Given $X$ or $M$ with its Sasakian structure
we can recover the orbifold $\mathcal{M}$, which has a unique toric $\ASD$ Einstein metric
by theorem (\ref{thm:asd-class}).  This uniquely determines the 3-Sasakian manifold by results
of section~\ref{subsec:3-Sasak}.
\end{proof}

\begin{thm}\label{thm:main-S-E}
For each odd $k\geq 3$ there is a countably infinite number of toric Sasaki-Einstein structures on
$\# k(S^2\times S^3)$.
\end{thm}
\begin{proof}
Recall from corollary (\ref{cor:3-Sasak-coh}) there are infinitely homotopically distinct smooth
simply connected 3-Sasakian manifolds $\mathcal{S}$ with $b_2(\mathcal{S})=k$ for $k>0$.
From theorem (\ref{thm:main-smooth}) associated to each $\mathcal{S}$ is a distinct
Sasakian-Einstein manifold diffeomorphic to $\# m(S^2\times S^3)$,
where $m=2k+1$.
\end{proof}
The Sasaki-Einstein structures $(g,\eta,\xi,\Phi)$ of Theorem~\ref{thm:main-S-E} have the property of being
isomorphic to the conjugate structure $(g,-\eta,-\xi,-\Phi)$.  This is because the K\"{a}hler-Einstein
orbifold $(X,h)$ has an anti-holomorphic involution $\varsigma:X\rightarrow X$.  Using a real embedding
$X\subset\mathcal{Z}$ in Theorem~\ref{thm:twistor-div} one gets a K\"{a}hler metric with K\"{a}hler
form $\omega$ with $\omega\in 2\pi c^{orb}_1(X)$ and $\varsigma^*\omega =-\omega$.  Then in solving
(\ref{eq:Monge-Amp}) one restricts to functions in $C^\infty (X)^G$ which are $\varsigma$-invariant.

The restriction of $k$ to be odd is merely a limitation on the techniques used.  Subsequent to
these examples appearing in the author's Ph.D. thesis, it was proved~\cite{ChoFutOn08} that there are toric Sasaki-Einstein
structures on $\# k(S^2\times S^3)$ for all $k$.

If a simply connected 5-manifold has two Sasakian-Einstein structures with, non-proportional Reeb vector fields,
for the same metric $g$, then it is $S^5$.
\begin{cor}
For each odd $k\geq 3$ there is a countably infinite number of cohomogeneity 2 Einstein metrics
on $\# k(S^2\times S^3)$.  In particular, the identity component of the isometry group is
$T^3$.
\end{cor}
These metrics have the following curious property.
\begin{prop}
For $M=\# k(S^2\times S^3)$ with $k>1$ odd, let $g_i$ be the sequence of Einstein metrics in the
theorem normalized so that $\Vol_{\sm g_i}(M)=1$.  Then we have $\Ric_{g_i}=\lambda_i g_i$ with the
Einstein constants $\lambda_i\rightarrow 0$ as $i\rightarrow\infty$.
\end{prop}
\begin{proof}
We have
\[\Vol(M,g)=d \left(\frac{\pi}{3}\right)^3 \Vol(\Sigma_{-k}),\]
for the volume of a Sasakian-Einstein manifold with toric leaf space $X$ the anti-canonical polytope
$\Sigma_{-k}$.  This is because an argument in~\cite{Gui94} shows that $\Vol(X_{\Sigma_{-k}})=\Vol(\Sigma_{-k})$.
We have $d=1$ or $2$.  The above Sasakian-Einstein manifolds have leaf spaces
$X_i$, where $X_i =X_{\Delta^*_i}$.  Observe that the polygons $\Delta^*_i$ get arbitrarily large,
and the anti-canonical polytopes $(\Sigma_{-k})_i$ satisfy
\[ \Vol((\Sigma_{-k})_i)\rightarrow 0,\text{  as  }i\rightarrow\infty.\]
\end{proof}

This implies the following.
\begin{thm}
The moduli space of Einstein structures, with a $T^3$ isometry group, on each of the manifolds $\# k(S^2\times S^3)$ for $k\geq 1$ odd
has infinitely many connected components.
\end{thm}
The case $k=1$ is covered by homogeneous examples by M. Wang and W. Ziller~\cite{WaZ}.

There are a couple of consequence of these examples following from some finiteness results.
There is a result of M. Gromov~\cite{Gro} that says that a manifold which admits a metric of nonnegative sectional
curvature satisfies a bound on the total Betti number depending only on the dimension.
Further, he proved that if the diameter is bounded, then as the total Betti number goes to infinity
the infimum of the sectional curvatures goes to $-\infty$.  For any $\kappa\leq 0$ and a fixed diameter $D >0$ there exists
$k_0$ so that, for $k>k_0$, $\# k(S^2\times S^3)$ does not admit a metric with sectional curvature
$K\geq\kappa$ and $\operatorname{diam}\leq D$.  We have the following.
\begin{thm}\label{thm:Gromov}
For any $\kappa\leq 0$ and fixed $\lambda >0$
there are infinitely many simply connected Einstein 5-manifolds with Einstein constant $\lambda$ which do not admit Einstein metrics
with Einstein constant $\lambda$ and sectional $K\geq\kappa$.
\end{thm}

One can also consider these examples in relation to a compactness result of M. Anderson~\cite{And}.
He showed that the space of Riemannian n-manifolds $(M,g)$,
$\mathscr{M}(\lambda, c, D)$ with $\Ric_{\sm g}=\lambda g$,
$\operatorname{inj}(g)\geq c>0$, and $\operatorname{diam}\leq D$ is compact in the $C^\infty$
topology.  For fixed $k>1$ odd in Theorem~\ref{thm:main-S-E} the Sasakian-Einstein metrics $g_i$
on $M=\# k(S^2\times S^3)$ have $\lambda=4$.  We have
$\Vol_{\sm g_i}(M)\rightarrow 0$ as $i\rightarrow\infty$, so no subsequence converges.
We have the following.
\begin{thm}\label{thm:Anderson}
For the the sequence of Einstein manifolds $(M,g_i)$ we have $\inj(g_i)\rightarrow 0$ as $i\rightarrow\infty$.
Also, take any sequence $k_i>1$ of odd integers and examples from Theorem~\ref{thm:main-S-E}
$(\# k_i(S^2\times S^3), g_i)$, then we have $\inj(g_i)\rightarrow 0$ as $i\rightarrow\infty$.
\end{thm}
Examples of Einstein 7-manifolds with properties as in Theorem~\ref{thm:Gromov} and in the second statement
of Theorem~\ref{thm:Anderson} have been given in~\cite{BGMR98}.  These are the toric
3-Sasakian 7-manifolds $\mathcal{S}_\Omega$ considered here.

\subsection{Space of toric $\ASD$ structures}\label{sec:asd-space}

We will compute the dimension of $H^1(\mathcal{Z},\Theta_\mathcal{Z})$ for a twistor space of a toric $\ASD$
Einstein orbifold $(\mathcal{M},g)$.  This will give both
the dimension of the local deformation space of $\mathcal{Z}$ and the local deformation space of $[g]$ as an $\ASD$
conformal class. Recall that the infinitesimal deformations of $[g]$ as an $\ASD$ conformal class correspond to
$\re H^1(\mathcal{Z},\Theta_\mathcal{Z})$.  Although, $\mathcal{Z}$ is not smooth here, both the conformal class $[g]$
and the twistor space structure on $\mathcal{Z}$ are defined on local uniformizing charts of $\mathcal{M}$ and uniformizing
charts these induce on $\mathcal{Z}$.  Thus the twistor correspondence of~\cite{AtiHitSin78} applies here.

In the following $X\subset Z$ will denote a relatively smooth divisor as in Theorem~\ref{thm:twistor-div} with
$\mathbf{K}^{-1}_\mathcal{Z} =[2X]$.  All of the sheaves, divisors and bundles are in the orbifold sense, but the
methods used below carry over to this case.

\begin{prop}\label{prop:cohom}
$h^1(\mathcal{Z},\Theta_\mathcal{Z}(-X))=h^{1,1}(\mathcal{Z})-1$.
\end{prop}
\begin{proof}
By Serre duality
\[ H^1(\mathcal{Z},\Theta_\mathcal{Z}(-X))=H^1(\mathcal{Z},\Theta_\mathcal{Z}\otimes\mathcal{O}(\mathbf{K}_{Z}^{\frac{1}{2}}))
=H^2(\mathcal{Z},\Omega^1\otimes\mathcal{O}(\mathbf{K}_{\mathcal{Z}}^{\frac{1}{2}}))=H^2(\mathcal{Z},\Omega^1(\mathbf{L}^{-1})). \]
From the restriction
\[ 0\rightarrow\Omega_\mathcal{Z}^1(-X)\rightarrow\Omega_\mathcal{Z}^1 \rightarrow\Omega^1_\mathcal{Z}|_X \rightarrow 0, \]
we obtain the exact sequence
\begin{equation}\label{eq:exact1}
0\rightarrow H^{1,1}(\mathcal{Z})\rightarrow H^1(X,\Omega^1_\mathcal{Z})\rightarrow H^2(\mathcal{Z},\Omega_\mathcal{Z}^1(-X))\rightarrow 0,
\end{equation}
because $H^1(\mathcal{Z},\Omega^1(-X))=0$ by Kodaira-Nakano vanishing and $H^{1,2}(\mathcal{Z})=0$.
It is proved in~\cite{Sal82,BoyGal97} that all the cohomology of $\mathcal{Z}$ vanishes besides $H^{k,k}(\mathcal{Z})$.

Next, consider the conormal sequence
\begin{equation}\label{eq:conorm}
0\rightarrow \mathcal{O}_\mathcal{Z}(-X)|_X\rightarrow\Omega^1_\mathcal{Z}|_X \rightarrow\Omega^1_X \rightarrow 0.
\end{equation}
From Theorem~\ref{thm:twistor-div} two distinct $X_{s_1},\ X_{s_2}\subset\mathcal{Z}$, with
$s_1,s_2\in P\setminus\{t_1,t_2,\ldots,t_{n}\}$, have
$X_{s_1} \cap X_{s_2} =\cup_{j=1}^{2k} C_k$, the anti-canonical divisor of $X_{s_1}$ and $X_{s_2}$.
Then
\begin{equation}\label{eq:vanish1}
H^1(X,\mathcal{O}_X (-X)) =H^1(X,\mathcal{O}(\mathbf{K}_X))=0,\text{ and }H^2(X,\mathcal{O}_X(-X))=H^2(X,\Omega_X^2)=\C
\end{equation}
And note that
\begin{equation}\label{eq:vanish2}
H^2(X,\Omega_\mathcal{Z}^1|_X) =H^0(X,\Theta_\mathcal{Z} \otimes\mathcal{O}_X(\mathbf{K}_X))=0.
\end{equation}
This can be seen as follows.  Consider
\begin{equation}\label{eq:theta-res}
0\rightarrow\Theta_\mathcal{Z}(-2X)\rightarrow\Theta_\mathcal{Z}(-X)\rightarrow\Theta_\mathcal{Z} \otimes\mathcal{O}_X(\mathbf{K}_X)\rightarrow 0.
\end{equation}
Note that $H^0(\mathcal{Z},\Theta_\mathcal{Z}(-X))=0$.  If $\beta\in H^0(\mathcal{Z},\Theta_\mathcal{Z}(-X))$,
then restricting $\beta$ to the normal bundle of a generic twistor line $P_x$ gives a section of
$N\cong\mathcal{O}(1)\oplus\mathcal{O}(1)$ vanishing to order 2, which therefore must vanish.
Thus $\beta$ must be tangent to the twistor lines.  But by the definition of the complex structure on $\mathcal{Z}$,
that is impossible.  And by Serre duality,
$H^1(\mathcal{Z},\Theta_\mathcal{Z}(-2X))=H^2(\mathcal{Z},\Omega_\mathcal{Z}^1)=0$.  Then (\ref{eq:vanish2}) follows from the cohomology
sequence of (\ref{eq:theta-res}).

The long exact sequence of (\ref{eq:conorm}) gives
\begin{equation}\label{eq:exact2}
0\rightarrow H^1(X,\Omega^1_\mathcal{Z})\rightarrow H^1(X,\Omega_X^1)\rightarrow\C\rightarrow 0.
\end{equation}
Since $h^{1,1}(X)=2h^{1,1}(\mathcal{Z})$ the proposition follows from (\ref{eq:exact1}) and (\ref{eq:exact2}).
\end{proof}

The arguments of~\cite{YanNag59} applied here prove the following.
\begin{lem}\label{lem:conf-Kill}
Let $(\mathcal{M},g)$ be a compact Einstein orbifold.  If $g$ admits a conformal-Killing vector field which
is not Killing, then $(\mathcal{M},g)$ is isometric to $S^n/\Gamma$, where $S^n$ has the constant curvature
metric and $\Gamma$ is a linear group of isometries fixing $(\pm 1,0,\ldots,0)\in S^n$.
\end{lem}

\begin{lem}
If $X$ is a toric Fano orbifold surface, then $H^1(X,\Theta_X)=0$.
\end{lem}
\begin{proof}
Let $C=\sum_i C_i$ be the anti-canonical divisor.
For each $\sigma\in\Delta^*$, let $\sigma'$ be the corresponding cone in the sublattice $N'$ as in
Proposition~\ref{prop:toric-orb}.

Suppose, for the moment, that $X$ be a nonsingular toric variety, in particular $X=U_{\sigma^\prime}$.
If $\Omega^1_X$ denotes
the algebraic sheaf of differential forms and $\Omega^1_{X}(\log C)$ the sheaf of differential
forms with logarithmic poles along $C=\sum_i C_i$, then
\begin{equation}\label{eq:toric-dif}
0\rightarrow\Omega^1_X \longrightarrow\Omega^1_{X}(\log C)\longrightarrow\bigoplus_{i=1}^{d}\mathcal{O}_{C_i}
\rightarrow 0,
\end{equation}
and
\begin{equation}\label{eq:dif-ident}
\Omega^1_{X}(\log C)\cong\mathcal{O}_X\oplus\mathcal{O}_X.
\end{equation}
See~\cite{Ful93} for a proof.  Consider (\ref{eq:toric-dif}) on each uniformizing neighborhood $U_{\sigma^\prime}$
of $X$.  It is easy to see that (\ref{eq:toric-dif}) and the identification (\ref{eq:dif-ident}) are compatible
with the identifications of the $U_{\sigma^\prime}$.  Extending to the structure sheaf of analytic functions
we have
\begin{equation}\label{eq:dif-orb}
0\rightarrow\Omega^1_X\longrightarrow\mathcal{O}_X\oplus\mathcal{O}_X\longrightarrow
\bigoplus_{i=1}^{d}\mathcal{O}_{C_i}\rightarrow 0,
\end{equation}
where $\Omega^1_C$ denotes the coherent analytic sheaf associated to the orbifold bundle of 1-forms.

We have $H^1(X,\Theta_X)\cong H^1(X,\Omega^1(\mathbf{K}_X))$.  Tensor (\ref{eq:dif-orb}) with
$\mathcal{O}(\mathbf{K}_X)$ and take the long exact cohomology sequence
\[\cdots\rightarrow\bigoplus_{i}H^0(C_i,\mathcal{O}_{C_i}(\mathbf{K}_X))\rightarrow H^1(X,\Omega^1(\mathbf{K}_X))
\rightarrow H^1(X,\mathcal{O}(\mathbf{K}_X))^{\oplus 2}\rightarrow\cdots. \]
Since $\mathbf{K}_{X}<0$, Kodaira-Nakano vanishing shows that
\[ H^0(C_i,\mathcal{O}_{C_i}(\mathbf{K}_X))=
H^1(X,\mathcal{O}(\mathbf{K}_X))=0,\]
thus $H^1(X,\Omega^1(\mathbf{K}_X))=0$.
\end{proof}

Let $\mathbf{N}_{X/\mathcal{Z}}$ be the normal orbifold bundle to $X$ in $\mathcal{Z}$ and $\mathcal{N}_{X/\mathcal{Z}}$ its sheaf
of sections.  Note that
$\mathcal{N}_{X/\mathcal{Z}}=\mathcal{O}_X(X)=\mathcal{O}_X(\mathbf{L})$.  Let $\Theta_{\mathcal{Z},X}$
be the sheaf of sections of $T\mathcal{Z}$ tangent to
$X$.  We will make use of the following exact sequences
\begin{gather}
0\rightarrow\Theta_{\mathcal{Z},X}\longrightarrow\Theta_\mathcal{Z}\longrightarrow\mathcal{N}_{X/\mathcal{Z}}\rightarrow 0\label{sheaf1} \\
0\rightarrow\Theta_\mathcal{Z}(-X)\longrightarrow\Theta_{\mathcal{Z},X}\longrightarrow\Theta_X\rightarrow 0.\label{sheaf2}
\end{gather}
Since $H^1(X,\mathcal{O}_X(\mathbf{L}))=H^1(X,\Omega_X^2(\mathbf{K}_X^{-2}))=0$ by Kodaira-Nakano vanishing,
(\ref{sheaf1}) gives
\begin{equation}\label{eq:coh-ident1}
h^1(\mathcal{Z},\Theta_\mathcal{Z})=h^1(\mathcal{Z},\Theta_{\mathcal{Z},X})-h^0(X,\mathcal{O}_X(X))+h^0(\mathcal{Z},\Theta_\mathcal{Z})
-h^0(\mathcal{Z},\Theta_{\mathcal{Z},X}).
\end{equation}
Since $H^0(\mathcal{Z},\Theta_{\mathcal{Z}}(-X))=0$ and $H^1(X,\Theta_X)=0$, from (\ref{sheaf2}) we have
\begin{equation}\label{eq:coh-ident2}
0\rightarrow H^0(\mathcal{Z},\Theta_{\mathcal{Z},X})\rightarrow H^0(X,\Theta_X)\rightarrow
H^1(\mathcal{Z},\Theta_\mathcal{Z}(-X))\rightarrow H^1(\mathcal{Z},\Theta_{\mathcal{Z},X})\rightarrow 0.
\end{equation}
Combining Proposition~\ref{prop:cohom}, (\ref{eq:coh-ident1}) and (\ref{eq:coh-ident2}) we get
\begin{equation}\label{eq:coh-ident3}
h^1(\mathcal{Z},\Theta_\mathcal{Z})=h^{1,1}(\mathcal{Z})-1 +h^0(\mathcal{Z},\Theta_\mathcal{Z})+h^0(X,\Theta_X)-h^0(X,\mathcal{O}_X(X)).
\end{equation}

If $\mathcal{M}$ is not conformally flat, i.e. $W_g^- \not\equiv 0$, then by Lemma~\ref{lem:conf-Kill} every conformal-Killing
vector field of $(\mathcal{M},g)$ is Killing and $H^0(\mathcal{Z},\Theta_\mathcal{Z})\cong H^0(\mathcal{Z},\mathcal{O}(\mathbf{L}))$.
Recall that $\re H^0(\mathcal{Z},\Theta_\mathcal{Z})$ is isomorphic to the space of conformal-Killing vector fields of
$(\mathcal{M},g)$, and $\re H^0(\mathcal{Z},\mathcal{O}(\mathbf{L}))$ the space of Killing vector fields.

One easily checks that
$h^0(\mathcal{Z},\mathcal{O}(X))=h^0(X,\mathcal{O}_X(X)) +1$, and from (\ref{eq:coh-ident3}) we have
\begin{equation}\label{eq:def-coh1}
h^1(\mathcal{Z},\Theta_\mathcal{Z}) =h^{1,1}(\mathcal{Z})-h^0(X,\Theta_X).
\end{equation}

Note that the compact torus $T^2$ acts on the above cohomology, and denote by $H^1(\mathcal{Z},\Theta_\mathcal{Z})^{T^2}$, etc.,
the fixed set.  One can see that, as $X$ is a toric variety, $H^0(X,\Theta_X)^{T^2} =\mathfrak{t}_{\C}^2$ as follows.
Suppose $\beta\in H^0(X,\Theta_X)^{T^2}$.  Then $\beta$ is also invariant under $\C^*\times\C^*$.
Let $\rho(t)$ be the one parameter group of transformations
generated by $\beta$.  Let $x\in U:=X\setminus\cup_{i=1}^{2k}C_i$, the open orbit isomorphic to $\C^*\times\C^*$.
Fix $t_0\in\C$ close to zero, and let $g\in\C^*\times\C^*$ be the element such that $g\rho(t_0)x=x$.
Since $g\rho(t_0)$ commutes with $\C^*\times\C^*$, it fixes all of $U$ and thus $X$.  Thus
$\rho(t)\in\C^*\times\C^*$.

Since $T^2$ acts on all the sheafs and preserves the exact sequences above,
\begin{equation}\label{eq:def-coh2}
h^1(\mathcal{Z},\Theta_\mathcal{Z})^{T^2}=h^{1,1}(\mathcal{Z})-h^0(X,\Theta_X)^{T^2}=h^{1,1}(\mathcal{Z})-2.
\end{equation}
Therefore from (\ref{eq:def-coh1}) and (\ref{eq:def-coh2}) we must have
\begin{equation}
h^1(\mathcal{Z},\Theta_\mathcal{Z})=h^1(\mathcal{Z},\Theta_\mathcal{Z})^{T^2}=h^{1,1}(\mathcal{Z})-2=b_2(\mathcal{M})-1.
\end{equation}
Note that this also proves that $h^0(X,\Theta_X)=2$ when $W_g^- \not\equiv 0$, which is always the case when $b_2(\mathcal{M})\geq 1$.

We have
\begin{equation}
H^3(\mathcal{Z},\Theta_\mathcal{Z})=H^0(\mathcal{Z},\Omega^1(\mathbf{K}_\mathcal{Z}))=0,
\end{equation}
where the second equality holds for any twistor space.  Also,
\begin{equation}
H^2(\mathcal{Z},\Theta_\mathcal{Z})=H^1(\mathcal{Z},\Omega^1(\mathbf{K}_\mathcal{Z}))=0,
\end{equation}
by Kodaira-Nakano vanishing.

We have proved the following.
\begin{prop}
Let $\mathcal{Z}$ be the twistor space of a compact toric $\ASD$ Einstein orbifold $(\mathcal{M},g)$.  Then
$h^2(\mathcal{Z},\Theta_\mathcal{Z})=h^3(\mathcal{Z},\Theta_\mathcal{Z})=0$.  If $b_2(\mathcal{M})\geq 1$, then
\[ h^1(\mathcal{Z},\Theta_\mathcal{Z})=h^1(\mathcal{Z},\Theta_\mathcal{Z})^{T^2} =b_2(\mathcal{M})-1.\]
If $b_2(\mathcal{M})=0$, then $\mathcal{M}$ has an orbifold covering by $S^4$, which has the round metric.
\end{prop}

We show that the local space of $\ASD$ metrics coincides with the $T^2$-invariant $\ASD$ conformal metrics given by
D. Joyce~\cite{Joy95}.  An explicit description of the toric $\ASD$ Einstein metrics $(\mathcal{M},g)$ was given
by D. Calderbank and M. Singer~\cite{CalSin06}, which made use of the description in~\cite{CalPed02} of
toric $\ASD$ Einstein metrics in terms of an eigen function potential on the hyperbolic plane.  The conformal
classes of these metrics are always given by the Joyce equation.

\begin{thm}
Let $(\mathcal{M},g)$ be a compact toric $\ASD$ Einstein orbifold, then locally the space of $\ASD$ conformal
classes near $[g]$ are those given by the Joyce ansatz.  It is therefore a space of dimension $b_2(\mathcal{M})-1$.
\end{thm}
\begin{proof}
We first prove a lemma.

\begin{lem}\label{lem:orb-isom}
Suppose $(\mathcal{M},g)$ is a compact toric $\ASD$ Einstein orbifold and \mbox{$W_g^- \not\equiv 0$}.  If $(\mathcal{M},g)$ is not
homogeneous, then the connected component
of the identity of the isometry group $\Isom(g)_0$ is either $T^2$ or $U(2)$ up to finite cover.
In the second case $b_2(\mathcal{M})=1$.
If $(\mathcal{M},g)$ is homogeneous, then $(\mathcal{M},g)$ is isometric to $\overline{\cps}^2$ with the Fubini-Study metric.
\end{lem}
\begin{proof}
We may assume that $b_2>0$, for otherwise the structure of toric 4-orbifolds~\cite{OrlRay70,HaeSal91}) implies that $\mathcal{M}$
is diffeomorphic to $S^4/\Gamma$ where $\Gamma$ is a finite group acting as in Proposition~\ref{lem:conf-Kill}.

Let $G=\Isom(g)_{0}$ be the connected component of the identity.  Let $\alpha\in\mathcal{H}^2_{g\,-}\cong H^2(M,\R)$
be non-zero.  Note that standard Bochner techniques show that $\mathcal{H}^2_{g\, +}=0$.
Let $x\in\mathcal{M}$ be a point in the open dense subset of principle orbits where $\alpha_x\neq 0$.
Since $G$ fixes $\alpha$, we have $H_x\subseteq U(2)$, where $H_x$ is the isotropy subgroup at $x$.  We have
$\dim G\leq 3+\dim H_x\leq 7$.  Note that $G$ cannot contain a 3-torus, because that 3-torus would give
a cohomogeneity one action.  And by the theory of such actions $\mathcal{M}$ would not be simply connected.
Considering possible compact groups of rank 2 of these dimensions we see that $G$ is $T^2, \U(2)$ or $\SO(4)$
up to finite coverings.  If $G=\SO(4)$ up to coverings, then $\dim H_x\geq 3$ and is disjoint from $T^2$ which is impossible.
Suppose $G$ is $\U(2)$.  Then the generic orbit $Gx$ is 3 dimensional and $H_{x}=S^1$.  For if $\dim Gx=2$,
then $Gx\cong T^2$ and $H_x\cong T^2$ which is impossible, since $G/T^2\cong S^2$ for any $T^2\subset G$.

So $\mathcal{M}$ is of cohomogeneity 1.  The orbifold isotropy group of a point $\Gamma_x$ is preserved by $G$.
Standard arguments show that the set of smooth
points of $\mathcal{M}$ is diffeomorphic to $(0,1)\times G/H$, where $H_0 =S^1$.  Since otherwise the orbit space
$\mathcal{M}/G \cong S^1$, would contradict $\pi_1(\mathcal{M})=e$.
Adding the two, possibly singular orbits at 0 and 1 gives a close dense subset of $\mathcal{M}$, so it must be all of $\mathcal{M}$.
Thus the orbit space is $\pi:\mathcal{M}\rightarrow\mathcal{M}/G=[0,1]$.  The each of the orbits $\pi^{-1}(0)$ and
$\pi^{-1}(1)$ is either a point or $\cps^1$.  Since $b_2>0$, at least one is $\cps^1$.  It is easy to see
that both orbits cannot be $\cps^1$, because in this case the isotropy subgroups of $T^2$ on
these two $\cps^1$ must be equal, which is impossible by the results in\cite{CalSin06}.
Therefore $b_2(\mathcal{M})=1$.

If $\mathcal{M}$ is homogeneous then it must be smooth.  Since $b_2>0$, the last statement of the proposition follows
from~\cite{Hit81}.
\end{proof}

We review the D. Joyce construction~\cite{Joy95} of $\ASD$ conformal metrics with a surface orthogonal action of $T^2$
by conformal transformations.  By \emph{surface orthogonal} we mean that the orthogonal distribution to the $T^2$ orbits
is integrable.  Locally all toric surface orthogonal $\ASD$ metrics are of this form.  The metrics are defined by
linearly independent solutions to a linear equation on the spinor bundle $\mathcal{W}\rightarrow\mathcal{H}^2$
over the hyperbolic plane.  See also~\cite{CalPed02,CalSin04} for more details on the following.

Let $\mathbb{V}$ be a real 2-dimensional vector space with a symplectic form $\varepsilon(\cdot,\cdot)$.  We
consider bundle isomorphisms $\Phi:\mathcal{W}\rightarrow\mathcal{H}\times\mathbb{V}$, and define a $\mathbb{V}$-invariant
metric on $\mathcal{H}\times\mathbb{V}$ in terms of $(\mathcal{H},h)$ and $\Phi$.  We define a family of metrics on $\mathbb{V}$ by
\[ (u,v)_{\Phi} =h(\Phi^{-1}(u),\Phi^{-1}(v)), \]
and on $\mathcal{H}\times\mathbb{V}$ by
\begin{equation}\label{eq:Joyce-asd}
 g_\Phi =\Omega^2(h +(\cdot,\cdot)_\Phi).
\end{equation}

Fix the spinor bundle $\mathcal{W}$ by $\mathcal{W}\otimes_\C \mathcal{W} =TN$.  We consider the half-space model
of $\mathcal{H}$ with coordinates $(\eta,\rho),\ \rho>0$, with metric $h=(d\rho^2 + d\eta^2)/\rho^2$.
Given a smooth section $\Phi\in C^\infty(\mathcal{W}^*)$ we have the \emph{Joyce equation}
\begin{equation}\label{eq:Joyce-eq}
\ol{\partial}\Phi =\frac{1}{2}\ol{\Phi}.
\end{equation}
We clarify the identifications made in (\ref{eq:Joyce-eq}).  On the left-hand side
\[\ol{\partial}\Phi\in\Gamma(\overline{T^*\mathcal{H}}\otimes_{\C}\mathcal{W}^*),\]
while using the induce Hermitian metrics on the spinor bundles
\[  \overline{T^*\mathcal{H}}\otimes_{\C}\mathcal{W}^* =\ol{\mathcal{W}^*}\otimes_{\C}\ol{\mathcal{W}^*}\otimes_{\C}\mathcal{W}^* =\ol{\mathcal{W}^*}.\]
Then $g_\Phi$ is an $\ASD$ metric if $\Phi\in C^\infty(\mathcal{H},\mathcal{W}^*\otimes\mathbb{V})$ is a linearly independent
solution to (\ref{eq:Joyce-eq}).

Considering $\mathcal{W}^*$ as a real bundle we can identify $S^2_0(\mathcal{W}^*)$ with $T^*\mathcal{H}$, in which there
is an orthonormal frame $\lambda_0,\lambda_1$ of $\mathcal{W}^*$ and identifications
$\lambda_0^2 -\lambda_1^2 =d\rho/\rho$ and $2\lambda_0 \lambda_1 =d\eta/\rho$.
Then a solution $\Phi\in C^\infty(\mathcal{H},\mathcal{W}^*\otimes\mathbb{V})$ can be written
\[ \Phi=\lambda_0 \otimes v_0 +\lambda_1 \otimes v_1 ,\]
with $v_0,v_1 \in C^\infty(\mathbb{V})$ satisfying the equations
\begin{equation}\label{eq:Joyce-expl}
\rho\partial_\rho v_0 +\rho\partial_\eta v_1 =v_0,\quad \rho\partial_\eta v_0 -\rho\partial_\rho v_1 =0.
\end{equation}
Then if $\mu_0, \mu_1$ is a dual frame to $\lambda_0,\lambda_1$,
\[ \Phi^{-1} =\frac{\varepsilon(v_0,\cdot)\otimes\mu_1 -\varepsilon(v_1,\cdot)\oplus\mu_0}{\varepsilon(v_0,v_1)}.\]
D. Joyce made the observation that $-\lambda_1$ is obviously a solution to (\ref{eq:Joyce-expl}) and acting on it by
$\SL(2,\R)$ gives the family of \emph{fundamental solutions} to (\ref{eq:Joyce-eq})
\begin{equation}
\phi(\rho,\eta,x) =\frac{\rho\lambda_0 +(\eta-x)\lambda_1}{\sqrt{\rho^2 +(\eta -x)^2}},
\end{equation}
where $x\in\partial\mathcal{H}$.

Conditions were given in~\cite{CalPed02,CalSin06} for the $\ASD$ structure $g_\Phi$ to be conformally Einstein.
The condition is that the solution $\Phi$ to (\ref{eq:Joyce-eq}) comes from an eigenfunction of $\Delta_{h}$, the negative
Laplacian on $(\mathcal{H},h)$,
\[ \Delta_{h} F =\frac{3}{4}F,  \quad F\in C^\infty(\mathcal{H}). \]
We define $f(\rho,\eta)=\sqrt{\rho}F(\rho,\eta)$, then with
\[ v_0 =(f_\rho, \eta f_\rho -\rho f_\eta),\quad v_1 =(f_\eta, \rho f_\rho +\eta f_\eta -f) \]
$\Phi =\lambda_0 \otimes v_0 +\lambda_1 \otimes v_1$ is a solution to (\ref{eq:Joyce-eq}) and
\begin{equation}\label{eq:ASD-Einst}
g_F =\frac{|F^2-4|dF|^2|}{4F^2}\Bigl(\frac{d\rho^2 + d\eta^2}{\rho^2} +\frac{\varepsilon(v_0,\cdot)^2 +\varepsilon(v_1,\cdot)^2}{\varepsilon(v_0,v_1)^2}  \Bigr)
\end{equation}
is an $\ASD$ Einstein equation with positive scalar curvature where $F^2 >4|dF|^2$ and negative scalar curvature where
$4|dF|^2 >F^2 >0$.

By Theorem~\ref{thm:asd-class} the isotropy data $(m_1,n_1),(m_2,n_2),\ldots,(m_{k+2},n_{k+2})$, with $(m_0,n_0)=(m_{k+2},n_{k+2})$,
of the $\ASD$ Einstein space $\mathcal{M}$ can be arranged as follows, after possibly changing signs and acting by $\SL(2,\Z)$.
If we define $(a_i,b_i),\ i=0,\ldots, k+2$, by
\[ 2(a_i,b_i) =(m_i,n_i)-(m_{i-1},n_{i-1}), \]
then $a_i >0$ and $x_i :=b_i /a_i$ are increasing for $i=0,\ldots, k+2$, where we set $x_0=-\infty$.  It was proved
in~\cite{CalSin06} that the potential
\[ F= \sum_{i=1}^{k+2} \frac{\sqrt{a_i \rho^2 +(a_i \eta-b_i)^2}}{\sqrt{\rho}},\]
gives uniquely the $\ASD$ Einstein metric on $\mathcal{M}$ by (\ref{eq:ASD-Einst}).  Then an easy computation gives
\begin{equation}\label{eq:asd-sol}
\Phi =\frac{1}{2}\sum_{i=1}^{k+2} \bigl(\phi(\rho,\eta,x_i)-\phi(\rho,\eta,x_{i-1}) \bigr)\otimes(m_i,n_i).
\end{equation}
Here the $x_1 < x_2<\cdots <x_{k+2}$ are points on $\partial\mathcal{H}$ corresponding to the points on the boundary
of $\mathcal{Q}$ fixed by $T^2$.  Conversely, given a sequence of points $x_1 < x_2<\cdots <x_{k+2}$ on $\partial\mathcal{H}$
the arguments in~\cite{Joy95} show that (\ref{eq:asd-sol}) gives a solutions which compactifies on the toric orbifold
$\mathcal{M}$ with the given isotropy data to a complete metric.  In particular the same arguments there show that
$\bigwedge^2 \Phi$ is non-vanishing on $\mathcal{H}$.  Clearly, $\PSL(2,\R)$ acts on (\ref{eq:asd-sol}) giving
isomorphic solutions and acts on the formula (\ref{eq:asd-sol}) by shifting the $x_1 < x_2<\cdots <x_{k+2}$.
Thus we have a $k-1$-parameter space of $\ASD$ structures.  See Figure~\ref{fig:hyp}.

Suppose two metrics in this family are conformally isomorphic.  So we have $\Phi_1$ defined by $x_1 < x_2<\cdots <x_{k+2}$
and $\Phi_2$ defined by $z_1 < z_2<\cdots <z_{k+2}$, metrics $g_1$ and $g_2$ defined in (\ref{eq:Joyce-asd}), and a
diffeomorphism $\psi:\mathcal{M}\rightarrow\mathcal{M}$ so that $\psi^* g_2 =e^2f g_1$.  By changing the conformal factor
$e^2f$ we may assume that $\psi$ is an isomorphism.   The arguments in Lemma~\ref{lem:orb-isom} show that
$\Isom(g_1)_0 =\Isom(g_2) =T^2$ unless $b_2(\mathcal{M})=1$, so we may assume that $\Isom(g_1)_0 =\Isom(g_2) =T^2$.
So $\psi$ must map the Lie algebra $\mathfrak{t}$ of $T^2$ to itself, and $\psi$ is equivariant up to an automorphism
$\psi_{T} \in\GL(2,\Z)$.  Since $\psi$ preserves the vector fields generated by $T^2$ it must preserve the orthogonal
distributions to the torus orbits.  Therefore $\phi$ descends to a conformal automorphism $\phi_{\mathcal{H}}$
of $\mathcal{H}$.  Thus $\phi_{\mathcal{H}}\in\PSL^*(2,\R)$, the isometry group of $(\mathcal{H},h)$ generated by
$\PSL(2,\R)$ and the orientation reversing $\eta+\sqrt{-1}\rho\mapsto -\eta+\sqrt{-1}\rho$.
Also, $\psi$ must map the fixed points $p_1,p_2,\ldots p_{k+2}$ of $T^2$ to themselves.
If $\psi$ permutes the fixed points, then it must permute the edges of $\mathcal{Q}$ in such a way that
the edge with stabilizer data $(m_i,n_i)$ goes to the edge with stabilizer data $\pm (\psi_{T})_* (m_i,n_i)$.
Thus $\psi_{T}$ is an automorphism of data of the toric orbifold.  This is obviously a finite group.

Suppose that $\phi$ fixes the points $p_1,p_2,\ldots p_{k+2}$.  Then $\phi_{\mathcal{H}} \in\PSL(2,\R)$ and
$(\psi_{T})_* (m_i,n_i) =\pm (m_i,n_i)$.  It is easy to see that $\psi_{T} =\pm\mathbb{1}$.
And $\phi_{\mathcal{H}}$ maps $x_1 < x_2<\cdots <x_{k+2}$ to $z_1 < z_2<\cdots <z_{k+2}$.
Therefore the map
\[ \Bigl\{ \{x_1 <x_2 <\ldots <x_{k+2}\}\subset\rps^1 \Bigr\}/\PSL(2,\R) \rightarrow\Bigl\{\text{Conf. }\ASD
\text{ str. on }\mathcal{M}\Bigr\} \]
is finite to one.

\begin{figure}[t!]
\labellist
\hair 2pt
\pinlabel $\mathcal{H}$ at 225 76
\pinlabel $(m_0,n_0)$ [l] at 95 25
\pinlabel $(m_1,n_1)$ at 240 25
\pinlabel {$(m_2,n_2)$} [r] at 380 25
\endlabellist
\centering
\includegraphics[scale=.6]{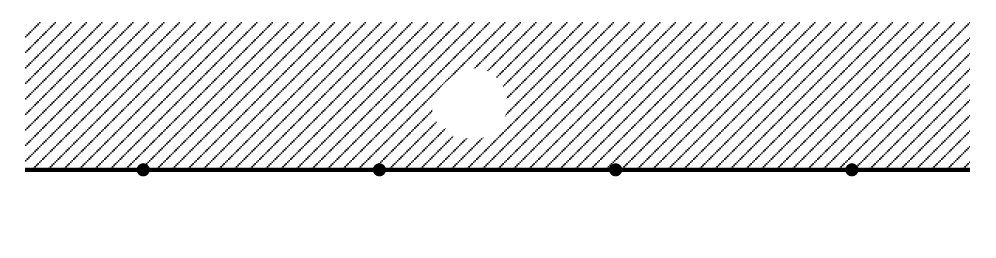}
\caption{orbit space of $\mathcal{M}$}
\label{fig:hyp}
\end{figure}
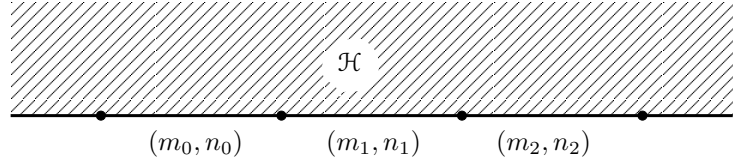

Let $g_0$ denote the $\ASD$ Einstein metric.  Recall that $\ASD$ conformal classes $[g]$ on $\mathcal{M}$, with its fixed
orbifold structure, are in correspondence with its twistor space $\mathcal{Z}$, with properties (a),(b) and (c)
at the beginning of Section~\ref{subsec:twist-div}.  Let $\mathcal{Z}_{[g]}$ denotes the twistor space of
$(\mathcal{M},[g])$ for $[g],\ \ASD$.  Then by the orbifold Riemann-Roch theorem
of~\cite{Kaw79}
\begin{equation}\label{eq:orb-R-R}
\chi(\Theta_{\mathcal{Z}_{[g]}})=h^0(\Theta_{\mathcal{Z}_{[g]}})-h^1(\Theta_{\mathcal{Z}_{[g]}})+h^2(\Theta_{\mathcal{Z}_{[g]}}),
\end{equation}
is independent of $[g]$.  By Lemma~\ref{lem:orb-isom} we may assume that $\Isom(g_0)_0 =T^2$.
Let $\mathcal{U}$ be a neighborhood of $[g_0]$ in the $C^\infty$ topology with $h^2(\Theta_{\mathcal{Z}_{[g]}})=0$.
Let $\mathcal{J}$ be the space of toric $\ASD$ structures constructed above.  Then for $[g]\in\mathcal{J}\cap\mathcal{U}$
we have $h^0(\Theta_{\mathcal{Z}_{[g]}})=2$ by the above assumption, and by (\ref{eq:orb-R-R}) $h^1(\Theta_{\mathcal{Z}_{[g]}})
=b_2(\mathcal{M})-1$.  Therefore $H^1(\Theta_{\mathcal{Z}_{[g]}}) =H^1(\Theta_{\mathcal{Z}_{[g]}})^{T^2}$ and
the twistor spaces $\mathcal{Z}_{[g]}$ for $[g]\in\mathcal{J}\cap\mathcal{U}$ provide a real subspace of
the deformation space of $\mathcal{Z}_{[g]}$ for $[g]\in\mathcal{J}\cap\mathcal{U}$.
The total local deformation space of $\mathcal{Z}_{[g]}$ is a complex thickening of the real deformations.
\end{proof}

\begin{rmk}
This is in contrast to the case of $\ASD$ structures on $\#\ell\overline{\cps^2}$.  There are many examples of toric $\ASD$ structures
on $\#\ell\overline{\cps^2},\ \ell\geq 3,$ for which most deformations are not toric. It is a result of A. Fujiki~\cite{Fuj00}
that the toric $\ASD$ conformal metrics on $\#\ell\overline{\cps^2}$ are the Joyce metrics, so each is in a $\ell-1$ dimensional
family.  But $\chi(\Theta_{Z_{[g]}})=\frac{1}{2}(15\chi +29\tau)=15-7\ell$ by a calculation originally due to N. Hitchin and I. Singer.
See A. King and D. Kotschick~\cite{KinKot92} for more details on the moduli of $\ASD$ conformal metrics.
\end{rmk}

\section{Examples}

We consider some of the examples obtained starting with the simplest.  In particular
we can determine some of the spaces in diagram (\ref{diag:intro}) associated to a smooth
toric 3-Sasakian 7-manifold more explicitly in some cases.

\subsection{Smooth examples}

It is well known that there exists only two complete examples of positive scalar curvature
anti-self-dual Einstein manifolds~\cite{Hit81,FriKur82}, $S^4$ and $\overline{\cps}^2$ with the round and
Fubini-Study metrics respectively, where $\overline{\cps}^2$ denotes $\cps^2$ with the opposite of the usual orientation.
\\

\noindent
\emph{$\mathcal{M}= S^4$}

Considering the spaces in diagram (\ref{diag:intro}) we have:
$\mathcal{M}= S^4$ with the round metric; its twistor space $\mathcal{Z}=\cps^3$ with the Fubini-study metric;
the quadratic divisor $X\subset\mathcal{Z}$ is $\cps^1\times\cps^1$ with the
homogeneous K\"{a}hler-Einstein metric; $M = S^2\times S^3$ with the homogeneous Sasakian-Einstein structure;
and $\mathcal{S}= S^7$ has the round metric.
In this case diagram (\ref{diag:intro}) becomes the following.

\begin{equation}
\beginpicture
\setcoordinatesystem units <1pt, 1pt> point at 0 30
\put {$S^2\times S^3$} at -30 60
\put {$S^7$} at 30 60
\put {$\cps^1\times\cps^1$} at -30 30
\put {$\cps^3$} at 30 30
\put {$S^4$} at 30 0
\arrow <2pt> [.3, 1] from -1 60 to 11 60
\arrow <2pt> [.3, 1] from -1 30 to 11 30
\arrow <2pt> [.3, 1] from -30 51 to -30 39
\arrow <2pt> [.3, 1] from 30 51 to 30 39
\arrow <2pt> [.3, 1] from 30 21 to 30 9
\endpicture
\end{equation}

This is the only example, I am aware of, for which the horizontal maps are isometric immersions
when the toric surface and Sasakian space are equipped with the Einstein metrics.
\\

\noindent
\emph{$\mathcal{M}=\cps^2$}

In this case $\mathcal{M}=\cps^2$ with the Fubini-Study metric; its twistor
space is $\mathcal{Z}=F_{1,2}$, the manifold of flags $V\subset W\subset\C^3$ with $\dim V =1$ and $\dim W=2$,
with the homogeneous K\"{a}hler-Einstein metric.
The projection $\pi:F_{1,2}\rightarrow\cps^2$ is as follows.  If $(p,l)\in F_{1,2}$ so
$l$ is a line in $\cps^2$ and $p\in l$, then $\pi(p,l)=p^\perp \cap l$, where $p^\perp$ is the
orthogonal compliment with respect to the standard hermitian inner product.
We can define $F_{1,2}\subset\cps^2\times(\cps^2)^*$ by
\[ F_{1,2}=\{([p_0:p_1:p_2],[q^0:q^1:q^2])\in\cps^2\times(\cps^2)^*: \sum p_iq^i =0\}.\]
And the complex contact structure is given by $\theta=q^idp_i -p_idq^i$.
Fix the action of $T^2$ on $\cps^2$ by
\[ (e^{i\theta},e^{i\phi})[z_0:z_1:z_2]=[z_0:e^{i\theta}z_1:e^{i\phi}z_2].\]
Then this induces the action on $F_{1,2}$
\[ (e^{i\theta},e^{i\phi})([p_0:p_1:p_2],[q^0:q^1:q^2])=([p_0:e^{i\theta}p_1:e^{i\phi}p_2],[q^0:e^{-i\theta}q^1:e^{-i\phi}q^2]). \]
Given $[a,b]\in\cps^1$ the one parameter group $(e^{ia\tau},e^{ib\tau})$ induces the holomorphic vector field
$W_\tau\in\Gamma(T^{1,0}F_{1,2})$ and the quadratic divisor $X_\tau=(\theta(W_\tau))$ given by
\[ X_\tau =(ap_1q^1+bp_2q^2=0,\quad p_iq^i=0).\]
One can check directly that $X_\tau$ is smooth for $\tau\in\cps^1\setminus\{[1,0],[0,1],[1,1]\}$ and
$X_\tau =\cps^2_{(3)}$, the equivariant blow-up of $\cps^2$ at 3 points.  For
$\tau\in\{[1,0],[0,1],[1,1]\}$, $X_\tau =D_\tau +\ol{D}_\tau$ where
both $D_\tau,\ol{D}_\tau$ are isomorphic to the Hirzebruch surface
$F_1 =\mathbb{P}(\mathcal{O}_{\cps^1}\oplus\mathcal{O}_{\cps^1}(1))$

The Sasakian-Einstein space is $M=\#3(S^2\times S^3)$.  And we have
$\mathcal{S}=\mathcal{S}(1,1,1)=SU(3)/U(1)$ with the homogeneous 3-Sasakian structure.
This case has the following diagram.

\begin{equation}
\beginpicture
\setcoordinatesystem units <1pt, 1pt> point at 0 30
\put {$\#3(S^2\times S^3)$} at -35 60
\put {$SU(3)/U(1)$} at 35 60
\put {$\cps^2_{(3)}$} at -35 30
\put {$F_{1,2}$} at 35 30
\put {$\cps^2$} at 35 0
\arrow <2pt> [.3, 1] from -6 60 to 6 60
\arrow <2pt> [.3, 1] from -6 30 to 6 30
\arrow <2pt> [.3, 1] from -35 51 to -35 39
\arrow <2pt> [.3, 1] from 35 51 to 35 39
\arrow <2pt> [.3, 1] from 35 21 to 35 9
\endpicture
\end{equation}

\subsection{Galicki-Lawson quotients}

The simplest examples of quaternionic-K\"{a}hler quotients are the Galicki-Lawson examples first appearing in~\cite{GalLaw88}
and further considered in~\cite{BGMR98}.
These are circle quotients of $\qps^2$.  In this case the weight matrices are of the form
$\Omega=\mathbf{p}=(p_1,p_2,p_3)$ with the admissible set
\[\{\mathcal{A}_{1,3}(Z)=\{\mathbf{p}\in\Z^3| p_i\neq 0 \text{ for }i=1,2,3 \text{ and } \gcd(p_i,p_j)=1\text{ for }i\neq j\} \]
We may take $p_i>0$ for $i=1,2,3$.
The zero locus of the 3-Sasakian moment map $N(\mathbf{p})\subset S^{11}$ is diffeomorphic to the
Stiefel manifold $V_{2,3}^\C$ of complex 2-frames in $\C^3$ which can be identified as
$V_{2,3}^\C\cong U(3)/U(1)\cong SU(3)$.
Let $f_{\mathbf{p}}: U(1)\rightarrow U(3)$ be
\[ f_{\mathbf{p}}(\tau)=\begin{bmatrix}\tau^{p_1} & 0 & 0 \\ 0 & \tau^{p_2} & 0 \\ 0 & 0 & \tau^{p_3}\end{bmatrix}. \]
Then the 3-Sasakian space $\mathcal{S}(\mathbf{p})$ is diffeomorphic to the quotient of $SU(3)$ by the
action of $U(1)$
\[\tau\cdot W=f_{\mathbf{p}}(\tau)Wf_{(0,0,-p_1 -p_2 -p_3)}(\tau)\text{  where  }\tau\in U(1)\text{ and }W\in SU(3). \]
Thus $\mathcal{S}(\mathbf{p})\cong SU(3)/U(1)$ is a biquotient similar to the examples considered by
Eschenburg in~\cite{Es}.

The action of the group $SU(2)$ generated by $\{\xi^1,\xi^2,\xi^3\}$ on
$N(\mathbf{p})\cong SU(3)$ commutes with the action of $U(1)$.
We have $N(\mathbf{p})/SU(2)\cong SU(3)/SU(2)\cong S^5$ with $U(1)$ acting by
\[\tau\cdot v=f_{(-p_2-p_3,-p_1-p_3,-p_1-p_2)}v \text{  for  }v\in S^5\subset\C^3. \]
We see that $\mathcal{M}_{\Omega}\cong\cps^2_{a_1,a_2,a_3}$ where $a_1=p_2+p_3, a_2=p_1+p_3, a_3=p_1+p_2$
and the quotient metric is anti-self-dual with the reverse of usual orientation.
If $p_1,p_2,p_3$ are all odd then the generic leaf of the 3-Sasakian foliation $\mathcal{F}_3$ is
$SO(3)$.  If exactly one is even, then the generic leaf is $Sp(1)$.
Denote by $X_{p_1,p_2,p_3}$ the toric Fano divisor, which can be considered as a generalization of
$\cps^2_{(3)}$.
We have the following spaces and embeddings.

\begin{equation}
\beginpicture
\setcoordinatesystem units <1pt, 1pt> point at 0 30
\put {$\#3(S^2\times S^3)$} at -38 60
\put {$\mathcal{S}(p_1,p_2,p_3)$} at 35 60
\put {$X_{p_1,p_2,p_3}$} at -30 30
\put {$\mathcal{Z}(p_1,p_2,p_3)$} at 35 30
\put {$\cps^2_{a_1,a_2,a_3}$} at 45 0
\arrow <2pt> [.3, 1] from -8 60 to 4 60
\arrow <2pt> [.3, 1] from -10 30 to 2 30
\arrow <2pt> [.3, 1] from -30 51 to -30 39
\arrow <2pt> [.3, 1] from 35 51 to 35 39
\arrow <2pt> [.3, 1] from 35 21 to 35 9
\endpicture
\end{equation}

A simple series of examples can be obtained by taking $\mathbf{p}=(2q-1,1,1)$ for any $q\geq 1$.
Then the anti-self-dual Einstein space is $\mathcal{M}=\cps^2_{1,q,q}$ which is homeomorphic to $\cps^2$,
but its metric is ramified along a $\cps^1$ to order $q$.  For the toric divisor $X\subset\mathcal{Z}$ we have
$X=\cps^2_{(3)}$ with the metric ramified along two $\cps^1$'s to order $q$.  We get an sequence of
distinct Sasakian-Einstein structures on $M\cong\# 3(S^2\times S^3)$.

\begin{figure}
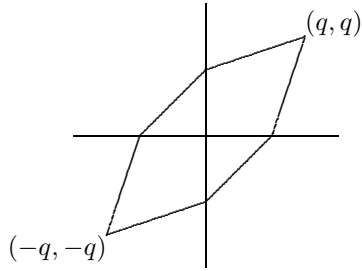

\centering
\mbox{
\beginpicture
\setcoordinatesystem units <.5pt,.5pt> point at 0 0
\setplotarea x from -100 to 100, y from -100 to 100
\axis bottom shiftedto x=0 /
\axis left shiftedto y=0 /
\setlinear
\plot 0 50  75 75  50 0  0 -50  -75 -75  -50 0  0 50 /
\put {$(q,q)$} [lb] at 75 75
\put {$(-q,-q)$} [rt] at -75 -75
\endpicture}
\caption{infinite Fano orbifold structures on $\cps^2_{(3)}$}
\end{figure}

\bibliographystyle{amsplain}

\providecommand{\bysame}{\leavevmode\hbox to3em{\hrulefill}\thinspace}
\providecommand{\MR}{\relax\ifhmode\unskip\space\fi MR }
\providecommand{\MRhref}[2]{%
  \href{http://www.ams.org/mathscinet-getitem?mr=#1}{#2}
}
\providecommand{\href}[2]{#2}

\end{document}